\documentclass{amsart}
\usepackage{amsmath}
\usepackage{amssymb}
\usepackage{mathrsfs}
\usepackage{pgf,tikz}
\usetikzlibrary{shadings}
\usetikzlibrary{arrows,automata}
\usepackage{graphicx}
\usepackage{systeme,mathtools}
\usepackage{amsthm,enumitem}
\usepackage{systeme}
\input xy
\xyoption{all}

\newtheorem{theorem}{Theorem}[section]

\newtheorem{conjecture}{Conjecture}[section]
\newtheorem{lemma}[theorem]{Lemma}
\newtheorem{corollary}[theorem]{Corollary}
\newtheorem{proposition}[theorem]{Proposition}
\theoremstyle{definition}
\newtheorem{definition}[theorem]{Definition}
\newtheorem{remark}{Remark}
\newtheorem{example}{Example}

\begin{document}

\title[Multiplicative groups of Leavitt path algebras]{Multiplicative groups of Leavitt path algebras}

\maketitle

\begin{center}	Bui Xuan Hai\footnote{Corresponding author}\end{center} 

\begin{center}
\tiny{	(1) \textit{University of Science, Ho Chi Minh City, Vietnam}\\ (2)\textit{Vietnam National University, Ho Chi Minh City, Vietnam}\\ \textit{e-mail: bxhai@hcmus.edu.vn}}
\end{center}

\begin{center}
	and
\end{center}

\begin{center}
	Huynh Viet Khanh
\end{center}

\begin{center}
\tiny{	 \textit{Department of Mathematics and Informatics, HCMC University of Education,\\ 280 An Duong Vuong Str., Dist. 5, Ho Chi Minh City, Vietnam}\\ \textit{e-mail: khanhhv@hcmue.edu.vn}}
\end{center}
\begin{abstract}
In this paper, we prove that the multiplicative group of a unital non-commutative Leavitt path algebra $L_K(E)$ and Cohn path algebra $C_K(E)$ contain a non-cyclic free subgroup, provided $K$ is a non-absolute field. We also provide a description of the generators of  free subgroups in term of the graph $E$. Finally, we determine  multiplicative groups  of Leavitt path algebras of some special types.
\end{abstract}
\section{Introduction and preliminaries}
Leavitt path algebras were firstly introduced and substantially studied in \cite{Pa_abrams-pino_2005}. Although various aspects of Leavitt path algebras were intensively studied by many mathematicians, the multiplicative groups of unital Leavitt path algebras have seemed to have gone relatively unexplored. The main purpose of this paper is to make the efforts to touch this study. More precisely, we mainly focus our attention on the existence of  non-cyclic free subgroups in the multiplicative groups of  non-commutative unital Leavitt path algebras. It will turn out that a non-commutative unital Leavitt path algebra always contains a non-cyclic free subgroup provided  the base field is non-absolute. We also give a description of the generators of such a free subgroup in term of the base graph. Roughly speaking, the generators may be determined by a sink, an infinite emitter or infinite paths.

Throughout this paper, the basic notations and conventions we use are standard and they are taken from the book of Gene Arbrams, Pere Ara and Mercedes Siles Molina \cite{Bo_abrams_2017}. In particular, a directed graph $E=(E^0, E^1, r, s)$ consists of two sets $E^0$ and $E^1$ together with maps $r, s: E^1\to E^0$. The elements of $E^0$ are called \textit{vertices} and the elements of $E^1$ \textit{edges}. For the convenience, in the following, we outline some concepts and notations that will be used substantially in this paper.

A vertex $v$ is a \textit{sink} if it emits no edges, while it is an \textit{infinite emitter }if it emits infinitely many edges. A vertex $v$ is said to be \textit{regular} if it is neither a sink nor an infinite emitter. The set of all regular vertices in a graph $E$ is denoted by ${\rm Reg}(E)$. For each $e\in E^1$, a \textit{ghost edge} $e^*$ is defined by setting $r(e^*):=s(e)$ and $s(e^*):=r(e)$. The set of all ghost edges is denoted by $(E^1)^*$. A \textit{finite path} $\mu$ of length $\ell(\mu):=n\ge 1$ is a finite sequence of edges $\mu=e_1e_2\dots e_n$ with $r(e_i) =s(e_{i+1})$ for all $1\leq i\leq n-1$. We set $s(\mu):=s(e_1)$ and $r(\mu):=r(e_n)$. The set of all finite paths in $E$ is denoted by ${\rm Path}(E)$. Let $\mu=e_1\dots e_n\in {\rm Path}(E)$. If $v=s(\mu)=r(\mu)$, then we say that $\mu$ is a \textit{closed path based at} $v$. If moreover, $s(e_j)\ne v$ for every $j>1$ then $\mu$ is a \textit{simple closed path based at} $v$. If $\mu$ is  a closed path based at $v$ and $s(e_i)\ne s(e_j)$ for every $i\ne j$,
then $\mu$ is called a \textit{cycle based at} $v$. Assume that $\mu=e_1\dots e_n$ is a cycle based at the vertex $v$. Then for each $1\le i\le n$, the path $e_ie_{i+1}\dots e_ne_1\dots e_{i-1}$ is a cycle based at $s(e_i)$. We call the collection of cycles $\{\mu_i\}$ based at $s(e_i)$ the \textit{cycle of} $\mu$. A \textit{cycle} $c$ is a set of paths consisting of the cycle of $\mu$ for  $\mu$ some cycle based at a vertex $v$. The \textit{length of a cycle} $c$ is the length of any path in $c$. In particular, a cycle of length $1$ is called a \textit{loop}. 

Also, $\mu^*=e_n^*\dots e_2^*e_1^*$ is the corresponding ghost path of the path $\mu=e_1\dots e_n$. The set of all vertices on a path $\mu$ is denoted by $\mu^0$. 

An \textit{exit} for a path $\mu=e_1\dots e_n$ is an edge $e$ such that $s(e)=s(e_i)$ for some $i$ and $e\ne e_i$. A graph $E$ is said to satisfy the \textit{Condition (L)} if every cycle in $E$ has an exit.

For $u,v\in E^0$, if there exists a path from $u$ to $v$, then we write $u\geq v$. A subset $S$ of $E^0$ is said to have the \textit{Countable Separation Property} with respect to a set $C$, if $C$ is a countable  set of $E^0$ with the property that for each $u\in S$ there is a $v\in C$ such that $u\geq v$.

For each $v\in E^0$ and each infinite path $p$, we define
$$ M(v)=\{w\in E^0\;|\; w\geq v\} \text{ and } M(p)=\{w\in E^0\;|\; w\geq v \text{ for some } v\in p^0\}.$$

Let $H$ be a subset of $E^0$. We say that $H$ is \textit{hereditary}  if whenever $u\in H$ and $u\geq v$ for some vertex $v$, then $v\in H$; and $H$ is \textit{saturated} if, for any regular vertex $v$, $r(s^{-1}(v))\subseteq H$ implies $v\in H$. 

The next two definitions are taken from \cite{Pa_tom_2007}.
\begin{definition} Assume that $H$ is a hereditary saturated subset of $E^0$. A vertex $w$ is called a \textit{breaking vertex} of $H$ if $w\in E^0\backslash H$ is an infinite emitter such that 
	$$
	1\leq |s^{-1}(w)\cap r^{-1}(E^0\backslash H)|<\infty.
	$$
	The set of all breaking vertices of $H$ is denoted by $B_H$. For each $w\in B_H$, we set 
	$$
	w^H=w-\sum_{\substack{s(e)=w,\\r(e)\not\in H}}ee^*.
	$$ 
\end{definition}
\begin{example}\label{example_1}
	Consider the graph:
	\begin{center}
	\tikzset{every picture/.style={line width=0.5pt}} 
		
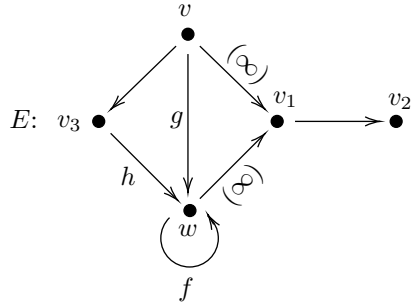
		\begin{tikzpicture}[x=0.75pt,y=0.75pt,yscale=-1,xscale=1]
		
		\draw  [color={rgb, 255:red, 0; green, 0; blue, 0 }  ,draw opacity=1 ][fill={rgb, 255:red, 0; green, 0; blue, 0 }  ,fill opacity=1 ] (133,135) .. controls (133,133.34) and (134.34,132) .. (136,132) .. controls (137.66,132) and (139,133.34) .. (139,135) .. controls (139,136.66) and (137.66,138) .. (136,138) .. controls (134.34,138) and (133,136.66) .. (133,135) -- cycle ;
		
		\draw  [draw opacity=0] (144.92,137.71) .. controls (148.21,140.36) and (150.31,144.43) .. (150.31,148.99) .. controls (150.31,156.99) and (143.82,163.48) .. (135.82,163.48) .. controls (127.81,163.48) and (121.32,156.99) .. (121.32,148.99) .. controls (121.32,144.9) and (123.01,141.22) .. (125.72,138.58) -- (135.82,148.99) -- cycle ; \draw   (144.92,137.71) .. controls (148.21,140.36) and (150.31,144.43) .. (150.31,148.99) .. controls (150.31,156.99) and (143.82,163.48) .. (135.82,163.48) .. controls (127.81,163.48) and (121.32,156.99) .. (121.32,148.99) .. controls (121.32,144.9) and (123.01,141.22) .. (125.72,138.58) ;
		
		\draw    (148,142) -- (146.03,138.71) ;
		\draw [shift={(145,137)}, rotate = 419.03999999999996] [color={rgb, 255:red, 0; green, 0; blue, 0 }  ][line width=0.5]    (7.65,-2.3) .. controls (4.86,-0.97) and (2.31,-0.21) .. (0,0) .. controls (2.31,0.21) and (4.86,0.98) .. (7.65,2.3)   ;
		
		\draw  [color={rgb, 255:red, 0; green, 0; blue, 0 }  ,draw opacity=1 ][fill={rgb, 255:red, 0; green, 0; blue, 0 }  ,fill opacity=1 ] (132,46) .. controls (132,44.34) and (133.34,43) .. (135,43) .. controls (136.66,43) and (138,44.34) .. (138,46) .. controls (138,47.66) and (136.66,49) .. (135,49) .. controls (133.34,49) and (132,47.66) .. (132,46) -- cycle ;
		
		\draw  [color={rgb, 255:red, 0; green, 0; blue, 0 }  ,draw opacity=1 ][fill={rgb, 255:red, 0; green, 0; blue, 0 }  ,fill opacity=1 ] (177,90) .. controls (177,88.34) and (178.34,87) .. (180,87) .. controls (181.66,87) and (183,88.34) .. (183,90) .. controls (183,91.66) and (181.66,93) .. (180,93) .. controls (178.34,93) and (177,91.66) .. (177,90) -- cycle ;
		
		\draw  [color={rgb, 255:red, 0; green, 0; blue, 0 }  ,draw opacity=1 ][fill={rgb, 255:red, 0; green, 0; blue, 0 }  ,fill opacity=1 ] (87,90) .. controls (87,88.34) and (88.34,87) .. (90,87) .. controls (91.66,87) and (93,88.34) .. (93,90) .. controls (93,91.66) and (91.66,93) .. (90,93) .. controls (88.34,93) and (87,91.66) .. (87,90) -- cycle ;
		
		\draw    (129,52) -- (97.44,82.61) ;
		\draw [shift={(96,84)}, rotate = 315.88] [color={rgb, 255:red, 0; green, 0; blue, 0 }  ][line width=0.5]    (7.65,-2.3) .. controls (4.86,-0.97) and (2.31,-0.21) .. (0,0) .. controls (2.31,0.21) and (4.86,0.98) .. (7.65,2.3)   ;
		
		\draw    (141,52) -- (172.56,82.61) ;
		\draw [shift={(174,84)}, rotate = 224.12] [color={rgb, 255:red, 0; green, 0; blue, 0 }  ][line width=0.5]    (7.65,-2.3) .. controls (4.86,-0.97) and (2.31,-0.21) .. (0,0) .. controls (2.31,0.21) and (4.86,0.98) .. (7.65,2.3)   ;
		
		\draw    (96,96) -- (128.59,128.59) ;
		\draw [shift={(130,130)}, rotate = 225] [color={rgb, 255:red, 0; green, 0; blue, 0 }  ][line width=0.5]    (7.65,-2.3) .. controls (4.86,-0.97) and (2.31,-0.21) .. (0,0) .. controls (2.31,0.21) and (4.86,0.98) .. (7.65,2.3)   ;
		
		\draw    (141,129) -- (172.59,97.41) ;
		\draw [shift={(174,96)}, rotate = 495] [color={rgb, 255:red, 0; green, 0; blue, 0 }  ][line width=0.5]    (7.65,-2.3) .. controls (4.86,-0.97) and (2.31,-0.21) .. (0,0) .. controls (2.31,0.21) and (4.86,0.98) .. (7.65,2.3)   ;
		
		\draw    (135,56) -- (135,124) ;
		\draw [shift={(135,126)}, rotate = 270] [color={rgb, 255:red, 0; green, 0; blue, 0 }  ][line width=0.5]    (7.65,-2.3) .. controls (4.86,-0.97) and (2.31,-0.21) .. (0,0) .. controls (2.31,0.21) and (4.86,0.98) .. (7.65,2.3)   ;
		
		\draw    (189,90) -- (230,90) ;
		\draw [shift={(232,90)}, rotate = 180] [color={rgb, 255:red, 0; green, 0; blue, 0 }  ][line width=0.5]    (7.65,-2.3) .. controls (4.86,-0.97) and (2.31,-0.21) .. (0,0) .. controls (2.31,0.21) and (4.86,0.98) .. (7.65,2.3)   ;
		
		\draw  [color={rgb, 255:red, 0; green, 0; blue, 0 }  ,draw opacity=1 ][fill={rgb, 255:red, 0; green, 0; blue, 0 }  ,fill opacity=1 ] (237,90) .. controls (237,88.34) and (238.34,87) .. (240,87) .. controls (241.66,87) and (243,88.34) .. (243,90) .. controls (243,91.66) and (241.66,93) .. (240,93) .. controls (238.34,93) and (237,91.66) .. (237,90) -- cycle ;
		
		\draw (-27,165) node [anchor=north west][inner sep=0.5pt]   [align=left] {$ $};
		
		\draw (44,83) node [anchor=north west][inner sep=0.5pt]   [align=left] {$\displaystyle E$:};
		
		\draw (129,141) node [anchor=north west][inner sep=0.5pt]   [align=left] {$\displaystyle w$};
		
		\draw (129,168) node [anchor=north west][inner sep=0.5pt]   [align=left] {$\displaystyle f$};
		
		\draw (125,84) node [anchor=north west][inner sep=0.5pt]   [align=left] {$\displaystyle g$};
		
		\draw (100,112) node [anchor=north west][inner sep=0.5pt]   [align=left] {$\displaystyle h$};
		
		\draw (161.46,46.74) node [anchor=north west][inner sep=0.5pt]  [rotate=-45] [align=left] {$\displaystyle ( \infty $)};
		
		\draw (177.26,118.46) node [anchor=north west][inner sep=0.5pt]  [rotate=-135] [align=left] {$\displaystyle ( \infty $)};
		
		\draw (68,86) node [anchor=north west][inner sep=0.5pt]   [align=left] {$\displaystyle v_{3}$};
		
		\draw (177,74) node [anchor=north west][inner sep=0.5pt]   [align=left] {$\displaystyle v_{1}$};
		
		\draw (235,75) node [anchor=north west][inner sep=0.5pt]   [align=left] {$\displaystyle v_{2}$};
		
		\draw (129,30) node [anchor=north west][inner sep=0.5pt]   [align=left] {$\displaystyle v$};
		
		\end{tikzpicture}
\end{center}
	 	 
	We see that  $H=\{v_1,v_2\}$ is the hereditary saturated of $E^0$. Also, the set of all breaking vertices of $H$ is $B_H=\{v, w\}$.
\end{example}
\begin{definition}\label{definition_1.2}
	For a given hereditary saturated subset $H$ of $E^0$  and a subset $S\subseteq B_H$, we call $(H,S)$ an \textit{admissible pair} and denote by $I(H, S)$ the ideal of $L_K(E)$ generated by the set  $H\cup \{v^H: v\in S\}$. For each admissible pair $(H,S)$, the corresponding \textit{quotient graph} $E\backslash (H,S)$ is defined as follows:
	$$ 
	\begin{aligned}
		& (E\backslash (H,S))^0= (E^0\backslash H) \cup \{v': v\in B_H\backslash S\},\\
		& (E\backslash (H,S))^1=\{e\in E^1: r(e)\not\in H\}\cup \{e': e\in E^1, r(e)\in B_H\backslash S\},
	\end{aligned}
	$$
	and $r$ and $s$ are extended to $(E\backslash (H,S))^0$ by setting $s(e')=s(e)$ and $r(e')=r(e)'$. 
\end{definition}
\begin{example}\label{example_2}
	Consider the graph $E$ in Example \ref{example_1}, and put $S=\{v\}\subseteq B_H$. Then, $B_H\backslash S=\{w\}$ and the edges ending in $w$ are $f$, $g$, $h$. Therefore, the quotient graph $E\backslash (H,S)$ corresponding with the admissible pair $(H,S)$ can be described in detail as follows:
	$$ 
	\begin{aligned}
		& (E\backslash (H,S))^0= (E^0\backslash H) \cup \{w'\},\\
		& (E\backslash (H,S))^1=\{e\in E^1: r(e)\not\in H\}\cup \{f', g', h'\}, 
	\end{aligned}
	$$
	and $r$ and $s$ are extended to $(E\backslash (H,S))^0$ by setting $s(f')=s(f)=v_3$,  $s(g')=s(g)=v$, $s(h)=s(h')=v_3$, and $r(f')=r(g')=r(h')=w'$. 
		 
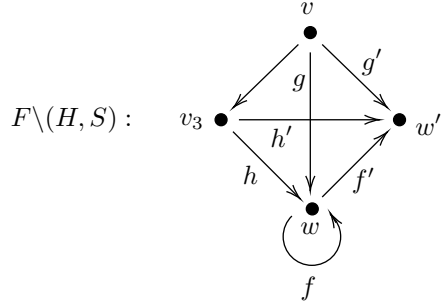
	 \begin{center}
	 	\tikzset{every picture/.style={line width=0.5pt}}       
	 	
	 	\begin{tikzpicture}[x=0.75pt,y=0.75pt,yscale=-1,xscale=1]
	 	
	 	\draw  [color={rgb, 255:red, 0; green, 0; blue, 0 }  ,draw opacity=1 ][fill={rgb, 255:red, 0; green, 0; blue, 0 }  ,fill opacity=1 ] (153,135) .. controls (153,133.34) and (154.34,132) .. (156,132) .. controls (157.66,132) and (159,133.34) .. (159,135) .. controls (159,136.66) and (157.66,138) .. (156,138) .. controls (154.34,138) and (153,136.66) .. (153,135) -- cycle ;
	 	
	 	\draw  [draw opacity=0] (164.92,137.71) .. controls (168.21,140.36) and (170.31,144.43) .. (170.31,148.99) .. controls (170.31,156.99) and (163.82,163.48) .. (155.82,163.48) .. controls (147.81,163.48) and (141.32,156.99) .. (141.32,148.99) .. controls (141.32,144.9) and (143.01,141.22) .. (145.72,138.58) -- (155.82,148.99) -- cycle ; \draw   (164.92,137.71) .. controls (168.21,140.36) and (170.31,144.43) .. (170.31,148.99) .. controls (170.31,156.99) and (163.82,163.48) .. (155.82,163.48) .. controls (147.81,163.48) and (141.32,156.99) .. (141.32,148.99) .. controls (141.32,144.9) and (143.01,141.22) .. (145.72,138.58) ;
	 	
	 	\draw    (168,142) -- (166.03,138.71) ;
	 	\draw [shift={(165,137)}, rotate = 419.03999999999996] [color={rgb, 255:red, 0; green, 0; blue, 0 }  ][line width=0.5]    (7.65,-2.3) .. controls (4.86,-0.97) and (2.31,-0.21) .. (0,0) .. controls (2.31,0.21) and (4.86,0.98) .. (7.65,2.3)   ;
	 	
	 	\draw  [color={rgb, 255:red, 0; green, 0; blue, 0 }  ,draw opacity=1 ][fill={rgb, 255:red, 0; green, 0; blue, 0 }  ,fill opacity=1 ] (152,46) .. controls (152,44.34) and (153.34,43) .. (155,43) .. controls (156.66,43) and (158,44.34) .. (158,46) .. controls (158,47.66) and (156.66,49) .. (155,49) .. controls (153.34,49) and (152,47.66) .. (152,46) -- cycle ;
	 	
	 	\draw  [color={rgb, 255:red, 0; green, 0; blue, 0 }  ,draw opacity=1 ][fill={rgb, 255:red, 0; green, 0; blue, 0 }  ,fill opacity=1 ] (197,90) .. controls (197,88.34) and (198.34,87) .. (200,87) .. controls (201.66,87) and (203,88.34) .. (203,90) .. controls (203,91.66) and (201.66,93) .. (200,93) .. controls (198.34,93) and (197,91.66) .. (197,90) -- cycle ;
	 	
	 	\draw  [color={rgb, 255:red, 0; green, 0; blue, 0 }  ,draw opacity=1 ][fill={rgb, 255:red, 0; green, 0; blue, 0 }  ,fill opacity=1 ] (107,90) .. controls (107,88.34) and (108.34,87) .. (110,87) .. controls (111.66,87) and (113,88.34) .. (113,90) .. controls (113,91.66) and (111.66,93) .. (110,93) .. controls (108.34,93) and (107,91.66) .. (107,90) -- cycle ;
	 	
	 	\draw    (149,52) -- (117.44,82.61) ;
	 	\draw [shift={(116,84)}, rotate = 315.88] [color={rgb, 255:red, 0; green, 0; blue, 0 }  ][line width=0.5]    (7.65,-2.3) .. controls (4.86,-0.97) and (2.31,-0.21) .. (0,0) .. controls (2.31,0.21) and (4.86,0.98) .. (7.65,2.3)   ;
	 	
	 	\draw    (161,52) -- (192.56,82.61) ;
	 	\draw [shift={(194,84)}, rotate = 224.12] [color={rgb, 255:red, 0; green, 0; blue, 0 }  ][line width=0.5]    (7.65,-2.3) .. controls (4.86,-0.97) and (2.31,-0.21) .. (0,0) .. controls (2.31,0.21) and (4.86,0.98) .. (7.65,2.3)   ;
	 	
	 	\draw    (116,96) -- (148.59,128.59) ;
	 	\draw [shift={(150,130)}, rotate = 225] [color={rgb, 255:red, 0; green, 0; blue, 0 }  ][line width=0.5]    (7.65,-2.3) .. controls (4.86,-0.97) and (2.31,-0.21) .. (0,0) .. controls (2.31,0.21) and (4.86,0.98) .. (7.65,2.3)   ;
	 	
	 	\draw    (161,129) -- (192.59,97.41) ;
	 	\draw [shift={(194,96)}, rotate = 495] [color={rgb, 255:red, 0; green, 0; blue, 0 }  ][line width=0.5]    (7.65,-2.3) .. controls (4.86,-0.97) and (2.31,-0.21) .. (0,0) .. controls (2.31,0.21) and (4.86,0.98) .. (7.65,2.3)   ;
	 	
	 	\draw    (155,56) -- (155,124) ;
	 	\draw [shift={(155,126)}, rotate = 270] [color={rgb, 255:red, 0; green, 0; blue, 0 }  ][line width=0.5]    (7.65,-2.3) .. controls (4.86,-0.97) and (2.31,-0.21) .. (0,0) .. controls (2.31,0.21) and (4.86,0.98) .. (7.65,2.3)   ;
	 	
	 	\draw    (119,90) -- (188,90) ;
	 	\draw [shift={(190,90)}, rotate = 180] [color={rgb, 255:red, 0; green, 0; blue, 0 }  ][line width=0.5]    (7.65,-2.3) .. controls (4.86,-0.97) and (2.31,-0.21) .. (0,0) .. controls (2.31,0.21) and (4.86,0.98) .. (7.65,2.3)   ;
	 	
	 	\draw (-7,165) node [anchor=north west][inner sep=0.5pt]   [align=left] {$ $};
	 	
	 	\draw (3,82) node [anchor=north west][inner sep=0.5pt]   [align=left] {$\displaystyle F\backslash ( H,S) :$};
	 	
	 	\draw (149,141) node [anchor=north west][inner sep=0.5pt]   [align=left] {$\displaystyle w$};
	 	
	 	\draw (149,168) node [anchor=north west][inner sep=0.5pt]   [align=left] {$\displaystyle f$};
	 	
	 	\draw (145,66) node [anchor=north west][inner sep=0.5pt]   [align=left] {$\displaystyle g$};
	 	
	 	\draw (120,112) node [anchor=north west][inner sep=0.5pt]   [align=left] {$\displaystyle h$};
	 	
	 	\draw (88,86) node [anchor=north west][inner sep=0.5pt]   [align=left] {$\displaystyle v_{3}$};
	 	
	 	\draw (207,86) node [anchor=north west][inner sep=0.5pt]   [align=left] {$\displaystyle w'$};
	 	
	 	\draw (149,30) node [anchor=north west][inner sep=0.5pt]   [align=left] {$\displaystyle v$};
	 	
	 	\draw (180,53) node [anchor=north west][inner sep=0.5pt]   [align=left] {$\displaystyle g'$};
	 	
	 	\draw (175,113) node [anchor=north west][inner sep=0.5pt]   [align=left] {$\displaystyle f'$};
	 	
	 	\draw (134,92) node [anchor=north west][inner sep=0.5pt]   [align=left] {$\displaystyle h'$};
	 	\end{tikzpicture}
	 \end{center}
	 
\end{example}
\begin{definition}[Leavitt path algebra]\label{definition_1.3}
	Let $E$ be a directed  graph and $K$ a field. We difine a set $\left( E^1\right)^*$ consisting of symbols of the form $\{e^*| e\in E^1\}$.   The \textit{Leavitt path algebra of $E$ with coefficients in $K$}, denoted by $L_K(E)$, is the free associative $K$-algebra generated by $E^0\cup E^1\cup (E^1)^*$, subject to the following relations:
	\begin{enumerate}[]
		\item[(V)] $vv'=\delta_{v,v'}v$ for all $v,v'\in E^0$,
		\item[(E1)] $s(e)e=er(e)=e$ for all $e\in E^1$,
		\item[(E2)] $r(e)e^*=e^*s(e)=e^*$ for all $e\in E^1$,
		\item[(CK1)] $e^*f=\delta_{e,f}r(e)$ for all $e,f\in E^1$, and
		\item[(CK2)] $v=\sum_{\{e\in E^1|s(e)=v\}}ee^*$ for every $v\in{\rm Reg}(E)$.
	\end{enumerate}
\end{definition}
The following definition is taken from \cite[Definition 1.5.1]{Bo_abrams_2017}.
\begin{definition}[Cohn path algebra]
	Let $E$ be an arbitrary graph and $K$ any field. The \textit{Cohn path algebra} of $E$ with coefficients in $K$, denoted by $C_K(E)$, is the free associative $K$-algebra generated by the set $E^0\cup E^1\cup (E^1)^*$, subject to the relations given in (V), (E1), (E2), and (CK1) of Definition \ref{definition_1.3}.
\end{definition}
In general, a Leavitt path algebra $L_K(E)$ does not have the identity. However, it is the case if and only if the vertex set $E^0$ is finite (see \cite[Lemma 1.2.12 (iv)]{Bo_abrams_2017}). 
The following result is deduced from \cite[Theorem 5.7]{Pa_tom_2007} as some its part. For the convenience of use, we restate it here as a lemma.
\begin{lemma}[See {\cite[Theorem 5.7]{Pa_tom_2007}}]\label{lemma_1.4}
	Let $E$ be a graph and $K$ a field.
	\begin{enumerate}[font=\normalfont]
		\item[(1)] For a given  admissible pair $(H,S)$ in the graph $E$, the mapping $$\varphi: L_K(E)\to L_K(E\backslash(H,S))$$ 
		defined by
		\begin{align*}
			&	\varphi(v) =    \begin{cases} 
				v+v'        & \text{ if } v\not\in H \text{ and } v\in B_H\backslash S \\
				v             & \text{ if } v\not\in H \text{ and } v\not\in B_H\backslash S \\
				0             & \text{ if } v\in H
			\end{cases},\\
			&	\varphi(e) =    \begin{cases} 
				e+e'        & \text{ if } r(e)\not\in H \text{ and } r(e)\in B_H\backslash S \\
				e             &  \text{ if } r(e)\not\in H \text{ and } r(e)\not\in B_H\backslash S \\
				0             & \text{ if } r(e)\in H
			\end{cases},  \\                     					
			&		\varphi(e^*) =    \begin{cases} 
				e^*+(e')^*      & \text{if } \text{ if } r(e)\not\in H \text{ and } r(e)\in B_H\backslash S \\
				e^*        			&  \text{ if } r(e)\not\in H \text{ and } r(e)\not\in B_H\backslash S \\
				0        			  & \text{ if } r(e)\in H
			\end{cases}
		\end{align*}
		is an epimorphism with $\ker\varphi = I(H,S)$.
		\item[(2)] The graded ideals of $L_K(E)$ are precisely ideals of the form $I(H,S)$, for some hereditary saturated subset $H$ and some subset $S\subseteq B_H$. Moreover,  $$
		I(H,S)\cap E^0=H \text{ and } \{v\in B_H: v^H\in I\}=S.
		$$
	\end{enumerate}
\end{lemma}

The Leavitt path algebra of a graph provides a rich supply of well-known algebraic structures. For instance, the classical Leavitt K-algebra $L_K(1,n)$ for $n\geq 2$, the full $n\times n$ matrix algebra ${\rm M}_n(K)$, and the Toeplitz $K$-algebra $\mathscr{T}_K$ are, respectively, the Leavitt path algebra of the ``rose with $n$ petals'' graph $R_n$ ($n\geq 2$), the oriented line graph $A_n$ with $n$ vertices, and the Toeplitz graph $E_T$, which are pictured as follows:\\
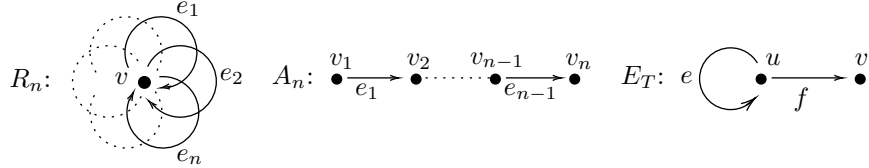
\begin{center}
	\tikzset{every picture/.style={line width=0.5pt}}  
	
	\begin{tikzpicture}[x=0.75pt,y=0.75pt,yscale=-1,xscale=1]
	
	\draw  [color={rgb, 255:red, 0; green, 0; blue, 0 }  ,draw opacity=1 ][fill={rgb, 255:red, 0; green, 0; blue, 0 }  ,fill opacity=1 ] (79.75,42.29) .. controls (79.75,40.67) and (81.06,39.37) .. (82.69,39.37) .. controls (84.31,39.37) and (85.63,40.67) .. (85.63,42.29) .. controls (85.63,43.9) and (84.31,45.21) .. (82.69,45.21) .. controls (81.06,45.21) and (79.75,43.9) .. (79.75,42.29) -- cycle ;
	
	\draw  [draw opacity=0][dash pattern={on 0.84pt off 2.51pt}] (81.05,49.3) .. controls (78.19,55.51) and (71.88,59.83) .. (64.56,59.83) .. controls (54.54,59.83) and (46.43,51.76) .. (46.43,41.8) .. controls (46.43,31.85) and (54.54,23.78) .. (64.56,23.78) .. controls (71.88,23.78) and (78.19,28.09) .. (81.05,34.31) -- (64.56,41.8) -- cycle ; \draw  [dash pattern={on 0.84pt off 2.51pt}] (81.05,49.3) .. controls (78.19,55.51) and (71.88,59.83) .. (64.56,59.83) .. controls (54.54,59.83) and (46.43,51.76) .. (46.43,41.8) .. controls (46.43,31.85) and (54.54,23.78) .. (64.56,23.78) .. controls (71.88,23.78) and (78.19,28.09) .. (81.05,34.31) ;
	
	\draw  [draw opacity=0][dash pattern={on 0.84pt off 2.51pt}] (58.32,37.24) .. controls (58.1,36.9) and (57.88,36.56) .. (57.67,36.2) .. controls (52.67,27.58) and (55.64,16.56) .. (64.31,11.58) .. controls (72.98,6.6) and (84.07,9.56) .. (89.08,18.18) .. controls (92.75,24.5) and (92.13,32.12) .. (88.12,37.69) -- (73.38,27.19) -- cycle ; \draw  [dash pattern={on 0.84pt off 2.51pt}] (58.32,37.24) .. controls (58.1,36.9) and (57.88,36.56) .. (57.67,36.2) .. controls (52.67,27.58) and (55.64,16.56) .. (64.31,11.58) .. controls (72.98,6.6) and (84.07,9.56) .. (89.08,18.18) .. controls (92.75,24.5) and (92.13,32.12) .. (88.12,37.69) ;
	
	\draw  [draw opacity=0] (74.64,34.93) .. controls (72.17,29.76) and (72.22,23.51) .. (75.32,18.18) .. controls (80.32,9.56) and (91.41,6.6) .. (100.08,11.58) .. controls (108.75,16.56) and (111.73,27.58) .. (106.72,36.2) .. controls (103.97,40.94) and (99.38,43.96) .. (94.36,44.9) -- (91.02,27.19) -- cycle ; \draw   (74.64,34.93) .. controls (72.17,29.76) and (72.22,23.51) .. (75.32,18.18) .. controls (80.32,9.56) and (91.41,6.6) .. (100.08,11.58) .. controls (108.75,16.56) and (111.73,27.58) .. (106.72,36.2) .. controls (103.97,40.94) and (99.38,43.96) .. (94.36,44.9) ;
	
	\draw  [draw opacity=0] (84.32,34.31) .. controls (87.19,28.09) and (93.49,23.78) .. (100.82,23.78) .. controls (110.83,23.78) and (118.95,31.85) .. (118.95,41.8) .. controls (118.95,51.76) and (110.83,59.83) .. (100.82,59.83) .. controls (95.09,59.83) and (89.98,57.19) .. (86.66,53.07) -- (100.82,41.8) -- cycle ; \draw   (84.32,34.31) .. controls (87.19,28.09) and (93.49,23.78) .. (100.82,23.78) .. controls (110.83,23.78) and (118.95,31.85) .. (118.95,41.8) .. controls (118.95,51.76) and (110.83,59.83) .. (100.82,59.83) .. controls (95.09,59.83) and (89.98,57.19) .. (86.66,53.07) ;
	
	\draw  [draw opacity=0] (90.22,39.45) .. controls (97.08,38.78) and (104.03,42.06) .. (107.7,48.38) .. controls (112.71,57) and (109.74,68.02) .. (101.06,73) .. controls (92.39,77.98) and (81.3,75.02) .. (76.29,66.4) .. controls (73.4,61.42) and (73.17,55.63) .. (75.17,50.67) -- (92,57.39) -- cycle ; \draw   (90.22,39.45) .. controls (97.08,38.78) and (104.03,42.06) .. (107.7,48.38) .. controls (112.71,57) and (109.74,68.02) .. (101.06,73) .. controls (92.39,77.98) and (81.3,75.02) .. (76.29,66.4) .. controls (73.4,61.42) and (73.17,55.63) .. (75.17,50.67) ;
	
	\draw  [draw opacity=0][dash pattern={on 0.84pt off 2.51pt}] (88.12,46.89) .. controls (92.13,52.46) and (92.75,60.08) .. (89.08,66.4) .. controls (84.07,75.02) and (72.98,77.97) .. (64.31,73) .. controls (55.64,68.02) and (52.67,57) .. (57.68,48.38) .. controls (58.13,47.59) and (58.64,46.85) .. (59.19,46.16) -- (73.38,57.39) -- cycle ; \draw  [dash pattern={on 0.84pt off 2.51pt}] (88.12,46.89) .. controls (92.13,52.46) and (92.75,60.08) .. (89.08,66.4) .. controls (84.07,75.02) and (72.98,77.97) .. (64.31,73) .. controls (55.64,68.02) and (52.67,57) .. (57.68,48.38) .. controls (58.13,47.59) and (58.64,46.85) .. (59.19,46.16) ;
	
	\draw  [color={rgb, 255:red, 0; green, 0; blue, 0 }  ,draw opacity=1 ][fill={rgb, 255:red, 0; green, 0; blue, 0 }  ,fill opacity=1 ] (177.25,40.61) .. controls (177.25,39.26) and (178.35,38.16) .. (179.7,38.16) .. controls (181.06,38.16) and (182.16,39.26) .. (182.16,40.61) .. controls (182.16,41.97) and (181.06,43.07) .. (179.7,43.07) .. controls (178.35,43.07) and (177.25,41.97) .. (177.25,40.61) -- cycle ;
	
	\draw  [color={rgb, 255:red, 0; green, 0; blue, 0 }  ,draw opacity=1 ][fill={rgb, 255:red, 0; green, 0; blue, 0 }  ,fill opacity=1 ] (217.25,40.61) .. controls (217.25,39.26) and (218.35,38.16) .. (219.7,38.16) .. controls (221.06,38.16) and (222.16,39.26) .. (222.16,40.61) .. controls (222.16,41.97) and (221.06,43.07) .. (219.7,43.07) .. controls (218.35,43.07) and (217.25,41.97) .. (217.25,40.61) -- cycle ;
	
	\draw  [color={rgb, 255:red, 0; green, 0; blue, 0 }  ,draw opacity=1 ][fill={rgb, 255:red, 0; green, 0; blue, 0 }  ,fill opacity=1 ] (257.25,40.61) .. controls (257.25,39.26) and (258.35,38.16) .. (259.7,38.16) .. controls (261.06,38.16) and (262.16,39.26) .. (262.16,40.61) .. controls (262.16,41.97) and (261.06,43.07) .. (259.7,43.07) .. controls (258.35,43.07) and (257.25,41.97) .. (257.25,40.61) -- cycle ;
	
	\draw  [color={rgb, 255:red, 0; green, 0; blue, 0 }  ,draw opacity=1 ][fill={rgb, 255:red, 0; green, 0; blue, 0 }  ,fill opacity=1 ] (297.25,40.61) .. controls (297.25,39.26) and (298.35,38.16) .. (299.7,38.16) .. controls (301.06,38.16) and (302.16,39.26) .. (302.16,40.61) .. controls (302.16,41.97) and (301.06,43.07) .. (299.7,43.07) .. controls (298.35,43.07) and (297.25,41.97) .. (297.25,40.61) -- cycle ;
	
	\draw  [color={rgb, 255:red, 0; green, 0; blue, 0 }  ,draw opacity=1 ][fill={rgb, 255:red, 0; green, 0; blue, 0 }  ,fill opacity=1 ] (391.25,40.61) .. controls (391.25,39.26) and (392.35,38.16) .. (393.7,38.16) .. controls (395.06,38.16) and (396.16,39.26) .. (396.16,40.61) .. controls (396.16,41.97) and (395.06,43.07) .. (393.7,43.07) .. controls (392.35,43.07) and (391.25,41.97) .. (391.25,40.61) -- cycle ;
	
	\draw  [color={rgb, 255:red, 0; green, 0; blue, 0 }  ,draw opacity=1 ][fill={rgb, 255:red, 0; green, 0; blue, 0 }  ,fill opacity=1 ] (441.25,40.61) .. controls (441.25,39.26) and (442.35,38.16) .. (443.7,38.16) .. controls (445.06,38.16) and (446.16,39.26) .. (446.16,40.61) .. controls (446.16,41.97) and (445.06,43.07) .. (443.7,43.07) .. controls (442.35,43.07) and (441.25,41.97) .. (441.25,40.61) -- cycle ;
	
	\draw    (185,40) -- (209,40) ;
	\draw [shift={(211,40)}, rotate = 180] [color={rgb, 255:red, 0; green, 0; blue, 0 }  ][line width=0.5]    (4.37,-1.32) .. controls (2.78,-0.56) and (1.32,-0.12) .. (0,0) .. controls (1.32,0.12) and (2.78,0.56) .. (4.37,1.32)   ;
	
	\draw    (265,40) -- (291,40) ;
	\draw [shift={(293,40)}, rotate = 180] [color={rgb, 255:red, 0; green, 0; blue, 0 }  ][line width=0.5]    (4.37,-1.32) .. controls (2.78,-0.56) and (1.32,-0.12) .. (0,0) .. controls (1.32,0.12) and (2.78,0.56) .. (4.37,1.32)   ;
	
	\draw    (399,40) -- (433,40) ;
	\draw [shift={(435,40)}, rotate = 180] [color={rgb, 255:red, 0; green, 0; blue, 0 }  ][line width=0.5]    (4.37,-1.32) .. controls (2.78,-0.56) and (1.32,-0.12) .. (0,0) .. controls (1.32,0.12) and (2.78,0.56) .. (4.37,1.32)   ;
	
	\draw  [draw opacity=0] (390.93,49.68) .. controls (388.08,53.52) and (383.5,56) .. (378.33,56) .. controls (369.68,56) and (362.67,49.03) .. (362.67,40.43) .. controls (362.67,31.82) and (369.68,24.85) .. (378.33,24.85) .. controls (384.08,24.85) and (389.1,27.93) .. (391.83,32.51) -- (378.33,40.43) -- cycle ; \draw   (390.93,49.68) .. controls (388.08,53.52) and (383.5,56) .. (378.33,56) .. controls (369.68,56) and (362.67,49.03) .. (362.67,40.43) .. controls (362.67,31.82) and (369.68,24.85) .. (378.33,24.85) .. controls (384.08,24.85) and (389.1,27.93) .. (391.83,32.51) ;
	
	\draw    (96.28,44.17) -- (92.98,44.69) ;
	\draw [shift={(91,45)}, rotate = 351.07] [color={rgb, 255:red, 0; green, 0; blue, 0 }  ][line width=0.5]    (4.37,-1.32) .. controls (2.78,-0.56) and (1.32,-0.12) .. (0,0) .. controls (1.32,0.12) and (2.78,0.56) .. (4.37,1.32)   ;
	
	\draw    (85.88,52.66) ;
	\draw [shift={(84.79,51.57)}, rotate = 405] [color={rgb, 255:red, 0; green, 0; blue, 0 }  ][line width=0.5]    (4.37,-1.32) .. controls (2.78,-0.56) and (1.32,-0.12) .. (0,0) .. controls (1.32,0.12) and (2.78,0.56) .. (4.37,1.32)   ;
	
	\draw    (75.09,51.47) -- (76.31,48.97) ;
	\draw [shift={(77.19,47.17)}, rotate = 476.04] [color={rgb, 255:red, 0; green, 0; blue, 0 }  ][line width=0.5]    (4.37,-1.32) .. controls (2.78,-0.56) and (1.32,-0.12) .. (0,0) .. controls (1.32,0.12) and (2.78,0.56) .. (4.37,1.32)   ;
	
	\draw  [dash pattern={on 0.84pt off 2.51pt}]  (220,40) -- (259,40) ;
	
	\draw    (386,54) -- (389.44,51.25) ;
	\draw [shift={(391,50)}, rotate = 501.34] [color={rgb, 255:red, 0; green, 0; blue, 0 }  ][line width=0.5]    (7.65,-2.3) .. controls (4.86,-0.97) and (2.31,-0.21) .. (0,0) .. controls (2.31,0.21) and (4.86,0.98) .. (7.65,2.3)   ;
	
	\draw (-7,165) node [anchor=north west][inner sep=0.5pt]   [align=left] {$ $};
	
	\draw (13.65,34.48) node [anchor=north west][inner sep=0.5pt]   [align=left] {$\displaystyle R_{n}$:$ $};
	
	\draw (146,33) node [anchor=north west][inner sep=0.5pt]   [align=left] {$\displaystyle A_{n}$:};
	
	\draw (322,33) node [anchor=north west][inner sep=0.5pt]   [align=left] {$\displaystyle E_{T}$:};
	
	\draw (97.76,0.32) node [anchor=north west][inner sep=0.5pt]   [align=left] {$\displaystyle e_{1}$};
	
	\draw (119.74,34.38) node [anchor=north west][inner sep=0.5pt]   [align=left] {$\displaystyle e_{2}$};
	
	\draw (97.81,75.04) node [anchor=north west][inner sep=0.5pt]   [align=left] {$\displaystyle e_{n}$};
	
	\draw (175,27) node [anchor=north west][inner sep=0.5pt]   [align=left] {$\displaystyle v_{1}$};
	
	\draw (214,27) node [anchor=north west][inner sep=0.5pt]   [align=left] {$\displaystyle v_{2}$};
	
	\draw (294,27) node [anchor=north west][inner sep=0.5pt]   [align=left] {$\displaystyle v_{n}$};
	
	\draw (188,42) node [anchor=north west][inner sep=0.5pt]   [align=left] {$\displaystyle e_{1}$};
	
	\draw (263,42) node [anchor=north west][inner sep=0.5pt]   [align=left] {$\displaystyle e_{n-1}$};
	
	\draw (352,36) node [anchor=north west][inner sep=0.5pt]   [align=left] {$\displaystyle e$};
	
	\draw (409,44) node [anchor=north west][inner sep=0.5pt]   [align=left] {$\displaystyle f$};
	
	\draw (395,27) node [anchor=north west][inner sep=0.5pt]   [align=left] {$\displaystyle u$};
	
	\draw (440,27) node [anchor=north west][inner sep=0.5pt]   [align=left] {$\displaystyle v$};
	
	\draw (66.75,36.58) node [anchor=north west][inner sep=0.5pt]   [align=left] {$\displaystyle v$};
	
	\draw (246,26) node [anchor=north west][inner sep=0.5pt]   [align=left] {$\displaystyle v_{n-1}$};
	\end{tikzpicture}
\end{center}
\vskip -2cm
The paper is organized as follows. Section 2 is devoted to a study of some types of simple modules over Leavitt path algebras and is made in preparation for further use. Section 3 presents the main results of this paper. In this section, we shall prove in Theorem \ref{theorem_3.6} that the non-commutative Leavitt path algebra $L_K(E)$ always contains non-cyclic free subgroups if $K$ is a non-absolute field. As a consequence of this theorem, we show that the same result also holds for Cohn path algebra $C_K(E)$ (Corollary~\ref{corollary_3.8}). The next three sections, Sections 4, 5 and 6, provide the descriptions of the generators of the non-cyclic free subgroups in $L_K(E)$ in term of the base graph $E$. Such the descriptions are obtained by examining three types of primitive ideals in $L_K(E)$. Finally, in Section 7, we shall determine the multiplicative group of some special types of Leavitt path algebras including commutative, artinian, noetherian Leavitt path algebras, and the Toeplitz $K$-algebra $\mathscr{T}_K$.

In the end, let $R$ be a ring or algebra with identity. Then, the symbol $R^{\times}$ stands for the multiplicative group of $R$. For a field $K$, we denote by ${\rm M}_n(K)$ the $n\times n$ full matrix algebra over $K$, and by $\{f_{ij}\;|\;1\leq i,j\leq n\}$ the standard matrix units in ${\rm M}_n(K)$.

\section{Simple modules over Leavitt path algebras}
Before going to the main study of this paper, we need to examine some classes of simple modules over  Leavitt path algebras introduced by Chen and Rangaswamy.

(A): \textit{Chen's module $V_{[p]}$ defined by an infinite path $p$}:
Let $p=e_1e_2\dots e_n\dots$ be an infinite path, and $n\geq1$ an integer. In  \cite{Pa_chen_2012}, Chen defines $\tau_{\leq n}(p):=e_1\dots e_n$ and $\tau_{> n}(p):=e_{n+1}e_{n+2}\dots$. In addition, we set $\tau_{\leq 0}(p):=s(p)$ and $\tau_{> 0}(p):=p$. Two infinite paths $p$ and $q$ are said to be \textit{tail-equivalent} if there exist positive integers $m,n\geq1$ such that $\tau_{>m}(p)=\tau_{>n}(q)$. A simple checking shows that this is an equivalence relation. For an infinite path $p$, the symbol $[p]$ stands for the equivalence class of all paths which are  tail-equivalent to $p$. If $p=g^{\infty}:=ggg\dots$, where $g$ is some  closed path, then we say that $p$ is \textit{rational}. Following Rangaswamy \cite{Pa_ara_rangas_2014}, if an infinite path $p$ is tail-equivalent to the rational path $c^\infty$, where $c$ is a cycle in $E$, then we say that $p$ \textit{ends in a cycle.} Also, a cycle $c$ is said to be an \textit{exclusive cycle} if no vertex on $c$ is the base of a different cycle (other than the cyclic permutation of $c$ based at that vertex).

For an infinite path $p$, Chen defines 
$$
V_{[p]}=\bigoplus_{q\in[p]}Kq,
$$ 
a $K$-vector space which has a $K$-basis $\{q:q\in [p]\}$. It was shown by Chen in \cite{Pa_chen_2012} that the following rules make $V_{[p]}$ into a simple left $L_K(E)$-module.

\begin{enumerate}
	\item[] $v\cdot q=q$ or 0 according to $v=s(q)$ or not;
	\item[] $e\cdot q =eq$ or 0 according to $r(e)=s(q)$ or not;
	\item[] $e^*\cdot q=\tau_{>1}(q)$ or $0$ according to $q=eq'$ or not.
\end{enumerate}

(B): \textit{Chen's module $\mathbf{N}_w$ defined by a sink $w$}: Let $w$ be a sink in the graph $E$. Following Chen, we denote by $\mathbf{N}_w$ the $K$-vector space whose basis consists of all finite paths ending in $w$; that is,
$$
\mathbf{N}_w=\bigoplus_{\substack{q\in{\rm Path}(E),\\ r(q)=w}}Kq.
$$
Then, $\mathbf{N}_w$ becomes a simple left $L_K(E)$-module by defining an action on its basis elements in the same way we have done for $V_{[p]}$ (with addition that $e^*\cdot w=0$ for all $e\in E^1$). From  \cite[Theorem 3.3(1) and Theorem 3.7(1)]{Pa_chen_2012}, we can state the following result.
\begin{lemma}\label{lemma_2.1}
	Let $p$ be an infinite path, and $w$ a sink in $E$. Let $V_{[p]}$ and $\mathbf{N}_w$ be the $L_K(E)$-modules defined above. Then, 
	$${\rm End}_{L_K(E)}(V_{[p]})\cong K \hbox{ and } \;\;{\rm End}_{L_K(E)}(\mathbf{N}_w)\cong K.$$
\end{lemma}

(C): \textit{Rangaswamy's twisted module $V_{[p]}^f$}: Let $A$ be an algebra and $M$ an $A$-module. Given an automorphism $\sigma$ of $A$, we can define the twisted $A$-module $M^\sigma$ as follows.
\begin{enumerate}
	\item[] $M^\sigma=M$ as a vector space over $A$,
	\item[] $a\cdot m^\sigma=(\sigma(a)m)^\sigma$.
\end{enumerate}
(Here, for an element $m\in M$, the symbol $m^{\sigma}$ stands for the corresponding element in $M^\sigma$.) One can easily check that ${\rm End}_A(M) = {\rm End}_A(M^\sigma)$. Moreover, the twisted module $M^\sigma$ is simple if and only if $M$ is simple. In \cite{Pa_chen_2012}, Chen has studied intensively the twisted modules $V_{[p]}^{\mathbf{a}}$ and $\mathbf{N}_w^{\mathbf{a}}$ of $V_{[p]}$ and $\mathbf{N}_w$ respectively, where ${\mathbf{a}}$ is an appropriate automorphism of $L_K(E)$.\\

As a modification of Chen's construction, in \cite{Pa_rangas_2013}, Rangaswamy introduced the new simple $L_K(E)$-module $V_{[p]}^f$, which can be described briefly as follows. Let $f(x)=1+a_1x+\dots+a_nx^n$, $n\geq1$, be an irreducible polynomial in $K[x]$ and let $c=e_1e_2\dots e_m$ be an exclusive cycle in $E$. Set $p=c^\infty$. Let $K'=K[x]/(f(x))$, which is a field extension of $K$. Then, we can form the Leavitt path algebra $L_{K'}(E)$ over the field $K'$. Let $\bar{x}$ be the image of $x$ in $K'$. It is clear that $\bar{x}\ne 0$ in $K'$, and so $\bar{x}$ is invertible in $K'$. Therefore, we can define a map $\sigma: L_{K'}(E)\to L_{K'}(E)$ which sends 

\begin{enumerate}
	\item[] $v\to v$ for $v\in E^0$,
	\item[] $e\to e$ and $e^*\to e^*$ for $e\in E^1$ with $e\ne e_1$,
	\item[] $e_1\to \bar{x}e_1$ and $e_1^*\to \bar{x}^{-1}e_1^*$.
\end{enumerate}
It is a simple matter to check that $\sigma$ is an automorphism of $L_{K'}(E)$, and so we can define the twisted $L_{K'}(E)$-module $V_{[p]}^\sigma$ of the $L_{K'}(E)$-module $V_{[p]}$. In \cite{Pa_ara_rangas_2014}, Rangaswamy denotes by $V_{[p]}^f$ the $L_K(E)$-module obtained by restricting scalars on $V_{[p]}$ from $L_{K'}(E)$ to $L_K(E)$. Also, it was proved that $V_{[p]}^f$ is a simple $L_K(E)$-module. Now, for further using, we need to conculate ${\rm End}_{L_K(E)}(V_{[p]}^f)$.
\begin{lemma}\label{lemma_2.2}
	Let $c=e_1e_2\dots e_m$ be an exclusive cycle in a graph $E$ and $p=c^\infty$. Let $f(x)=1+a_1x+\dots+a_nx^n$, $n\geq1$, be an irreducible polynomial in $K[x]$ and let $K'=K[x]/(f(x))$. Consider the simple $L_K(E)$-module $V_{[p]}^f$. Then, 
	$${\rm End}_{L_K(E)}(V_{[p]}^f)\cong K'.$$
\end{lemma}
\begin{proof}
	For a sake of simplicity, we set $R=L_K(E)$ and $R'=L_{K'}(E)$. Let us consider the simple $R'$-module $V_{[p]}$ defined by the infinite path $p$. In view of Lemma~\ref{lemma_2.1}, we conclude that ${\rm End}_{R'}(V_{[p]})\cong K'.$ Let $\sigma$ be the automorphism of $R'$ defined as above.  Then  ${\rm End}_{R'}(V_{[p]}^\sigma)\cong K'.$ It follows from $R\subseteq R'$ that ${\rm End}_{R'}(V_{[p]}^\sigma)
	\subseteq {\rm End}_{R}(V_{[p]}^f)$. For the inverse inclusion, we take $\varphi\in{\rm End}_R(V_{[p]}^f)$ and write
	$$
	K'=K\oplus \bar{x}\oplus\dots\oplus K\bar{x}^{n-1},$$ where $\bar{x}$ is the image of $x$ in $K'$. Then $\varphi$ can be viewed as an element of ${\rm End}_{R'}(V_{[p]}^\sigma)$ by extending the scalars on $V_{[p]}^\sigma$ from $R$ to $R'$ with $\varphi(\bar{x})=\bar{x}$. It follows that ${\rm End}_R(V_{[p]}^f)
	\subseteq {\rm End}_{R'}(V_{[p]}^\sigma)$. Therefore, ${\rm End}_R(V_{[p]}^f)={\rm End}_{R'}(V_{[p]}^\sigma)\cong K'$.
\end{proof}

(D): \textit{Rangaswamy's module $\mathbf{S}_{v\infty}$ defined by an infinite emitter $v$}: Let $v$ be an infinite emitter in $E$. In \cite{Pa_rangas_2015}, Rangaswamy defines $\mathbf{S}_{v\infty}$ to be the $K$-vector space whose basis is the set $B=\{p: p \hbox{ is a path in } E \hbox{ with } r(p)=v\}$; that is,
$$
\mathbf{S}_{v\infty}=\bigoplus_{\substack{p\in{\rm Path}(E), \\ r(p)=v}}Kp.
$$
For each vertex $u$ and each edge $e$ in $E$, he defines linear transformations $P_u$, $S_e$ and $S_{e^*}$ on $\mathbf{S}_{v\infty}$ as follows. For any $p\in B$, set

\begin{align*}
	P_u(p) =    \begin{cases} 
						p        & \text{if } u=s(p) \\
						0        & \text{otherwise }
					  \end{cases} ,                        					
& \hspace{1cm} S_e(p) =    \begin{cases} 
												ep      & \text{if } r(e)=s(p) \\
												0        & \text{otherwise }
				   							\end{cases},\\
S_{e^*}(p) =    \begin{cases} 
							p'       & \text{if } p=ep' \\
							0        & \text{otherwise }
						\end{cases},										
& \hspace{1cm} S_{e^*}(v) =   0.
\end{align*}
The mapping $\phi: L_K(E) \to {\rm End}_K(\mathbf{S}_{v\infty})$, which sends $u$ to $P_u$, $e$ to $S_e$ and $e^*$ to $S_{e^*}$, is an algebra homomorphism. Therefore, $\mathbf{S}_{v\infty}$ can be viewed as a left module over $L_K(E)$ via $\phi$.  In \cite{Pa_rangas_2015}, Rangaswamy proved that $\mathbf{S}_{v\infty}$ is a simple left $L_K(E)$-module. For further use, we need to compute ${\rm End}_K(\mathbf{S}_{v\infty})$.
\begin{lemma}\label{lemma_2.3}
	If  $v$ is an infinite emitter in a graph $E$, then $
	{\rm End}_{L_K(E)}\left( \mathbf{S}_{v\infty}\right) \cong K.$
\end{lemma}
\begin{proof}
	Assume that $p, q\in B$. It is easy to check that $p^*q=\delta_{pq}v$.  Consider a nonzero element $\varphi\in {\rm End}_{L_K(E)}(S_{v\infty})$. Then, $\varphi$ is an automorphism because $S_{v\infty}$ is a simple $L_K(E)$-module. For an arbitrary $p\in B$, let
	$$
	\varphi(p)=\sum_{i=1}^\ell k_iq_i.
	$$
	If $q_i\ne p$ for all $i$, then
	$$
	0\ne\varphi(v)=\varphi(p^*p)=p^*\varphi(p)=p^*\varphi(p)=p^*\Big(\sum_{i=1}^\ell k_iq_i\Big)=0,
	$$
	a contradiction. Hence, $\varphi(p)=kp$ for some $k\in K$. We claim that $k$ depends only on $\varphi$, not on $p$. Indeed, assume that $q\in B$, $q\ne p$. Then, $\varphi(q)=k'q$ for some $k'\in K$. Since $r(p)=r(q)=v$, we have
	$$
	kp=\varphi(p)=\varphi(pv)=\varphi(pq^*q)=pq^*\varphi(q)=pq^*(k'q)=k'p(q^*q)=k'p,
	$$
	and it follows that $k=k'$. The claim is shown. Therefore, we conclude that for every nonzero element $\varphi\in {\rm End}_{L_K(E)}(S_{v\infty})$, there exists a unique element, say, $\lambda_{\varphi}\in K$ such that $\varphi(p)=\lambda_{\varphi}p$ for all $p\in B$. If we set additionally  $\lambda_0=0$, then, clearly, the assignment $\varphi\mapsto \lambda_{\varphi}$ defines an isomorphism between ${\rm End}_{L_K(E)}(S_{v\infty})$ and $K$. 
\end{proof}

\begin{proposition}\label{proposition_2.4}
	Let $E$ be an arbitrary graph, $K$ an arbitrary field, and $P$ any primitive ideal of $L_K(E)$. Then, there exists a simple left module $V$ which is one of the types $V_{[p]}$, $\mathbf{N}_w$, $V_{[p]}^f$, and $\mathbf{S}_{v\infty}$ such that ${\rm Ann}_{L_k(E)}(V)=P$ and  ${\rm End}_{L_K(E)} V$ is a field containing $K$. 
\end{proposition}
\begin{proof}
	The existence of $V$ follows from \cite[Propositions 2.6 and 2.7]{Pa_rangas_2015} and \cite[Theorem~ 3.9]{Pa_ara_rangas_2014}. For the second assertion, we note that, according to Lemmas \ref{lemma_2.1}, \ref{lemma_2.2} and \ref{lemma_2.3}, we conclude that ${\rm End}_{L_K(E)} V$ is isomorphic to $K$ or $K[x]/(p(x))$, where $p(x)$ is an irreducible polynomial in $K[x]$.
\end{proof}

\section{Non-cyclic free subgroups in Leavitt path algebras }

For a graph $E$ and an arbitrary field $K$, recall that the Leavitt path algebra $L_K(E)$ is unital if and only if the vertex set $E^0$ is finite \cite[Lemma 1.2.12 (iv)]{Pa_abrams-pino_2005}. The main aim of this section is to prove that if $K$ is not a locally finite field and $E^0$ is finite, then the multiplicative group of the Leavitt path algebra $L_K(E)$ contains a non-cyclic free subgroup in case $L_K(E)$ is non-commutative. Recall that a field $K$ is \textit{locally finite} (or \textit{absolute}) if every its finite subset generates over the prime subfield a finite field. It is clear that a field $K$ with the prime subfield $\mathbb{P}$ is locally finite if and only if $\mathbb{P}$ is finite and $K$ is algebraic over $\mathbb{P}$. In this paper, we prefer to use the terminology ``absolute field'' than ``locally finite field''. Fields of zero characteristic are trivial examples of non-absolute fields. Other examples of such a kind are uncountable fields (see \cite[Corollary B-2.41, p. 343]{Bo_Rotman_2015}). To ensure the existence of the identity in $L_K(E)$, in the remaining part of this paper, we always assume that $E$ has a finite number of vertices; that is, $E^0$ is a finite set. The following lemma is somewhat elementary (so its proof will be omitted), but plays an important role in solving the problem of this section.

\begin{lemma}\label{lemma_3.1}
	Let $R$ and $S$ be unital rings, and let $\varphi: R\to S$ be a surjective  homomorphism. Assume that $\left\langle x, y\right\rangle $ is a non-cyclic free subgroup of $S^\times$ generated by $x$ and $y$. If $a\in \varphi^{-1}(x)$ and $b\in \varphi^{-1}(y)$ such that $a, b \in R^\times$, then $\left\langle a, b\right\rangle $ is a non-cyclic free subgroup of $R^\times$ and $\left\langle a, b\right\rangle \cong \left\langle x, y\right\rangle$.
\end{lemma}

\begin{proposition}\label{proposition_3.2}
	Let $E$ be a graph, and $K$ an arbitrary field. If $L_K(E)$ is non-commutative, then there exists a primitive ideal $P$ of $L_K(E)$ such that $L_K(E)/P$ is a (unital) non-commutative primitive ring.
\end{proposition}
\begin{proof}
	Recall that for any graph $E$ and any field $K$, the Leavitt path algebra $L_K(E)$ is semiprimitive \cite[Proposition 2.3.2]{Pa_abrams-pino_2005}.	
	Let $\{P_i\mid i\in I\}$ be the set of all primitive ideals in $L_K(E)$.	Then, $\bigcap_{i\in I}P_i=0$, and there exists an injective ring homomorphism
	$$\varepsilon: L_K(E)\longrightarrow\prod_{i\in I}(L_K(E)/P_i).$$
	
	Since $L_K(E)$ is non-commutative, so is $\prod_{i\in I}(L_K(E)/P_i)$. It
	follows that there exists $i_0\in I$ such that $L_K(E)/P_{i_0}$ is non-commutative. By setting $P=P_{i_0}$, the proof of the lemma is now complete.
\end{proof}

In view of Proposition \ref{proposition_3.2}, we know that if $L_K(E)$ is a non-commutative Leavitt path algebra, then there always exists a primitive ideal $P$ for which $L_K(E)/P$ is non-commutative. One would like to know that what type of $P$ is. According to 
\cite[Theorem 4.3]{Pa_rangas_2013}, there are exactly three types of non-zero primitive ideals in a Leavitt path algebra $L_K(E)$. For the convenience of our use in this paper, we restate this result in the following theorem.
\begin{theorem}\label{theorem_3.3}
	Let $E$ be a graph,  and $K$ a field. For a given ideal $P$ of $L_K(E)$, set $H=P\cap E^0$. Then, $P$ is primitive if and only if $P$ is one of the following types:
	
	\begin{enumerate}[font=\normalfont]
		\item[I:] $P$ is a  graded ideal of the form $I(H,B_H\backslash\{w\})$, which is necessarily non-zero, where $w\in B_H$ and $M(w)=E^0\backslash H$.
		\item[II:] $P$ is a graded ideal of the form $I(H,B_H)$, which is possibly zero, and $E\backslash(H,B_H)$ is downward directed, satisfying the Condition (L) and the Countable Separation Property.
		\item[III:] $P=I(H, B_H, f(c))$, which is necessarily  non-zero, where $c$ is an exclusive cycle based at a vertex $u$, $E^0\backslash H=M(u)$, and $f(x)$ is an irreducible polynomial in $K[x, x^{-1}]$.
	\end{enumerate}
\end{theorem}

The following examples show that we can always find an appropriate graph $E$ for which $L_K(E)$ contains a primitive ideal $P$ of either types given in Theorem \ref{theorem_3.3} such that $L_{K(E)}/P$ is non-commutative.
\begin{example}\label{example_3}
	Consider again the graph given in Example \ref{example_1}. The set $H=\{v_1,v_2\}$ is a hereditary saturated subset of $E^0$ and $B_H=\{v, w\}$. It is clear that $M(w)=E^0\backslash H$, which implies that $P:=I(H, B_H\backslash \{w\})$ is a primitive graded ideal of type I  of $L_K(E)$. Moreover, if we set $F:=E\backslash(H, B_H\backslash\{w\})$, then $F$ is exactly the graph $E\backslash(H,S)$, which was described in Example \ref{example_2}, where $S=\{v\}$. 
	
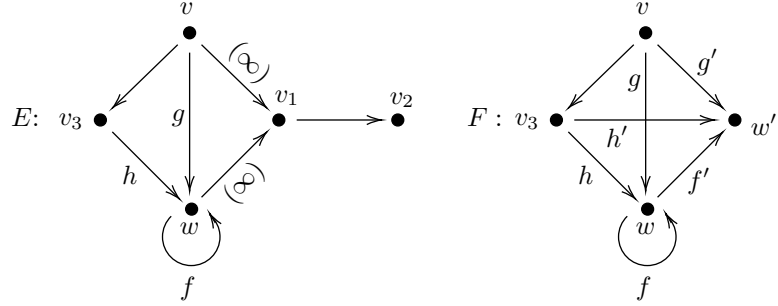
	\begin{center}
		\tikzset{every picture/.style={line width=0.5pt}} 
		
		\begin{tikzpicture}[x=0.75pt,y=0.75pt,yscale=-1,xscale=1]
		
		\draw  [color={rgb, 255:red, 0; green, 0; blue, 0 }  ,draw opacity=1 ][fill={rgb, 255:red, 0; green, 0; blue, 0 }  ,fill opacity=1 ] (93,135) .. controls (93,133.34) and (94.34,132) .. (96,132) .. controls (97.66,132) and (99,133.34) .. (99,135) .. controls (99,136.66) and (97.66,138) .. (96,138) .. controls (94.34,138) and (93,136.66) .. (93,135) -- cycle ;
		
		\draw  [draw opacity=0] (104.92,137.71) .. controls (108.21,140.36) and (110.31,144.43) .. (110.31,148.99) .. controls (110.31,156.99) and (103.82,163.48) .. (95.82,163.48) .. controls (87.81,163.48) and (81.32,156.99) .. (81.32,148.99) .. controls (81.32,144.9) and (83.01,141.22) .. (85.72,138.58) -- (95.82,148.99) -- cycle ; \draw   (104.92,137.71) .. controls (108.21,140.36) and (110.31,144.43) .. (110.31,148.99) .. controls (110.31,156.99) and (103.82,163.48) .. (95.82,163.48) .. controls (87.81,163.48) and (81.32,156.99) .. (81.32,148.99) .. controls (81.32,144.9) and (83.01,141.22) .. (85.72,138.58) ;
		
		\draw    (108,142) -- (106.03,138.71) ;
		\draw [shift={(105,137)}, rotate = 419.03999999999996] [color={rgb, 255:red, 0; green, 0; blue, 0 }  ][line width=0.5]    (7.65,-2.3) .. controls (4.86,-0.97) and (2.31,-0.21) .. (0,0) .. controls (2.31,0.21) and (4.86,0.98) .. (7.65,2.3)   ;
		
		\draw  [color={rgb, 255:red, 0; green, 0; blue, 0 }  ,draw opacity=1 ][fill={rgb, 255:red, 0; green, 0; blue, 0 }  ,fill opacity=1 ] (92,46) .. controls (92,44.34) and (93.34,43) .. (95,43) .. controls (96.66,43) and (98,44.34) .. (98,46) .. controls (98,47.66) and (96.66,49) .. (95,49) .. controls (93.34,49) and (92,47.66) .. (92,46) -- cycle ;
		
		\draw  [color={rgb, 255:red, 0; green, 0; blue, 0 }  ,draw opacity=1 ][fill={rgb, 255:red, 0; green, 0; blue, 0 }  ,fill opacity=1 ] (137,90) .. controls (137,88.34) and (138.34,87) .. (140,87) .. controls (141.66,87) and (143,88.34) .. (143,90) .. controls (143,91.66) and (141.66,93) .. (140,93) .. controls (138.34,93) and (137,91.66) .. (137,90) -- cycle ;
		
		\draw  [color={rgb, 255:red, 0; green, 0; blue, 0 }  ,draw opacity=1 ][fill={rgb, 255:red, 0; green, 0; blue, 0 }  ,fill opacity=1 ] (47,90) .. controls (47,88.34) and (48.34,87) .. (50,87) .. controls (51.66,87) and (53,88.34) .. (53,90) .. controls (53,91.66) and (51.66,93) .. (50,93) .. controls (48.34,93) and (47,91.66) .. (47,90) -- cycle ;
		
		\draw    (89,52) -- (57.44,82.61) ;
		\draw [shift={(56,84)}, rotate = 315.88] [color={rgb, 255:red, 0; green, 0; blue, 0 }  ][line width=0.5]    (7.65,-2.3) .. controls (4.86,-0.97) and (2.31,-0.21) .. (0,0) .. controls (2.31,0.21) and (4.86,0.98) .. (7.65,2.3)   ;
		
		\draw    (101,52) -- (132.56,82.61) ;
		\draw [shift={(134,84)}, rotate = 224.12] [color={rgb, 255:red, 0; green, 0; blue, 0 }  ][line width=0.5]    (7.65,-2.3) .. controls (4.86,-0.97) and (2.31,-0.21) .. (0,0) .. controls (2.31,0.21) and (4.86,0.98) .. (7.65,2.3)   ;
		
		\draw    (56,96) -- (88.59,128.59) ;
		\draw [shift={(90,130)}, rotate = 225] [color={rgb, 255:red, 0; green, 0; blue, 0 }  ][line width=0.5]    (7.65,-2.3) .. controls (4.86,-0.97) and (2.31,-0.21) .. (0,0) .. controls (2.31,0.21) and (4.86,0.98) .. (7.65,2.3)   ;
		
		\draw    (101,129) -- (132.59,97.41) ;
		\draw [shift={(134,96)}, rotate = 495] [color={rgb, 255:red, 0; green, 0; blue, 0 }  ][line width=0.5]    (7.65,-2.3) .. controls (4.86,-0.97) and (2.31,-0.21) .. (0,0) .. controls (2.31,0.21) and (4.86,0.98) .. (7.65,2.3)   ;
		
		\draw    (95,56) -- (95,124) ;
		\draw [shift={(95,126)}, rotate = 270] [color={rgb, 255:red, 0; green, 0; blue, 0 }  ][line width=0.5]    (7.65,-2.3) .. controls (4.86,-0.97) and (2.31,-0.21) .. (0,0) .. controls (2.31,0.21) and (4.86,0.98) .. (7.65,2.3)   ;
		
		\draw    (149,90) -- (190,90) ;
		\draw [shift={(192,90)}, rotate = 180] [color={rgb, 255:red, 0; green, 0; blue, 0 }  ][line width=0.5]    (7.65,-2.3) .. controls (4.86,-0.97) and (2.31,-0.21) .. (0,0) .. controls (2.31,0.21) and (4.86,0.98) .. (7.65,2.3)   ;
		
		\draw  [color={rgb, 255:red, 0; green, 0; blue, 0 }  ,draw opacity=1 ][fill={rgb, 255:red, 0; green, 0; blue, 0 }  ,fill opacity=1 ] (197,90) .. controls (197,88.34) and (198.34,87) .. (200,87) .. controls (201.66,87) and (203,88.34) .. (203,90) .. controls (203,91.66) and (201.66,93) .. (200,93) .. controls (198.34,93) and (197,91.66) .. (197,90) -- cycle ;
		
		\draw  [color={rgb, 255:red, 0; green, 0; blue, 0 }  ,draw opacity=1 ][fill={rgb, 255:red, 0; green, 0; blue, 0 }  ,fill opacity=1 ] (323,135) .. controls (323,133.34) and (324.34,132) .. (326,132) .. controls (327.66,132) and (329,133.34) .. (329,135) .. controls (329,136.66) and (327.66,138) .. (326,138) .. controls (324.34,138) and (323,136.66) .. (323,135) -- cycle ;
		
		\draw  [draw opacity=0] (334.92,137.71) .. controls (338.21,140.36) and (340.31,144.43) .. (340.31,148.99) .. controls (340.31,156.99) and (333.82,163.48) .. (325.82,163.48) .. controls (317.81,163.48) and (311.32,156.99) .. (311.32,148.99) .. controls (311.32,144.9) and (313.01,141.22) .. (315.72,138.58) -- (325.82,148.99) -- cycle ; \draw   (334.92,137.71) .. controls (338.21,140.36) and (340.31,144.43) .. (340.31,148.99) .. controls (340.31,156.99) and (333.82,163.48) .. (325.82,163.48) .. controls (317.81,163.48) and (311.32,156.99) .. (311.32,148.99) .. controls (311.32,144.9) and (313.01,141.22) .. (315.72,138.58) ;
		
		\draw    (338,142) -- (336.03,138.71) ;
		\draw [shift={(335,137)}, rotate = 419.03999999999996] [color={rgb, 255:red, 0; green, 0; blue, 0 }  ][line width=0.5]    (7.65,-2.3) .. controls (4.86,-0.97) and (2.31,-0.21) .. (0,0) .. controls (2.31,0.21) and (4.86,0.98) .. (7.65,2.3)   ;
		
		\draw  [color={rgb, 255:red, 0; green, 0; blue, 0 }  ,draw opacity=1 ][fill={rgb, 255:red, 0; green, 0; blue, 0 }  ,fill opacity=1 ] (322,46) .. controls (322,44.34) and (323.34,43) .. (325,43) .. controls (326.66,43) and (328,44.34) .. (328,46) .. controls (328,47.66) and (326.66,49) .. (325,49) .. controls (323.34,49) and (322,47.66) .. (322,46) -- cycle ;
		
		\draw  [color={rgb, 255:red, 0; green, 0; blue, 0 }  ,draw opacity=1 ][fill={rgb, 255:red, 0; green, 0; blue, 0 }  ,fill opacity=1 ] (367,90) .. controls (367,88.34) and (368.34,87) .. (370,87) .. controls (371.66,87) and (373,88.34) .. (373,90) .. controls (373,91.66) and (371.66,93) .. (370,93) .. controls (368.34,93) and (367,91.66) .. (367,90) -- cycle ;
		
		\draw  [color={rgb, 255:red, 0; green, 0; blue, 0 }  ,draw opacity=1 ][fill={rgb, 255:red, 0; green, 0; blue, 0 }  ,fill opacity=1 ] (277,90) .. controls (277,88.34) and (278.34,87) .. (280,87) .. controls (281.66,87) and (283,88.34) .. (283,90) .. controls (283,91.66) and (281.66,93) .. (280,93) .. controls (278.34,93) and (277,91.66) .. (277,90) -- cycle ;
		
		\draw    (319,52) -- (287.44,82.61) ;
		\draw [shift={(286,84)}, rotate = 315.88] [color={rgb, 255:red, 0; green, 0; blue, 0 }  ][line width=0.5]    (7.65,-2.3) .. controls (4.86,-0.97) and (2.31,-0.21) .. (0,0) .. controls (2.31,0.21) and (4.86,0.98) .. (7.65,2.3)   ;
		
		\draw    (331,52) -- (362.56,82.61) ;
		\draw [shift={(364,84)}, rotate = 224.12] [color={rgb, 255:red, 0; green, 0; blue, 0 }  ][line width=0.5]    (7.65,-2.3) .. controls (4.86,-0.97) and (2.31,-0.21) .. (0,0) .. controls (2.31,0.21) and (4.86,0.98) .. (7.65,2.3)   ;
		
		\draw    (286,96) -- (318.59,128.59) ;
		\draw [shift={(320,130)}, rotate = 225] [color={rgb, 255:red, 0; green, 0; blue, 0 }  ][line width=0.5]    (7.65,-2.3) .. controls (4.86,-0.97) and (2.31,-0.21) .. (0,0) .. controls (2.31,0.21) and (4.86,0.98) .. (7.65,2.3)   ;
		
		\draw    (331,129) -- (362.59,97.41) ;
		\draw [shift={(364,96)}, rotate = 495] [color={rgb, 255:red, 0; green, 0; blue, 0 }  ][line width=0.5]    (7.65,-2.3) .. controls (4.86,-0.97) and (2.31,-0.21) .. (0,0) .. controls (2.31,0.21) and (4.86,0.98) .. (7.65,2.3)   ;
		
		\draw    (325,56) -- (325,124) ;
		\draw [shift={(325,126)}, rotate = 270] [color={rgb, 255:red, 0; green, 0; blue, 0 }  ][line width=0.5]    (7.65,-2.3) .. controls (4.86,-0.97) and (2.31,-0.21) .. (0,0) .. controls (2.31,0.21) and (4.86,0.98) .. (7.65,2.3)   ;
		
		\draw    (289,90) -- (358,90) ;
		\draw [shift={(360,90)}, rotate = 180] [color={rgb, 255:red, 0; green, 0; blue, 0 }  ][line width=0.5]    (7.65,-2.3) .. controls (4.86,-0.97) and (2.31,-0.21) .. (0,0) .. controls (2.31,0.21) and (4.86,0.98) .. (7.65,2.3)   ;
		
		\draw (-27,165) node [anchor=north west][inner sep=0.5pt]   [align=left] {$ $};
		
		\draw (4,83) node [anchor=north west][inner sep=0.5pt]   [align=left] {$\displaystyle E$:};
		
		\draw (89,141) node [anchor=north west][inner sep=0.5pt]   [align=left] {$\displaystyle w$};
		
		\draw (89,168) node [anchor=north west][inner sep=0.5pt]   [align=left] {$\displaystyle f$};
		
		\draw (85,84) node [anchor=north west][inner sep=0.5pt]   [align=left] {$\displaystyle g$};
		
		\draw (60,112) node [anchor=north west][inner sep=0.5pt]   [align=left] {$\displaystyle h$};
		
		\draw (121.46,46.74) node [anchor=north west][inner sep=0.5pt]  [rotate=-45] [align=left] {$\displaystyle ( \infty $)};
		
		\draw (137.26,118.46) node [anchor=north west][inner sep=0.5pt]  [rotate=-135] [align=left] {$\displaystyle ( \infty $)};
		
		\draw (28,86) node [anchor=north west][inner sep=0.5pt]   [align=left] {$\displaystyle v_{3}$};
		
		\draw (137,74) node [anchor=north west][inner sep=0.5pt]   [align=left] {$\displaystyle v_{1}$};
		
		\draw (195,75) node [anchor=north west][inner sep=0.5pt]   [align=left] {$\displaystyle v_{2}$};
		
		\draw (89,30) node [anchor=north west][inner sep=0.5pt]   [align=left] {$\displaystyle v$};
		
		\draw (234,83) node [anchor=north west][inner sep=0.5pt]   [align=left] {$\displaystyle F:$};
		
		\draw (319,141) node [anchor=north west][inner sep=0.5pt]   [align=left] {$\displaystyle w$};
		
		\draw (319,168) node [anchor=north west][inner sep=0.5pt]   [align=left] {$\displaystyle f$};
		
		\draw (315,66) node [anchor=north west][inner sep=0.5pt]   [align=left] {$\displaystyle g$};
		
		\draw (290,112) node [anchor=north west][inner sep=0.5pt]   [align=left] {$\displaystyle h$};
		
		\draw (258,86) node [anchor=north west][inner sep=0.5pt]   [align=left] {$\displaystyle v_{3}$};
		
		\draw (377,86) node [anchor=north west][inner sep=0.5pt]   [align=left] {$\displaystyle w'$};
		
		\draw (319,30) node [anchor=north west][inner sep=0.5pt]   [align=left] {$\displaystyle v$};
		
		\draw (350,53) node [anchor=north west][inner sep=0.5pt]   [align=left] {$\displaystyle g'$};
		
		\draw (345,113) node [anchor=north west][inner sep=0.5pt]   [align=left] {$\displaystyle f'$};
		
		\draw (304,92) node [anchor=north west][inner sep=0.5pt]   [align=left] {$\displaystyle h'$};
		\end{tikzpicture}
	\end{center}
	
	In view of Lemma \ref{lemma_1.4}, we have  $L_K(E)/P\cong L_K(F)$ which is clearly non-commutative.
\end{example}
\begin{example}\label{exmaple_4}
	Consider the following graphs: 
	
	\begin{center}
		
		\tikzset{every picture/.style={line width=0.5pt}}
		
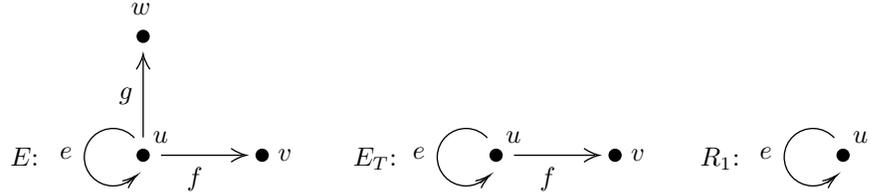
		\begin{tikzpicture}[x=0.75pt,y=0.75pt,yscale=-1,xscale=1]
		
		\draw  [color={rgb, 255:red, 0; green, 0; blue, 0 }  ,draw opacity=1 ][fill={rgb, 255:red, 0; green, 0; blue, 0 }  ,fill opacity=1 ] (80,135) .. controls (80,133.34) and (81.34,132) .. (83,132) .. controls (84.66,132) and (86,133.34) .. (86,135) .. controls (86,136.66) and (84.66,138) .. (83,138) .. controls (81.34,138) and (80,136.66) .. (80,135) -- cycle ;
		
		\draw  [color={rgb, 255:red, 0; green, 0; blue, 0 }  ,draw opacity=1 ][fill={rgb, 255:red, 0; green, 0; blue, 0 }  ,fill opacity=1 ] (140,135) .. controls (140,133.34) and (141.34,132) .. (143,132) .. controls (144.66,132) and (146,133.34) .. (146,135) .. controls (146,136.66) and (144.66,138) .. (143,138) .. controls (141.34,138) and (140,136.66) .. (140,135) -- cycle ;
		
		\draw    (92,135) -- (133,135) ;
		\draw [shift={(135,135)}, rotate = 180] [color={rgb, 255:red, 0; green, 0; blue, 0 }  ][line width=0.5]    (7.65,-3.43) .. controls (4.86,-1.61) and (2.31,-0.47) .. (0,0) .. controls (2.31,0.47) and (4.86,1.61) .. (7.65,3.43)   ;
		
		\draw  [color={rgb, 255:red, 0; green, 0; blue, 0 }  ,draw opacity=1 ][fill={rgb, 255:red, 0; green, 0; blue, 0 }  ,fill opacity=1 ] (80,75) .. controls (80,73.34) and (81.34,72) .. (83,72) .. controls (84.66,72) and (86,73.34) .. (86,75) .. controls (86,76.66) and (84.66,78) .. (83,78) .. controls (81.34,78) and (80,76.66) .. (80,75) -- cycle ;
		
		\draw    (83,126) -- (83,86) ;
		\draw [shift={(83,84)}, rotate = 450] [color={rgb, 255:red, 0; green, 0; blue, 0 }  ][line width=0.5]    (7.65,-3.43) .. controls (4.86,-1.61) and (2.31,-0.47) .. (0,0) .. controls (2.31,0.47) and (4.86,1.61) .. (7.65,3.43)   ;
		
		\draw  [draw opacity=0] (78.77,145.48) .. controls (76,148.67) and (71.87,150.63) .. (67.31,150.47) .. controls (59.31,150.19) and (53.05,143.48) .. (53.33,135.48) .. controls (53.61,127.48) and (60.32,121.22) .. (68.32,121.5) .. controls (72.4,121.64) and (76.03,123.46) .. (78.56,126.26) -- (67.82,135.99) -- cycle ; \draw   (78.77,145.48) .. controls (76,148.67) and (71.87,150.63) .. (67.31,150.47) .. controls (59.31,150.19) and (53.05,143.48) .. (53.33,135.48) .. controls (53.61,127.48) and (60.32,121.22) .. (68.32,121.5) .. controls (72.4,121.64) and (76.03,123.46) .. (78.56,126.26) ;
		
		\draw    (75,149) -- (78.44,146.25) ;
		\draw [shift={(80,145)}, rotate = 501.34] [color={rgb, 255:red, 0; green, 0; blue, 0 }  ][line width=0.5]    (7.65,-2.3) .. controls (4.86,-0.97) and (2.31,-0.21) .. (0,0) .. controls (2.31,0.21) and (4.86,0.98) .. (7.65,2.3)   ;
		
		\draw  [color={rgb, 255:red, 0; green, 0; blue, 0 }  ,draw opacity=1 ][fill={rgb, 255:red, 0; green, 0; blue, 0 }  ,fill opacity=1 ] (258,135) .. controls (258,133.34) and (259.34,132) .. (261,132) .. controls (262.66,132) and (264,133.34) .. (264,135) .. controls (264,136.66) and (262.66,138) .. (261,138) .. controls (259.34,138) and (258,136.66) .. (258,135) -- cycle ;
		
		\draw  [color={rgb, 255:red, 0; green, 0; blue, 0 }  ,draw opacity=1 ][fill={rgb, 255:red, 0; green, 0; blue, 0 }  ,fill opacity=1 ] (318,135) .. controls (318,133.34) and (319.34,132) .. (321,132) .. controls (322.66,132) and (324,133.34) .. (324,135) .. controls (324,136.66) and (322.66,138) .. (321,138) .. controls (319.34,138) and (318,136.66) .. (318,135) -- cycle ;
		
		\draw    (270,135) -- (311,135) ;
		\draw [shift={(313,135)}, rotate = 180] [color={rgb, 255:red, 0; green, 0; blue, 0 }  ][line width=0.5]    (7.65,-3.43) .. controls (4.86,-1.61) and (2.31,-0.47) .. (0,0) .. controls (2.31,0.47) and (4.86,1.61) .. (7.65,3.43)   ;
		
		\draw  [draw opacity=0] (256.77,145.48) .. controls (254,148.67) and (249.87,150.63) .. (245.31,150.47) .. controls (237.31,150.19) and (231.05,143.48) .. (231.33,135.48) .. controls (231.61,127.48) and (238.32,121.22) .. (246.32,121.5) .. controls (250.4,121.64) and (254.03,123.46) .. (256.56,126.26) -- (245.82,135.99) -- cycle ; \draw   (256.77,145.48) .. controls (254,148.67) and (249.87,150.63) .. (245.31,150.47) .. controls (237.31,150.19) and (231.05,143.48) .. (231.33,135.48) .. controls (231.61,127.48) and (238.32,121.22) .. (246.32,121.5) .. controls (250.4,121.64) and (254.03,123.46) .. (256.56,126.26) ;
		
		\draw    (253,149) -- (256.44,146.25) ;
		\draw [shift={(258,145)}, rotate = 501.34] [color={rgb, 255:red, 0; green, 0; blue, 0 }  ][line width=0.5]    (7.65,-2.3) .. controls (4.86,-0.97) and (2.31,-0.21) .. (0,0) .. controls (2.31,0.21) and (4.86,0.98) .. (7.65,2.3)   ;
		
		\draw  [color={rgb, 255:red, 0; green, 0; blue, 0 }  ,draw opacity=1 ][fill={rgb, 255:red, 0; green, 0; blue, 0 }  ,fill opacity=1 ] (433,135) .. controls (433,133.34) and (434.34,132) .. (436,132) .. controls (437.66,132) and (439,133.34) .. (439,135) .. controls (439,136.66) and (437.66,138) .. (436,138) .. controls (434.34,138) and (433,136.66) .. (433,135) -- cycle ;
		
		\draw  [draw opacity=0] (431.77,145.48) .. controls (429,148.67) and (424.87,150.63) .. (420.31,150.47) .. controls (412.31,150.19) and (406.05,143.48) .. (406.33,135.48) .. controls (406.61,127.48) and (413.32,121.22) .. (421.32,121.5) .. controls (425.4,121.64) and (429.03,123.46) .. (431.56,126.26) -- (420.82,135.99) -- cycle ; \draw   (431.77,145.48) .. controls (429,148.67) and (424.87,150.63) .. (420.31,150.47) .. controls (412.31,150.19) and (406.05,143.48) .. (406.33,135.48) .. controls (406.61,127.48) and (413.32,121.22) .. (421.32,121.5) .. controls (425.4,121.64) and (429.03,123.46) .. (431.56,126.26) ;
		
		\draw    (428,149) -- (431.44,146.25) ;
		\draw [shift={(433,145)}, rotate = 501.34] [color={rgb, 255:red, 0; green, 0; blue, 0 }  ][line width=0.5]    (7.65,-2.3) .. controls (4.86,-0.97) and (2.31,-0.21) .. (0,0) .. controls (2.31,0.21) and (4.86,0.98) .. (7.65,2.3)   ;
		
		\draw (-27,165) node [anchor=north west][inner sep=0.5pt]   [align=left] {$ $};
		
		\draw (15,130) node [anchor=north west][inner sep=0.5pt]   [align=left] {$\displaystyle E$:};
		
		\draw (87,122) node [anchor=north west][inner sep=0.5pt]   [align=left] {$\displaystyle u$};
		
		\draw (150,131) node [anchor=north west][inner sep=0.5pt]   [align=left] {$\displaystyle v$};
		
		\draw (76,57) node [anchor=north west][inner sep=0.5pt]   [align=left] {$\displaystyle w$};
		
		\draw (70,100) node [anchor=north west][inner sep=0.5pt]   [align=left] {$\displaystyle g$};
		
		\draw (104,140) node [anchor=north west][inner sep=0.5pt]   [align=left] {$\displaystyle f$};
		
		\draw (40,130) node [anchor=north west][inner sep=0.5pt]   [align=left] {$\displaystyle e$};
		
		\draw (188,130) node [anchor=north west][inner sep=0.5pt]   [align=left] {$\displaystyle E_{T}$:};
		
		\draw (265,122) node [anchor=north west][inner sep=0.5pt]   [align=left] {$\displaystyle u$};
		
		\draw (328,131) node [anchor=north west][inner sep=0.5pt]   [align=left] {$\displaystyle v$};
		
		\draw (282,140) node [anchor=north west][inner sep=0.5pt]   [align=left] {$\displaystyle f$};
		
		\draw (218,130) node [anchor=north west][inner sep=0.5pt]   [align=left] {$\displaystyle e$};
		
		\draw (363,130) node [anchor=north west][inner sep=0.5pt]   [align=left] {$\displaystyle R_{1}$:};
		
		\draw (440,122) node [anchor=north west][inner sep=0.5pt]   [align=left] {$\displaystyle u$};
		
		\draw (393,130) node [anchor=north west][inner sep=0.5pt]   [align=left] {$\displaystyle e$};
		\end{tikzpicture}
	\end{center}
	
	If we set $H_1=\{w\}$ and $H_2=\{v,w\}$, then $H_1$ and $H_2$ are hereditary saturated subsets of $E^0$. It is clear that $B_{H_1}=B_{H_2}=\varnothing$. Also, it is easy to see that $P_1:=I(H_1)$ and $P_2:=I(H_2)$ are primitive ideals of type II. We also have $E\backslash H_1=E_T$ and $E\backslash H_2=R_1$, which imply that $L_K(E)/ P_1\cong L_K(E_T) = \mathscr{T}_K$ and $L_K(E)/P_2\cong L_K(R_1) \cong K[x,x^{-1}]$. It follows that $L_K(E)/P_1$ is non-commutative and $L_K(E)/P_2$ is a commutative ring.
\end{example}
\begin{example}\label{exmaple_5}
	Consider the following graphs:\\
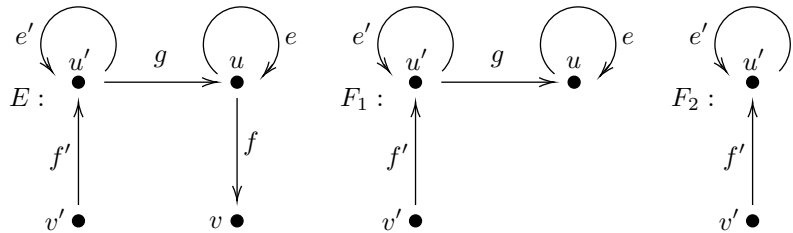
	\begin{center}
		
		\tikzset{every picture/.style={line width=0.5pt}}  
		
		\begin{tikzpicture}[x=0.75pt,y=0.75pt,yscale=-1,xscale=1]
		
		\draw  [color={rgb, 255:red, 0; green, 0; blue, 0 }  ,draw opacity=1 ][fill={rgb, 255:red, 0; green, 0; blue, 0 }  ,fill opacity=1 ] (113,135) .. controls (113,133.34) and (114.34,132) .. (116,132) .. controls (117.66,132) and (119,133.34) .. (119,135) .. controls (119,136.66) and (117.66,138) .. (116,138) .. controls (114.34,138) and (113,136.66) .. (113,135) -- cycle ;
		
		\draw  [draw opacity=0] (106.13,130.84) .. controls (101.78,127.36) and (99,122) .. (99,116) .. controls (99,105.51) and (107.51,97) .. (118,97) .. controls (128.49,97) and (137,105.51) .. (137,116) .. controls (137,121.25) and (134.87,126) .. (131.44,129.44) -- (118,116) -- cycle ; \draw   (106.13,130.84) .. controls (101.78,127.36) and (99,122) .. (99,116) .. controls (99,105.51) and (107.51,97) .. (118,97) .. controls (128.49,97) and (137,105.51) .. (137,116) .. controls (137,121.25) and (134.87,126) .. (131.44,129.44) ;
		
		\draw    (134,126) -- (130.28,130.46) ;
		\draw [shift={(129,132)}, rotate = 309.81] [color={rgb, 255:red, 0; green, 0; blue, 0 }  ][line width=0.5]    (7.65,-2.3) .. controls (4.86,-0.97) and (2.31,-0.21) .. (0,0) .. controls (2.31,0.21) and (4.86,0.98) .. (7.65,2.3)   ;
		
		\draw  [color={rgb, 255:red, 0; green, 0; blue, 0 }  ,draw opacity=1 ][fill={rgb, 255:red, 0; green, 0; blue, 0 }  ,fill opacity=1 ] (33,135) .. controls (33,133.34) and (34.34,132) .. (36,132) .. controls (37.66,132) and (39,133.34) .. (39,135) .. controls (39,136.66) and (37.66,138) .. (36,138) .. controls (34.34,138) and (33,136.66) .. (33,135) -- cycle ;
		
		\draw    (20,126) -- (21.8,128.4) ;
		\draw [shift={(23,130)}, rotate = 233.13] [color={rgb, 255:red, 0; green, 0; blue, 0 }  ][line width=0.5]    (7.65,-2.3) .. controls (4.86,-0.97) and (2.31,-0.21) .. (0,0) .. controls (2.31,0.21) and (4.86,0.98) .. (7.65,2.3)   ;
		
		\draw    (49,135) -- (104,135) ;
		\draw [shift={(106,135)}, rotate = 180] [color={rgb, 255:red, 0; green, 0; blue, 0 }  ][line width=0.5]    (7.65,-2.3) .. controls (4.86,-0.97) and (2.31,-0.21) .. (0,0) .. controls (2.31,0.21) and (4.86,0.98) .. (7.65,2.3)   ;
		
		\draw  [color={rgb, 255:red, 0; green, 0; blue, 0 }  ,draw opacity=1 ][fill={rgb, 255:red, 0; green, 0; blue, 0 }  ,fill opacity=1 ] (33,205) .. controls (33,203.34) and (34.34,202) .. (36,202) .. controls (37.66,202) and (39,203.34) .. (39,205) .. controls (39,206.66) and (37.66,208) .. (36,208) .. controls (34.34,208) and (33,206.66) .. (33,205) -- cycle ;
		
		\draw    (36,197) -- (36,147) ;
		\draw [shift={(36,145)}, rotate = 450] [color={rgb, 255:red, 0; green, 0; blue, 0 }  ][line width=0.5]    (7.65,-2.3) .. controls (4.86,-0.97) and (2.31,-0.21) .. (0,0) .. controls (2.31,0.21) and (4.86,0.98) .. (7.65,2.3)   ;
		
		\draw  [color={rgb, 255:red, 0; green, 0; blue, 0 }  ,draw opacity=1 ][fill={rgb, 255:red, 0; green, 0; blue, 0 }  ,fill opacity=1 ] (113,205) .. controls (113,203.34) and (114.34,202) .. (116,202) .. controls (117.66,202) and (119,203.34) .. (119,205) .. controls (119,206.66) and (117.66,208) .. (116,208) .. controls (114.34,208) and (113,206.66) .. (113,205) -- cycle ;
		
		\draw    (116,143) -- (116,193) ;
		\draw [shift={(116,195)}, rotate = 270] [color={rgb, 255:red, 0; green, 0; blue, 0 }  ][line width=0.5]    (7.65,-2.3) .. controls (4.86,-0.97) and (2.31,-0.21) .. (0,0) .. controls (2.31,0.21) and (4.86,0.98) .. (7.65,2.3)   ;
		
		\draw  [draw opacity=0] (20.68,127.24) .. controls (18.37,124.09) and (17,120.2) .. (17,116) .. controls (17,105.51) and (25.51,97) .. (36,97) .. controls (46.49,97) and (55,105.51) .. (55,116) .. controls (55,121.25) and (52.87,126) .. (49.44,129.44) -- (36,116) -- cycle ; \draw   (20.68,127.24) .. controls (18.37,124.09) and (17,120.2) .. (17,116) .. controls (17,105.51) and (25.51,97) .. (36,97) .. controls (46.49,97) and (55,105.51) .. (55,116) .. controls (55,121.25) and (52.87,126) .. (49.44,129.44) ;
		
		\draw  [color={rgb, 255:red, 0; green, 0; blue, 0 }  ,draw opacity=1 ][fill={rgb, 255:red, 0; green, 0; blue, 0 }  ,fill opacity=1 ] (283,135) .. controls (283,133.34) and (284.34,132) .. (286,132) .. controls (287.66,132) and (289,133.34) .. (289,135) .. controls (289,136.66) and (287.66,138) .. (286,138) .. controls (284.34,138) and (283,136.66) .. (283,135) -- cycle ;
		
		\draw  [draw opacity=0] (276.13,130.84) .. controls (271.78,127.36) and (269,122) .. (269,116) .. controls (269,105.51) and (277.51,97) .. (288,97) .. controls (298.49,97) and (307,105.51) .. (307,116) .. controls (307,121.25) and (304.87,126) .. (301.44,129.44) -- (288,116) -- cycle ; \draw   (276.13,130.84) .. controls (271.78,127.36) and (269,122) .. (269,116) .. controls (269,105.51) and (277.51,97) .. (288,97) .. controls (298.49,97) and (307,105.51) .. (307,116) .. controls (307,121.25) and (304.87,126) .. (301.44,129.44) ;
		
		\draw    (304,126) -- (300.28,130.46) ;
		\draw [shift={(299,132)}, rotate = 309.81] [color={rgb, 255:red, 0; green, 0; blue, 0 }  ][line width=0.5]    (7.65,-2.3) .. controls (4.86,-0.97) and (2.31,-0.21) .. (0,0) .. controls (2.31,0.21) and (4.86,0.98) .. (7.65,2.3)   ;
		
		\draw  [color={rgb, 255:red, 0; green, 0; blue, 0 }  ,draw opacity=1 ][fill={rgb, 255:red, 0; green, 0; blue, 0 }  ,fill opacity=1 ] (203,135) .. controls (203,133.34) and (204.34,132) .. (206,132) .. controls (207.66,132) and (209,133.34) .. (209,135) .. controls (209,136.66) and (207.66,138) .. (206,138) .. controls (204.34,138) and (203,136.66) .. (203,135) -- cycle ;
		
		\draw    (190,126) -- (191.8,128.4) ;
		\draw [shift={(193,130)}, rotate = 233.13] [color={rgb, 255:red, 0; green, 0; blue, 0 }  ][line width=0.5]    (7.65,-2.3) .. controls (4.86,-0.97) and (2.31,-0.21) .. (0,0) .. controls (2.31,0.21) and (4.86,0.98) .. (7.65,2.3)   ;
		
		\draw    (219,135) -- (274,135) ;
		\draw [shift={(276,135)}, rotate = 180] [color={rgb, 255:red, 0; green, 0; blue, 0 }  ][line width=0.5]    (7.65,-2.3) .. controls (4.86,-0.97) and (2.31,-0.21) .. (0,0) .. controls (2.31,0.21) and (4.86,0.98) .. (7.65,2.3)   ;
		
		\draw  [color={rgb, 255:red, 0; green, 0; blue, 0 }  ,draw opacity=1 ][fill={rgb, 255:red, 0; green, 0; blue, 0 }  ,fill opacity=1 ] (203,205) .. controls (203,203.34) and (204.34,202) .. (206,202) .. controls (207.66,202) and (209,203.34) .. (209,205) .. controls (209,206.66) and (207.66,208) .. (206,208) .. controls (204.34,208) and (203,206.66) .. (203,205) -- cycle ;
		
		\draw    (206,197) -- (206,147) ;
		\draw [shift={(206,145)}, rotate = 450] [color={rgb, 255:red, 0; green, 0; blue, 0 }  ][line width=0.5]    (7.65,-2.3) .. controls (4.86,-0.97) and (2.31,-0.21) .. (0,0) .. controls (2.31,0.21) and (4.86,0.98) .. (7.65,2.3)   ;
		
		\draw  [draw opacity=0] (190.68,127.24) .. controls (188.37,124.09) and (187,120.2) .. (187,116) .. controls (187,105.51) and (195.51,97) .. (206,97) .. controls (216.49,97) and (225,105.51) .. (225,116) .. controls (225,121.25) and (222.87,126) .. (219.44,129.44) -- (206,116) -- cycle ; \draw   (190.68,127.24) .. controls (188.37,124.09) and (187,120.2) .. (187,116) .. controls (187,105.51) and (195.51,97) .. (206,97) .. controls (216.49,97) and (225,105.51) .. (225,116) .. controls (225,121.25) and (222.87,126) .. (219.44,129.44) ;
		
		\draw  [color={rgb, 255:red, 0; green, 0; blue, 0 }  ,draw opacity=1 ][fill={rgb, 255:red, 0; green, 0; blue, 0 }  ,fill opacity=1 ] (373,135) .. controls (373,133.34) and (374.34,132) .. (376,132) .. controls (377.66,132) and (379,133.34) .. (379,135) .. controls (379,136.66) and (377.66,138) .. (376,138) .. controls (374.34,138) and (373,136.66) .. (373,135) -- cycle ;
		
		\draw    (360,126) -- (361.8,128.4) ;
		\draw [shift={(363,130)}, rotate = 233.13] [color={rgb, 255:red, 0; green, 0; blue, 0 }  ][line width=0.5]    (7.65,-2.3) .. controls (4.86,-0.97) and (2.31,-0.21) .. (0,0) .. controls (2.31,0.21) and (4.86,0.98) .. (7.65,2.3)   ;
		
		\draw  [color={rgb, 255:red, 0; green, 0; blue, 0 }  ,draw opacity=1 ][fill={rgb, 255:red, 0; green, 0; blue, 0 }  ,fill opacity=1 ] (373,205) .. controls (373,203.34) and (374.34,202) .. (376,202) .. controls (377.66,202) and (379,203.34) .. (379,205) .. controls (379,206.66) and (377.66,208) .. (376,208) .. controls (374.34,208) and (373,206.66) .. (373,205) -- cycle ;
		
		\draw    (376,197) -- (376,147) ;
		\draw [shift={(376,145)}, rotate = 450] [color={rgb, 255:red, 0; green, 0; blue, 0 }  ][line width=0.5]    (7.65,-2.3) .. controls (4.86,-0.97) and (2.31,-0.21) .. (0,0) .. controls (2.31,0.21) and (4.86,0.98) .. (7.65,2.3)   ;
		
		\draw  [draw opacity=0] (360.68,127.24) .. controls (358.37,124.09) and (357,120.2) .. (357,116) .. controls (357,105.51) and (365.51,97) .. (376,97) .. controls (386.49,97) and (395,105.51) .. (395,116) .. controls (395,121.25) and (392.87,126) .. (389.44,129.44) -- (376,116) -- cycle ; \draw   (360.68,127.24) .. controls (358.37,124.09) and (357,120.2) .. (357,116) .. controls (357,105.51) and (365.51,97) .. (376,97) .. controls (386.49,97) and (395,105.51) .. (395,116) .. controls (395,121.25) and (392.87,126) .. (389.44,129.44) ;
		
		\draw (-27,165) node [anchor=north west][inner sep=0.5pt]   [align=left] {$ $};
		
		\draw (0,137) node [anchor=north west][inner sep=0.5pt]   [align=left] {$\displaystyle E:$};
		
		\draw (111,121) node [anchor=north west][inner sep=0.5pt]   [align=left] {$\displaystyle u$};
		
		\draw (101,202) node [anchor=north west][inner sep=0.5pt]   [align=left] {$\displaystyle v$};
		
		\draw (139,109) node [anchor=north west][inner sep=0.5pt]   [align=left] {$\displaystyle e$};
		
		\draw (118,159) node [anchor=north west][inner sep=0.5pt]   [align=left] {$\displaystyle f$};
		
		\draw (3,104) node [anchor=north west][inner sep=0.5pt]   [align=left] {$\displaystyle e'$};
		
		\draw (30,117) node [anchor=north west][inner sep=0.5pt]   [align=left] {$\displaystyle u'$};
		
		\draw (20,165) node [anchor=north west][inner sep=0.5pt]   [align=left] {$\displaystyle f'$};
		
		\draw (18,198) node [anchor=north west][inner sep=0.5pt]   [align=left] {$\displaystyle v'$};
		
		\draw (73,118) node [anchor=north west][inner sep=0.5pt]   [align=left] {$\displaystyle g$};
		
		\draw (167,137) node [anchor=north west][inner sep=0.5pt]   [align=left] {$\displaystyle F_{1} :$};
		
		\draw (281,121) node [anchor=north west][inner sep=0.5pt]   [align=left] {$\displaystyle u$};
		
		\draw (309,109) node [anchor=north west][inner sep=0.5pt]   [align=left] {$\displaystyle e$};
		
		\draw (173,104) node [anchor=north west][inner sep=0.5pt]   [align=left] {$\displaystyle e'$};
		
		\draw (200,117) node [anchor=north west][inner sep=0.5pt]   [align=left] {$\displaystyle u'$};
		
		\draw (190,165) node [anchor=north west][inner sep=0.5pt]   [align=left] {$\displaystyle f'$};
		
		\draw (188,198) node [anchor=north west][inner sep=0.5pt]   [align=left] {$\displaystyle v'$};
		
		\draw (243,118) node [anchor=north west][inner sep=0.5pt]   [align=left] {$\displaystyle g$};
		
		\draw (334,137) node [anchor=north west][inner sep=0.5pt]   [align=left] {$\displaystyle F_{2} :$};
		
		\draw (343,104) node [anchor=north west][inner sep=0.5pt]   [align=left] {$\displaystyle e'$};
		
		\draw (370,117) node [anchor=north west][inner sep=0.5pt]   [align=left] {$\displaystyle u'$};
		
		\draw (360,165) node [anchor=north west][inner sep=0.5pt]   [align=left] {$\displaystyle f'$};
		
		\draw (358,198) node [anchor=north west][inner sep=0.5pt]   [align=left] {$\displaystyle v'$};
		\end{tikzpicture}
	\end{center}
	
	In the graph $E$, the set $H_1=\{v\}$ is a hereditary saturated subset of $E^0$, so for the quotient graph $F_1=E\backslash H$, we have  $L_K(E)/I(H)\cong L_K(F_1)$. Because $E^0\backslash H_1=M(u)$,  by Theorem \ref{theorem_3.3}, for each irreducible polynomial $f(x)\in K[x,x^{-1}]$, the ideal $P=I(H_1, f(e))$ is a non-graded primitive ideal of type III of $L_K(E)$. If we set $H_2=\{u,v\}$, then $H_2$ is also a hereditary saturated subset of $E^0$ and we set $F_2=E\backslash H_2$. It is clear that 
	$$
	I(H_2)=I(u,v)=I(u,v, f(e))=I(H_1, u),
	$$
	so $P\subseteq I(H_2)$. Since  $L_K(E)/I(H_2)\cong L_K(F_2)$, which is clearly non-commutative, we deduce that $L_K(E)/P$ is non-commutative too.
\end{example}

For the convenience of further use, we restate the following result from \cite[(11.9)]{Bo_lam_2001}.

\begin{theorem}[\textbf{Structure  Theorem for left primitive rings}]\label{theorem_3.4}
	Let $R$ be a left primitive ring and $V$  a faithful simple left $R$-module. Let $k={\rm End}_R(V)$. Then, $k$ is a division ring and $R$ is isomorphic to a dense ring of linear transformations on $V_k$. Moreover:
	\begin{enumerate}[font=\normalfont]
		\item[(1)] If $R$ is left artinian, then $n:=\dim_k(V)$ is finite and $R\cong {\rm M}_n(k)$.
		\item[(2)] If $R$ is not left artinian, then $\dim_k(V)$ is infinite, and for any positive integer $n$, there exists a subring $\Lambda_n$ of $R$ which admits a ring homomorphism onto $ {\rm M}_n(k)$.
	\end{enumerate}
\end{theorem}
The following lemma follows from \cite[Lemma 2.8 (p.30) and Exercises 2.2 and 2.3 (p.31)]{Bo_weh_1973}.
\begin{lemma}\label{lemma_3.5}
	Let $K$ be a non-absolute field. Then ${\rm GL}_2(K)$ contains a non-cyclic free subgroup of which the generators are as follows. Consider the following matrices in ${\rm M}_2(K)$:
	\begin{align*}
		&	A             =    \begin{pmatrix} 1 & 0  \\ \alpha & 1\\ \end{pmatrix}={\rm I}_2+\alpha f_{21},\\
		&	B             =    \begin{pmatrix} 1 & \alpha \\ 0 & 1 \\\end{pmatrix}={\rm I}_2+\alpha f_{12},\\
		&	C             =   \begin{pmatrix} \beta& 0  \\ \alpha & \beta^{-1} \\ \end{pmatrix} = A+(\beta-1)f_{11}+(\beta^{-1}-1)f_{22},\\
		&	D             =    \begin{pmatrix} \beta & \alpha  \\ 0 & \beta^{-1} \\ \end{pmatrix}= B+(\beta-1)f_{11}+(\beta^{-1}-1)f_{22}.
	\end{align*}
	\begin{enumerate}[font=\normalfont]
		\item[(1)] If ${\rm char} K =0$ and $\alpha$ is equal to $2$ or transcendental over the prime subfield ${\mathbb{Q}}$, then $A$ and $B$ generate a non-cyclic free subgroup of ${\rm GL}_2(K)$.
		\item[(2)] If ${\rm char} K = t>0$ and $\alpha, \beta$ are algebraically independent over the prime subfield $\mathbb{F}_t$, then  $C$ and $D$ generate a non-cyclic free subgroup of ${\rm GL}_2(K)$.
	\end{enumerate}
\end{lemma}
We are now in the position to prove the following theorem which is the main result of this section.
\begin{theorem}\label{theorem_3.6}
	Let $E$ be a graph, and $K$ a non-absolute field. If $L_K(E)$ is non-commutative, then $L_K(E)^\times$ contains a non-cyclic free subgroup.
\end{theorem}
\begin{proof}
	In view of Proposition \ref{proposition_3.2}, there exists a primitive ideal $P$ of $L_K(E)$ such that $L_K(E)/P$ is non-commutative. Also, Proposition \ref{proposition_2.4} implies that there is a left simple $L_K(E)$-module $V$ such that  $P={\rm Ann}_{L_K(E)}(V)$ and that $k:={\rm End}_{L_{K}(E)}(V)$ is a field containing $K$. Let $\varphi: L_K(E)\to L_K(E)/P$ be the natural epimorphism. Then, via $\varphi$, $V$ became $L_K(E)/P$-module.  A routine checking shows that $V$  is a left faithful simple $L_K(E)/P$-module, and  ${\rm End}_{L_K(E)/P}(V)={\rm End}_{L_{K}(E)}(V)=k$. By  Theorem \ref{theorem_3.4}, there are two possible cases:
	
	\bigskip
	
	\textit{\textbf{Case 1}: $L_K(E)/P$ is artinian.} If we set $n=\dim_kV$, then $L_K(E)/P\cong {\rm M}_n(k)$. Since $L_K(E)/P$ is non-commutative, we must have $n>1$. Since $K\subseteq k$, and $K$ is non-absolute, it follows that $k$ is also non-absolute. Because ${\rm GL}_n(k)$ contains a copy of ${\rm GL}_2(k)$, by Lemma \ref{lemma_3.5}, we conclude that ${\rm GL}_n(k)$ also contains a non-cyclic free subgroups with the generators $x, y$. It can be shown that we may choose some $a\in \varphi^{-1}(x), b\in \varphi^{-1}(y)$ such that $a, b\in L_K(E)^\times$. (This fact will be shown in the next three sections where we give in detail such the elements $a$ and $b$ in term of the base graph $E$.) Therefore, we may apply Lemma \ref{lemma_3.1} to conclude that $ \left\langle a,b\right\rangle $ is a non-cyclic free subgroup of   $L_K(E)^\times$. 
	
	\bigskip
	
	\textit{\textbf{Case 2}: $L_K(E)/P$ is not artinian.} Again, according to Theorem \ref{theorem_3.4}, for each $n>1$, there exists a subring $\Lambda_n$ of $L_K(E)/P$ and a surjective ring homomorphism $f: \Lambda_n\to{\rm M}_n(k)$. Let $x$, $y$ be the generators of a free subgroup of ${\rm GL}_n(k)$ given in Lemma \ref{lemma_3.5}. Since $\varphi$ and $f$ are surjective, there exist $a, b\in L_K(E)^\times$ such that $\varphi(f(a))=x$ and $\varphi(f(b))=y$. We may choose $a$ and $b$ such that $a, b\in L_K(E)^\times$. (Again, this fact will be shown in the next three sections.) Finally, it follows from Lemma  \ref{lemma_3.1} that $ \left\langle a,b\right\rangle $ is a non-cyclic free subgroup of   $L_K(E)^\times$. 
\end{proof}

In the following, we shall examine the existence of non-cyclic free subgroups in the Cohn path algebra $C_K(E)$. As a consequence of the above theorem, we shall prove that non-commutative Cohn path algebra always contains a non-cyclic free subgroup. Before proceeding to this purpose, we need to record the following auxiliary lemma. 
\begin{lemma}\label{lemma_cohn-path}
	Let $E$ be a graph and $K$ a field. Then, the Cohn path algebra $C_K(E)$ is commutative if and only if so is the Leavitt path algebra $L_K(E)$. 
\end{lemma}
\begin{proof}
	For the proof of the ``only if" part, we assume that $L_K(E)$ is commutative. For each edge $\in E^1$, we have $e=s(e)er(e)=s(e)r(e)e$, which implies that $s(e)=r(e)$. This says that $E$ contains either isolated vertices or isolated graphs $R_n$ (see the picture in the end of Introduction). The commutativity of $L_K(E)$ assures that $n=1$. Therefore $E$  consists of only isolated vertices or isolated graphs $R_1$. This implies that $C_K(E)$ is also commutative. The converse follows immediately in a similar way.
\end{proof}

\begin{corollary}\label{corollary_3.8}
	Let $E$ be a graph, and $K$ a non-absolute field. If $C_K(E)$ is non-commutative, then $C_K(E)^\times$ contains a non-cyclic free subgroup.
\end{corollary}
\begin{proof}
	Let $I$ be the ideal of $C_K(E)$ generated by the set $$\left\lbrace v-\sum_{e\in s^{-1}(v)}ee^* \;|\; v\in{\rm Reg}(E)\right\rbrace .$$ It follows from \cite[Corollary 1.5.5]{Bo_abrams_2017} that $C_K(E)/I\cong L_K(E)$, which is non-commutative by Lemma \ref{lemma_cohn-path}. With reference to Theorem \ref{theorem_3.6}, we conclude that $L_K(E)^\times$ contains a non-cyclic free subgroup. According to Lemma \ref{lemma_3.1}, $C_K(E)^\times$ also contains a non-cyclic free subgroup.
\end{proof}
\begin{corollary}
	Let $E$ be a graph, and $K$ a non-absolute field. If the Leavitt path algebra $L_K(E)$ is non-commutative, then it contains a non-commutative free $K$-subalgebra.
\end{corollary}
\begin{proof}
	Let $x$ and $y$ be the generators of the non-cyclic free subgroup in $L_K(E)$. Then the $K$-subalgebra $K\left\langle x, y\right\rangle $ of $L_K(E)$ generated by $x$ and $y$ over $K$  is of course a free algebra.
\end{proof}
\begin{corollary}
	Let $E$ be a graph, and $K$ a non-absolute field. If $(L_K(E))^\times$ is a locally solvable group, then $L_K(E)$ is commutative.
\end{corollary}
\begin{proof}
	The local solvability of $L_K(E)^\times$ implies that it does not contain a non-cyclic free subgroup. It follows from Theorem \ref{theorem_3.6} that $L_K(E)$ is commutative. 
\end{proof}
In the literature, there are many works devoted to the study of the multiplicative group of a semiprime ring (see e.g. \cite{Pa_lanski_1971}, \cite{Pa_lanski_1981} and references therein). For instance, in \cite{Pa_lanski_1971} and \cite{Pa_lanski_1981}, C. Lanski have intensively examined  the algebraic structure of such a group, and more generally, he also studied the solvable normal subgroups of the multiplicative group of a semiprime ring. Now, let $R$ be a semiprime ring containing a non-central idempotent, and $S$ the additive subgroup generated by all elements having square zero. Then, as it was pointed out by Lanski in  \cite[Point~1, p.313]{Pa_lanski_1971}, $R$ contains a non-central Lie ideal $L$ such that $L\subseteq S$. Therefore,  Theorem~9 of \cite{Pa_lanski_1971} may be rephrased as follows. 
\begin{theorem}\label{theorem_3.9}
	Let $R$ be a simiprime ring containing a non-central idempotent, and $S$ the additive subgroup of $R$ generated by all elements having square zero. Then, $R$ contains a non-central Lie ideal $L$ such that $L\subseteq S$. Moreover, if $R$ is $6$-torsion-free $($that is, for each $x\in R$, $6x=0$ implies $x=0$$)$, then every solvable normal subgroup $G$ of $R^\times$ commutes elementwise with $L$. In particular, $G$ is contained in the center of $R$.
\end{theorem} 

As a consequence of this theorem, we have:

\begin{corollary}\label{corollary_3.10}
	Let $E$ be a graph and $K$ a field. Assume that $L_K(E)$ is a $6$-torsion-free Leavitt path algebra. If $G$ is a solvable normal subgroup of $L_K(E)^\times$, then $G$ is contained in the center of $L_K(E)$.
\end{corollary}
\begin{proof}
	We know that for a graph $E$ and a field $K$, the Leavitt path algebra $L_K(E)$ is always semiprime (see, e.g. \cite[Proposition 2.3.1]{Pa_abrams-pino_2005}). Clearly, we may assume that $L_K(E)$ is non-commutative. If $E$ contains a non-central vertex, then this vertex is a non-central idempotent of $L_K(E)$. Assume that all vertices of $E$ are central. For any $e\in E^1$, because $s(e)$ and $r(e)$ are central, we conclude that $s(e)=r(e)$ and, in consequence, $e$ is a loop. This implies that all edges in $E$ are loops. If $E$ contains only either isolated vertices or isolated loops, $L_K(E)$ is commutative, which is a contradiction. The consequence is that $E$ contains at least an isolated graph $R_n$ (see the picture in the end of Introduction), the rose with $n\geq 2$ petals. 
	
	Now, it is simple to check that each $e_ie_i^*$ is a non-central idempotent. In any case, there always exists a non-central idempotent in $L_K(E)$. Finally, the result follows from Theorem \ref{theorem_3.9}.
\end{proof}

In view of Theorem \ref{theorem_3.6} and Corollary \ref{corollary_3.10}, it is reasonable to pose the following conjecture.

\begin{conjecture}\label{conjecture_3.1}
	Let $E$ be a graph and $K$ a non-absolute field. Then, every non-central normal subgroup of $L_K(E)^\times$ contains a non-cyclic free subgroup.
\end{conjecture}

Although we know that the multiplicative group of a non-commutative Leavitt path algebra always contains a non-cyclic free subgroup, it is quite hard to describe the generators of this free subgroup. Now, we devote a major remaining part of this paper to examining this problem. First, we record an easier case:
\begin{proposition}
	Let $A_n$ be the oriented $n$-line graph with $n$ vertices and $n-1$ edges $($see the picture in the end of Introduction$)$, $n\geq 2$, and $K$ a non-absolute field.
	
	If we set
	$$
	a=1_{L_K(E)}+\alpha e_ie_{i+1}\dots e_j,\;\;\; b=1_{L_K(E)}+\alpha e_j^*\dots e_{i+1}^*e_i^*,
	$$ 
	$$
	c=a +(\beta-1) v_i + (\beta^{-1}-1) v_j,\;\;\; d=b+(\beta-1)v_i + (\beta^{-1}-1) v_j ,
	$$ 
	where $\alpha, \beta \in K$ and $1\leq i\lneqq j\leq n$, then the following statements hold:
	\begin{enumerate}[font=\normalfont]
		\item[(1)] If ${\rm char} K =0$ and $\alpha$ is equal to $2$ or transcendental over the prime subfield ${\mathbb{Q}}$, then $\left\langle \alpha,\beta\right\rangle $ is a non-cyclic free subgroup of $L_K(A_n)^\times$.
		\item[(2)] If ${\rm char} K = t>0$ and $\alpha,\beta$ are algebraically independent over the prime subfield $\mathbb{F}_t$, then  $\left\langle c,d\right\rangle $ is a non-cyclic free subgroup of $L_K(A_n)^\times$.
	\end{enumerate}
\end{proposition}
\begin{proof}
	It was shown in \cite[Proposition 1.3.5]{Bo_abrams_2017} that the map $$\varphi: L_K(A_n)\longrightarrow {\rm M}_n(K),$$
	which is defined by $\varphi(v_i)=f_{ii}$, $\varphi(e_i)=f_{i(i+1)}$ and $\varphi(e_i^*)=f_{(i+1)i}$, is an isomorphism. Since $f_{ij}=f_{i(i+1)}\dots f_{j-1,j}$, a straightforward checking shows that   
	$$ 
	\varphi(a)={\rm I}_n+xf_{ij},
	$$
	$$
	\varphi(b)={\rm I}_n+yf_{ji},
	$$
	$$ 
	\varphi(c)= {\rm I}_n+xf_{ij} + (y-1)f_{ii}+(y^{-1}-1)f_{jj},
	$$
	$$
	\varphi(d)={\rm I}_n+xf_{ji}+(y-1)f_{ii}+(y^{-1}-1)f_{jj}.
	$$
	It is easy to see that $\varphi(a)$, $\varphi(b)$, $\varphi(c)$, $\varphi(d)$ are the copies of the matrices $A$, $B$, $C$, $D$ in Lemma \ref{lemma_3.5}, respectively. Hence, the conclusions follow.
\end{proof}
In view of Proposition \ref{proposition_3.2}, a non-commutative Leavitt path algebra $L_K(E)$ always contains a primitive ideal $P$ for which $L_K(E)/P$ is non-commutative. Also, the proof of Theorem \ref{theorem_3.6} shows that we may construct a non-cyclic free subgroup in $L_K(E)$ as follows. First,  we construct a non-cyclic free subgroup $G$ in $L_K(E)/P$, then the inverse image $\varphi^{-1}(G)$ of $G$ under the natural homomorphism $\varphi: L_K(E)\to L_K(E)/P$ is a non-cyclic free subgroup in $L_K(E)$. In view of Theorem \ref{theorem_3.3}, there are exactly three types of primitive ideals, which we named as I, II and III, in $L_K(E)$. Therefore, the non-cyclic free subgroups in $L_K(E)$ may be determined by these types of primitive ideals. In the sequel, we shall construct free subgroups in $L_K(E)$ by this way. 
 
\section{Free subgroups determined by the primitive ideals of type I}
First, we need to determine a $K$-basis for $L_K(E)$.
\begin{lemma}[{\cite[Corollary 1.5.12]{Bo_abrams_2017}}]\label{lemma_4.1}
	Let $E$ be a graph and $K$ a field. 
	\begin{enumerate}[font=\normalfont]
		\item[(1)] The set of monomials 
		$$
		\mathscr{A}=\{\lambda\nu^*:r(\lambda)=r(\nu)\}
		$$ 
		is a $K$-basis of the Cohn path algebra $C_K(E)$.
		\item[(2)] For each $v\in{\rm Reg}(E) $, let $\{e_1^v,\dots,e_{n_v}^v\}$ be an enumeration of the elements of $s^{-1}(v)$. Then, a $K$-basis of the Leavitt path algebra $L_K(E)$ is given by the family 
		$$\mathscr{B}=\mathscr{A}\backslash\{\lambda e_{n_v}^v(e_{n_v}^v)^*\nu^*: r(\lambda)=r(\nu)=v\in {\rm Reg}(E)\}.$$
	\end{enumerate}
\end{lemma}

The following lemma shows that a non-cyclic free subgroup in $L_K(E)$ may be determined by a sink in $E$.

\begin{lemma}\label{lemma_4.2}
	Let $E$ be a graph and $K$ a non-absolute field. Assume that there exist a sink $w$ and an edge $f$ such that $r(f)=w$. Then, the following assertions hold.
	\begin{enumerate}[font=\normalfont]
		\item[(i)] If ${\rm char}\; K =0$ and $\alpha$ is either equal to $2$ or transcendental over the prime subfield ${\mathbb{Q}}$, then 
		$$
		a:=1_{L_K(E)}+\alpha f^*,\hspace*{0,5cm} b:=1_{L_K(E)}+\alpha f
		$$
		generate a non-cyclic free subgroup in $L_K(E)$.
		\item[(ii)] If ${\rm char}\;K=t>0$ and $\alpha, \beta$ are algebraically independent over the prime subfield $\mathbb{F}_t$, then  
		$$c:=a +(\beta-1) s(f) + (\beta^{-1}-1) r(f),$$ $$d:=b+(\beta-1) s(f) + (\beta^{-1}-1) r(f) $$ generate a non-cyclic free subgroup in $L_K(E)$.
	\end{enumerate}
\end{lemma}
\begin{proof}
	If we set $V=Kw\oplus Kf$, then this is a non-zero $K$-space with a $K$-basis $\{w, f\}$ and $\dim_KV=2$. Put
	$$
	R=\{r\in L_K(E): rV\subseteq V\} \text{ and } I=\{r\in L_K(E): r V=0\}, 
	$$ 
	then $R$ is a subring of $L_K(E)$ and $I$ is an ideal of $R$, so $R/I$ can be viewed as a $K$-space. Now, we are going to seek a $K$-basis for $R/I$. Take $\gamma\lambda^*$ with $r(\gamma)= r(\lambda)$, and $v = a_1w+a_2f\in V$,  where $a_1, a_2\in K$. We first make the following computations
	$$
	a_1\gamma\lambda^* w =   \begin{cases} 
		a_1\gamma      & \text{ if }  \lambda=w,\\
		0      & \text{ otherwise}.
	\end{cases}
	$$
	
	$$
	a_2\gamma\lambda^* f =    \begin{cases} 
		a_2\gamma f      & \text{ if }  \lambda=s(f),\\
		a_2\gamma      & \text{ if }  \lambda=f,\\
		0      & \text{ otherwise}.
	\end{cases}
	$$
	It follows that 
	$$
	\gamma\lambda^* v = a_1\gamma\lambda^* w+a_2\gamma\lambda^*f 
	=  \begin{cases} 
		a_1\gamma      & \text{ if }  \lambda=w,\\
		a_2\gamma f      & \text{ if }  \lambda=s(f),\\
		a_2\gamma      & \text{ if }  \lambda=f,\\
		0      & \text{ otherwise}.
	\end{cases}\eqno(1)
	$$
	Now, take $r\in L_K(E)$ and write 
	$$
	r=k_1\gamma_1\lambda_1^*+k_2\gamma_1\lambda_2^*+\dots+k_n\gamma_n\lambda_n^*,
	$$
	where $k_i\in K$ and $ \gamma_i, \lambda_i \in{\rm Path}(E) $ with $r(\gamma_i)=r(\lambda_i)$ for all $1\leq i \leq n$. 
	In view of (1), $r \in R$ if and only if either $r v=0 $, or else $\lambda_i\in \{s(f), w, f\}$ and $\gamma_i \in \{s(f), w, f\}$. Equivalently, $r \in R$ if and only if either $r v=0 $ or $\gamma_i\lambda_i^*\in \{s(f), w, f, f^*, ff^*\}$. Since $ ff^*v=ff^*(a_1w+a_2f)=a_2f $, we conclude that $ff^*+I=a_2f+I$. This means that 	
	$$
	R/I=\left\lbrace k_1s(f)+k_2w+k_3f+k_4f^*+I \;|\; k_1, k_2, k_3, k_4\in K\right\rbrace.
	$$
	From this, it follows that $s(f)+I, r(f)+I, f+I, f^*+I$ spans $R/I$ as a $K$-space. Moreover, in view of Lemma \ref{lemma_4.1}, one can easily check  that these elements are linearly independent over $K$. Therefore, the set $\{s(f)+I, r(f)+I, f+I, f^*+I\}$ is a $K$-basis for $R/I$. Accordingly, the mapping $\varphi: R \to {\rm M}_2(K)$ given by 
	$$
	\begin{aligned}
		\varphi(s(f))=   f_{11}, & \hspace{1cm}
		\varphi(r(f))= f_{22},\\
		\varphi(f)= f_{12}, &  \hspace{1cm}
		\varphi(f^*)= f_{21}\\
	\end{aligned}
	$$
	is a surjective ring homomorphism with $\ker \varphi =I$. Consider the following matrices in ${\rm M}_2(K)$:
	\begin{align*}
		&	A             =    {\rm I}_2+\alpha f_{21},\\
		&	B             =    {\rm I}_2+\alpha f_{12},\\
		&	C             =    A+(\beta-1)f_{11}+(\beta^{-1}-1)f_{22},\\
		&	D             =    B+(\beta-1)f_{11}+(\beta^{-1}-1)f_{22}.
	\end{align*}
	It is easy to see that	
	\begin{align*}
		&	\varphi(a)             =    \varphi\left( 1_{L_K(E)}+\alpha f^*\right)  = A,\\
		&	\varphi(b)             =    \varphi\left( 1_{L_K(E)}+\alpha f\right)      = B,\\
		&	\varphi(c) 	            = 	\varphi\left( a +(\beta-1) s(f) + (\beta^{-1}-1) r(f)\right) = C,\\
		&	\varphi(d) 	            = 	\varphi\left( b+(\beta-1) s(f) + (\beta^{-1}-1) r(f)\right)= D.
	\end{align*}
	To prove (i), we assume that ${\rm char} K =0$ and $\alpha$ is either equal to $2$ or transcendental over the prime subfield ${\mathbb{Q}}$. According to Lemma \ref{lemma_3.5}, $A$ and $B$ generate a non-cyclic free subgroup of ${\rm GL}_2(K)$. A simple calculation shows that $a$, $b$, $c$, $d$ are invertible in $L_K(E)$ with 
	\begin{align*}
		&	a^{-1}     =   1_{L_K(E)}-\alpha f^*,\\
		&	b^{-1}     =   1_{L_K(E)}-\alpha f,\\
		&	c^{-1}     =   a^{-1} + (\beta^{-1}-1) s(f) + (\beta-1) r(f),\\
		&	d^{-1}     =   b^{-1} + (\beta^{-1}-1) s(f) + (\beta-1) r(f).
	\end{align*}
Since $\varphi: R\to {\rm M}_2(K)$ is surjective, we may apply Lemma \ref{lemma_3.1} to conclude that $\left\langle a, b\right\rangle $ is a  non-cyclic free subgroup of $R^\times \subseteq L_K(E)^\times$, proving (i). Finally, the proof of (ii) follows similarly. 
\end{proof}

\begin{lemma}\label{lemma_4.3}
	Let $E$ be a graph and $K$ a field. Let $(H,S)$ be an admissible pair in the graph $E$, $F=E\backslash(H,S)$ a quotient graph,   and $\varphi: L_K(E)\to L_K(E\backslash(H,S))$ the epimorphism given in Lemma~ \ref{lemma_1.4}. Assume that $f\in E^1$ with $w=r(f)\in B_H\backslash S$. Then we have 
	$$\varphi(s(f))=s(f'), \;\; \varphi(w^H)=w', \;\;  \varphi(fw^H)=f',  \;\; \varphi(w^Hf^*)=f'^*,$$
	where $w'$ and $f'$ are the corresponding vertex and edge in $F$ of $w$ and $f$ respectively given in Definition \ref{definition_1.2}.
\end{lemma}
\begin{proof}
	The Definition \ref{definition_1.2} confirms that $s(f')=s(f)$, so by Lemma \ref{lemma_1.4}, we have  $\varphi(s(f))=s(f')$. For the proof of the remaining relation, we recall that $$w^H=w-\sum_{s(e)=w,\;r(e)\not\in H}ee^*.$$
	For a sake of simplicity, we set $A=\sum_{s(e)=w,\;r(e)\not\in H}ee^*$. First, we claim that $\varphi(A)=w\in F$. Let $\{e_1, e_2,\dots,e_n\}$ be an enumeration of the elements of the set $\{e\in E^1| s(e)=w,r(e)\not\in H\}$. (Note that since $w$ is a breaking vertex of $H$, this set is non-empty and finite). If $r(e_i)\not\in B_H\backslash S$ for all $1\leq i\leq n$, then, by Lemma \ref{lemma_1.4}, we have $\varphi(e_i)=e_i$ and $\varphi(e_i^*)=e_i^*$ for all $i$. In this case, the set of all edges in $F$ having the common source $w$ is $\{e_1, e_2,\dots,e_n\}$. Therefore, the relation (CK2) in $F$ implies that $\varphi(A)=\sum_{i=1}^{n}e_ie_i^*=w\in F$. Now, assume that some of the $r(e_i)$'s belong to $B_H\backslash S$ and the remaining of them do not. Without lose of the generality, we may assume that $r(e_1),\dots, r(e_m) \in B_H\backslash S$ and $r(e_{m+1}),\dots, r(e_n)\not\in B_H\backslash S$, where $1\leq m\le n$. It follows from Lemma  \ref{lemma_1.4} that $\varphi(e_i)=e_i+e_i'$, $\varphi(e_i^*)=e_i^*+e_i'^*$ for $1\leq i\leq m$ and $\varphi(e_i)=e_i$, $\varphi(e_i^*)=e_i^*$ for $m+1\leq i\leq n$. Consequently, we have 
	$$
	\varphi(A)=\sum_{i=1}^{m}(e_i+e_i')(e_i^*+e_i'^*)+\sum_{i=m+1}^{n}ee_i^*.
	$$
	For each $1\leq i\leq m$, in $F$, we have the following relation $$(e_i+e_i')(e_i^*+e_i'^*)=e_ie_i^*+e_ie_i'^*+e_i'e_i^*+e_i'e_i'^*.$$
	Because $r(e')$ is a sink in $F$, we have $e_ie_i'^*=e_i'e_i^*=0$, from which it follows that $(e_i+e_i')(e_i^*+e_i'^*)=e_ie_i^*+e_i'e_i'^*$. This implies that 
	$$ \varphi(A)=\sum_{i=1}^{m}e_ie_i^*+\sum_{i=1}^{m}e_i'e_i'^*+\sum_{i=m+1}^{n}ee_i^*=\sum_{i=1}^{n}e_ie_i^*+\sum_{i=1}^{m}e_i'e_i'^*.
	$$
	Observe that the set of all  edges in $F$ having the common source $w$ is $$\{e_1, e_2,\dots,e_n\} \cup \{e_1',\dots,e_m'\}.$$ Again, the relation (CK2) in $F$ implies that $ \varphi(A)=w\in F$, and the claim is shown. Using the fact that $\varphi(A)=w$, we now easily make the following computations:
	\begin{align*}
		&	\varphi(w^H)\hspace*{0.3cm}			= 	\varphi(w)-\varphi(A) = w+w'-w=w'.\\
		&	\varphi(fw^H)\hspace*{0.1cm}			= 	\varphi(f)\varphi(w)-\varphi(f)\varphi(A)\\&\hspace*{1.4cm} = (f+f')(w+w')-(f+f')w=f+f'-f=f'.\\	
		&	\varphi(w^Hf^*)		  =   \varphi(w)\varphi(f^*)-\varphi(A)\varphi(f^*) = (w+w')(f^*+f'^*)-w(f^*+f'^*)\\&\hspace*{1.4cm}=f^*+f'^*-f^*=f'^*.	
	\end{align*}
	Therefore, the proof is now complete.
\end{proof}

The following theorem is the main result of this section.
\begin{theorem}
	Let $E$ be a graph and $K$ a non-absolute field. Assume that $L_K(E)$ contains a  primitive ideal $P$  of type I given in Theorem \ref{theorem_3.3} such that $L_K(E)/P$ is non-commutative. Let $w$ be the breaking vertex determining $P$. Then, for any $f\in E^1$ with $r(f)=w$, the following assertions hold:
	\begin{enumerate}[font=\normalfont]
		\item[(i)] If ${\rm char}\; K =0$ and $\alpha\in K$ is either equal to $2$ or transcendental over the prime subfield ${\mathbb{Q}}$, then the elements
		$$
		a:=1_{L_K(E)}+\alpha w^Hf^* \text{ and } b:=1_{L_K(E)}+\alpha fw^H
		$$ 
		generate a non-cyclic free subgroup in $L_K(E)$.
		\item[(ii)] If ${\rm char}\;K=t>0$ and $\alpha, \beta\in K$ are algebraically independent over the prime subfield $\mathbb{F}_t$, then  
		$$c:=a +(\beta-1) s(f) + (\beta^{-1}-1) w^H,$$ 
		$$d:=b+(\beta-1) s(f) + (\beta^{-1}-1) w^H$$
		generate a non-cyclic free subgroup in $L_K(E)$.
	\end{enumerate}
\end{theorem}
\begin{proof}
	Take a vertex $w\in B_H$ as in Theorem \ref{theorem_3.3}, namely $w\in B_H$ is  a vertex such that $P=I(H, B_H\backslash\{w\})$ and $M(w)=E^0\backslash H$. Let $F=E\backslash (H,B_H\backslash\{w\})$.  By Lemma \ref{lemma_1.4},  there is an epimorphism $\varphi: L_K(E) \to L_K(F)$ with $\ker\varphi = P$. The quotient graph $F$ is described as follows:
	$$ 
	\begin{aligned}
		& F^0= (E^0\backslash H) \cup \{w'\};\\
		& F^1=\{e\in E^1: r(e)\not\in H\}\cup \{f':  r(f)=w\},
	\end{aligned}
	$$
	and $r$, $s$ are extended to $F^0$ by setting $s(f')=s(f)$ and $r(f')=w'$ for all $f\in E^1$ with $r(f)=w$. Put $A=\{f\in E^1| r(f)=w\}$. If $A=\varnothing$, then $F^0$ would be the single point $w'$ which means that $L_K(F)$ is commutative, a contradiction. This guarantees that $A\ne\varnothing$. Now, take $f\in A$ and set $f'=\varphi(f)$. Then $f'$ is an edge in $F^1$ satisfying $r(f')=w'$. Because $w'$ is a sink in $F$, we conclude that $s(f')\ne r(f')={w'}$.  To prove (i), assume that ${\rm char} F=0$, $\alpha$ is either equal to $2$ or transcendental over the prime subfield ${\mathbb{Q}}$, and put 	
	$$
		a'   = 1_{L_K(F)}+\alpha f'^* , \;\;\; b'= 1_{L_K(F)}+\alpha f'.
	$$

	In view of Lemma \ref{lemma_4.3},  $\varphi(s(f))=s(f')$, $\varphi(w^H)=w'$, $\varphi(fw^H)=f'$ and $\varphi(w^Hf^*)=f'^*$.
	This implies that 
	\begin{align*}
		&	\varphi(a)   =    \varphi\left( 1_{L_K(E)}+\alpha w^Hf^* \right)  = 1_{L_K(F)}+\alpha f'^* = a',\\
		&	\varphi(b)   =    \varphi\left( 1_{L_K(E)}+\alpha fw^H \right) = 1_{L_K(F)}+\alpha f' = b'.
	\end{align*}

	It follows from Lemma  \ref{lemma_4.2}(i) that $\left\langle a',b'\right\rangle $ is a non-cyclic free subgroup of $L_K(F)^\times$. Moreover, the elements $a$ and $b$ are invertible in $L_K(E)$ with the inverses 
	\begin{align*}
		&	a^{-1}     =   1_{L_K(E)}-\alpha w^Hf^*,\\
		&	b^{-1}     =   1_{L_K(E)}-\alpha fw^H.
	\end{align*}
	Therefore, Lemma \ref{lemma_3.1} implies that $\left\langle a,b\right\rangle $ is a non-cyclic free subgroup of $L_K(E)^\times$, proving (i). 

	To prove (ii), assume that ${\rm char}\;K=t>0$, $\alpha, \beta\in K$ are algebraically independent over the prime subfield $\mathbb{F}_t$, and put
	$$
		c'= a' + (\beta-1) s(f')+ (\beta^{-1}-1) w' \text{ and } d'= b' + (\beta-1) s(f')+ (\beta^{-1}-1) w'.
	$$
	We have 
	\begin{align*}
		&	\varphi(c) 	 = 	\varphi\left( a +(\beta-1) s(f) + (\beta^{-1}-1) w^H \right) \\
		& \hspace*{0.7cm} = a' + (\beta-1) s(f')+ (\beta^{-1}-1) w' = c',\\
		&	\varphi(d) 	 = 	\varphi\left( b++(\beta-1) s(f) + (\beta^{-1}-1) w^H \right) \\
		&\hspace*{0.7cm} = b' + (\beta-1) s(f')+ (\beta^{-1}-1) w' = d'.
	\end{align*}
	The two elements $c$ and $d$ are invertible in $L_K(E)$ with 
	\begin{align*}
		&	c^{-1}     =   a^{-1} + (\beta^{-1}-1) s(f) + (\beta-1) w^H\\
		&	d^{-1}     =   b^{-1} + (\beta^{-1}-1) s(f) + (\beta-1) w^H.
	\end{align*}
	According to Lemma  \ref{lemma_4.2}(ii), we conclude that $\left\langle c',d'\right\rangle $ is a non-cyclic free subgroup of $L_K(F)^\times$. Thus, Lemma \ref{lemma_3.1} implies that $\left\langle c,d\right\rangle $ is a non-cyclic free subgroup of $L_K(E)^\times$. The proof is now finished.
\end{proof}

In the following example, we demontrate a certain situation in which we give in detail a description of the gernerators of the non-cyclic free subgroup determined by a primitive ideal of type I.
\begin{example}
	Let consider again the graph $E$ given in Example \ref{example_3}. 
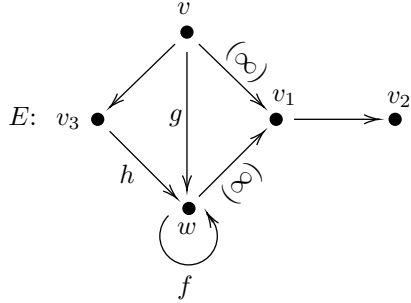
	\begin{center}
		
		\tikzset{every picture/.style={line width=0.5pt}} 
		
		\begin{tikzpicture}[x=0.75pt,y=0.75pt,yscale=-1,xscale=1]
		
		\draw  [color={rgb, 255:red, 0; green, 0; blue, 0 }  ,draw opacity=1 ][fill={rgb, 255:red, 0; green, 0; blue, 0 }  ,fill opacity=1 ] (133,135) .. controls (133,133.34) and (134.34,132) .. (136,132) .. controls (137.66,132) and (139,133.34) .. (139,135) .. controls (139,136.66) and (137.66,138) .. (136,138) .. controls (134.34,138) and (133,136.66) .. (133,135) -- cycle ;
		
		\draw  [draw opacity=0] (144.92,137.71) .. controls (148.21,140.36) and (150.31,144.43) .. (150.31,148.99) .. controls (150.31,156.99) and (143.82,163.48) .. (135.82,163.48) .. controls (127.81,163.48) and (121.32,156.99) .. (121.32,148.99) .. controls (121.32,144.9) and (123.01,141.22) .. (125.72,138.58) -- (135.82,148.99) -- cycle ; \draw   (144.92,137.71) .. controls (148.21,140.36) and (150.31,144.43) .. (150.31,148.99) .. controls (150.31,156.99) and (143.82,163.48) .. (135.82,163.48) .. controls (127.81,163.48) and (121.32,156.99) .. (121.32,148.99) .. controls (121.32,144.9) and (123.01,141.22) .. (125.72,138.58) ;
		
		\draw    (148,142) -- (146.03,138.71) ;
		\draw [shift={(145,137)}, rotate = 419.03999999999996] [color={rgb, 255:red, 0; green, 0; blue, 0 }  ][line width=0.5]    (7.65,-2.3) .. controls (4.86,-0.97) and (2.31,-0.21) .. (0,0) .. controls (2.31,0.21) and (4.86,0.98) .. (7.65,2.3)   ;
		
		\draw  [color={rgb, 255:red, 0; green, 0; blue, 0 }  ,draw opacity=1 ][fill={rgb, 255:red, 0; green, 0; blue, 0 }  ,fill opacity=1 ] (132,46) .. controls (132,44.34) and (133.34,43) .. (135,43) .. controls (136.66,43) and (138,44.34) .. (138,46) .. controls (138,47.66) and (136.66,49) .. (135,49) .. controls (133.34,49) and (132,47.66) .. (132,46) -- cycle ;
		
		\draw  [color={rgb, 255:red, 0; green, 0; blue, 0 }  ,draw opacity=1 ][fill={rgb, 255:red, 0; green, 0; blue, 0 }  ,fill opacity=1 ] (177,90) .. controls (177,88.34) and (178.34,87) .. (180,87) .. controls (181.66,87) and (183,88.34) .. (183,90) .. controls (183,91.66) and (181.66,93) .. (180,93) .. controls (178.34,93) and (177,91.66) .. (177,90) -- cycle ;
		
		\draw  [color={rgb, 255:red, 0; green, 0; blue, 0 }  ,draw opacity=1 ][fill={rgb, 255:red, 0; green, 0; blue, 0 }  ,fill opacity=1 ] (87,90) .. controls (87,88.34) and (88.34,87) .. (90,87) .. controls (91.66,87) and (93,88.34) .. (93,90) .. controls (93,91.66) and (91.66,93) .. (90,93) .. controls (88.34,93) and (87,91.66) .. (87,90) -- cycle ;
		
		\draw    (129,52) -- (97.44,82.61) ;
		\draw [shift={(96,84)}, rotate = 315.88] [color={rgb, 255:red, 0; green, 0; blue, 0 }  ][line width=0.5]    (7.65,-2.3) .. controls (4.86,-0.97) and (2.31,-0.21) .. (0,0) .. controls (2.31,0.21) and (4.86,0.98) .. (7.65,2.3)   ;
		
		\draw    (141,52) -- (172.56,82.61) ;
		\draw [shift={(174,84)}, rotate = 224.12] [color={rgb, 255:red, 0; green, 0; blue, 0 }  ][line width=0.5]    (7.65,-2.3) .. controls (4.86,-0.97) and (2.31,-0.21) .. (0,0) .. controls (2.31,0.21) and (4.86,0.98) .. (7.65,2.3)   ;
		
		\draw    (96,96) -- (128.59,128.59) ;
		\draw [shift={(130,130)}, rotate = 225] [color={rgb, 255:red, 0; green, 0; blue, 0 }  ][line width=0.5]    (7.65,-2.3) .. controls (4.86,-0.97) and (2.31,-0.21) .. (0,0) .. controls (2.31,0.21) and (4.86,0.98) .. (7.65,2.3)   ;
		
		\draw    (141,129) -- (172.59,97.41) ;
		\draw [shift={(174,96)}, rotate = 495] [color={rgb, 255:red, 0; green, 0; blue, 0 }  ][line width=0.5]    (7.65,-2.3) .. controls (4.86,-0.97) and (2.31,-0.21) .. (0,0) .. controls (2.31,0.21) and (4.86,0.98) .. (7.65,2.3)   ;
		
		\draw    (135,56) -- (135,124) ;
		\draw [shift={(135,126)}, rotate = 270] [color={rgb, 255:red, 0; green, 0; blue, 0 }  ][line width=0.5]    (7.65,-2.3) .. controls (4.86,-0.97) and (2.31,-0.21) .. (0,0) .. controls (2.31,0.21) and (4.86,0.98) .. (7.65,2.3)   ;
		
		\draw    (189,90) -- (230,90) ;
		\draw [shift={(232,90)}, rotate = 180] [color={rgb, 255:red, 0; green, 0; blue, 0 }  ][line width=0.5]    (7.65,-2.3) .. controls (4.86,-0.97) and (2.31,-0.21) .. (0,0) .. controls (2.31,0.21) and (4.86,0.98) .. (7.65,2.3)   ;
		
		\draw  [color={rgb, 255:red, 0; green, 0; blue, 0 }  ,draw opacity=1 ][fill={rgb, 255:red, 0; green, 0; blue, 0 }  ,fill opacity=1 ] (237,90) .. controls (237,88.34) and (238.34,87) .. (240,87) .. controls (241.66,87) and (243,88.34) .. (243,90) .. controls (243,91.66) and (241.66,93) .. (240,93) .. controls (238.34,93) and (237,91.66) .. (237,90) -- cycle ;
		
		\draw (-27,165) node [anchor=north west][inner sep=0.5pt]   [align=left] {$ $};
		
		\draw (44,83) node [anchor=north west][inner sep=0.5pt]   [align=left] {$\displaystyle E$:};
		
		\draw (129,141) node [anchor=north west][inner sep=0.5pt]   [align=left] {$\displaystyle w$};
		
		\draw (129,168) node [anchor=north west][inner sep=0.5pt]   [align=left] {$\displaystyle f$};
		
		\draw (125,84) node [anchor=north west][inner sep=0.5pt]   [align=left] {$\displaystyle g$};
		
		\draw (100,112) node [anchor=north west][inner sep=0.5pt]   [align=left] {$\displaystyle h$};
		
		\draw (161.46,46.74) node [anchor=north west][inner sep=0.5pt]  [rotate=-45] [align=left] {$\displaystyle ( \infty $)};
		
		\draw (177.26,118.46) node [anchor=north west][inner sep=0.5pt]  [rotate=-135] [align=left] {$\displaystyle ( \infty $)};
		
		\draw (68,86) node [anchor=north west][inner sep=0.5pt]   [align=left] {$\displaystyle v_{3}$};
		
		\draw (177,74) node [anchor=north west][inner sep=0.5pt]   [align=left] {$\displaystyle v_{1}$};
		
		\draw (235,75) node [anchor=north west][inner sep=0.5pt]   [align=left] {$\displaystyle v_{2}$};
		
		\draw (129,30) node [anchor=north west][inner sep=0.5pt]   [align=left] {$\displaystyle v$};
		
		\end{tikzpicture}
	\end{center}
	If $H=\{v_1,v_2\}$, then $B_H=\{v, w\}$. As we have point out in Example \ref{example_3}, $P:=I(H, B_H\backslash \{w\})$ is a primitive graded ideal of type I  of $L_K(E)$. It is clear that $w^H=w-ff^*$. Therefore, if $K$ is a field of characteristic $0$, for instance,  and $\alpha$ is either equal to $2$ or transcendental over the prime subfield ${\mathbb{Q}}$, then the following elements will generate a non-cyclic free subgroup of $L_K(E)^\times$:
	$$
	a=1_{L_K(E)}+\alpha(w-ff^*)f^* \text{ and } b=1_{L_K(E)}+\alpha f(w-ff^*).
	$$	
\end{example}

\section{Free subgroups determined by the primitive ideals of type II}
In this section, we study non-cyclic free subgroups in $L_K(E)$ determined by primitive ideals of type II. Recall that a primitive ideal $P$ of $L_K(E)$ is said to be of type II if $P$ is a graded ideal of the form $I(H,B_H)$, and $E\backslash(H,B_H)$ is downward directed, satisfying the Condition (L) and the Countable Separation Property (see Theorem \ref{theorem_3.3}). Before proceed further, we record an easy lemma.
\begin{lemma}\label{lemma_5.1}
	Let $E$ be a graph and $K$ a field. If $E$ contains two distinct tail-equivalent infinite paths $p'$ and $q'$, then there exist two infinite paths of the form $p$ and $fp$ which are tail-equivalent to $p'$ and $q'$, where $f\in E^1$ such that $r(f)=s(p)$.
\end{lemma}

\begin{proof}
	Since $p'$ and $q'$ are tail-equivalent infinite paths in $E$, there exist $m$, $n\geq 0$ which are minimal such that $\tau_{> m}(p')=\tau_{> n}(q')$. Write $p'=s(p')e_1\dots e_m\tau_{> m}(p')$ and $q'=s(q')f_1\dots f_n\tau_{> n}(q')$, where $e_i, f_j\in E^1$. Since $p'\ne q'$, we have either $m\geq 1$ or $n\geq 1$. Without loss of generality, we may assume that $n\geq 1$. Clearly, the paths $p=\tau_{> m}(p')=\tau_{> n}(q')$ and $q=f_np$ are infinite of desired form.
\end{proof}

The following lemma may be considered as the way to determine a non-cyclic free subgroup in $L_K(E)$ by an infinite path.
\begin{lemma}\label{lemma_5.3}
	Let $E$ be a graph and $K$ be a non-absolute field. Assume that $E$ contains an infinite path $p$ and an edge $f$ such that $s(f)\ne r(f)=s(p)$. Then, the following assertions hold.
	\begin{enumerate}[font=\normalfont]
		\item[(i)] If ${\rm char}\; K =0$ and $\alpha$ is either equal to $2$ or transcendental over the prime subfield ${\mathbb{Q}}$, then 
		$$
		a:=1_{L_K(E)}+\alpha f^*,\hspace*{0,5cm} b:=1_{L_K(E)}+\alpha f
		$$
		generate a non-cyclic free subgroup in $L_K(E)$.
		\item[(ii)] If ${\rm char}\;K=t>0$ and $\alpha, \beta$ are algebraically independent over the prime subfield $\mathbb{F}_t$, then  
		$$c:=a +(\beta-1) s(f) + (\beta^{-1}-1) r(f),$$ $$d:=b+(\beta-1) s(f) + (\beta^{-1}-1) r(f) $$ generate a non-cyclic free subgroup in $L_K(E)$.
	\end{enumerate}
\end{lemma}
\begin{proof}
		Set $V=Kp\oplus Kq$, where $q=fp$. Then, $V$ is a $K$-space with a basis $\{p, q\}$, so $\dim_KV=2$. Put
	$$
	R=\{r\in L_K(E): rV\subseteq V\} \text{ and } I=\{r\in L_K(E): r V=0\}. 
	$$ 
	Then, $R$ is a subring of $L_K(E)$ and $I$ is an ideal of $R$. Therefore, we can view $R/I$ as a $K$-space in an obvious way. Next, we are going to find a $K$-basis for $R/I$. To do this, take $\gamma\lambda^*$, where $\gamma, \lambda\in{\rm Path}(E)$ with $r(\gamma)=r(\lambda)$, and $v = a_1p+a_2q\in V$, where $a_1, a_2\in K$. We have
	$$
	a_1\gamma\lambda^* p = \begin{cases} 
		a_1\gamma\eta      & \text{ if }  p=\lambda\eta,\\
		0      & \text{ otherwise}.
	\end{cases}
	$$
	$$
	a_2\gamma\lambda^* q = \begin{cases} 
		a_2\gamma\beta   &  \text{ if }  q=\lambda\beta,\\
		0      & \text{ otherwise}.
	\end{cases}
	$$
	Since $s(p)\ne s(q)$, the cases $p=\lambda\eta$ and $q=\lambda\beta$ cannot happen at the same time. This implies that
	$$
	\gamma\lambda^*v= \gamma\lambda^* (a_1p+a_2q)=    \begin{cases} 
		a_1\gamma\eta      & \text{ if }  p=\lambda\eta,\\
		a_2\gamma\beta   &  \text{ if }  q=\lambda\beta,\\
		0      & \text{ otherwise}.
	\end{cases}\eqno(2)
	$$
	Therefore, $\gamma\lambda^*v\in V$ if and only if either $\gamma\lambda^*v=0$ or one of the following cases occurs:
	
	\bigskip 
	
	\textit{\textbf{Case 1}: $ p=\lambda\eta$ and $\gamma\eta\in V$}. But $\gamma\eta\in V$ implies that $\gamma\eta=p$ or $\gamma\eta=q$. The first case implies that $\gamma\lambda^*\in \{\tau_{\leq n}(p)\tau_{\leq n}(p)^*: n\geq0\}$ and $\gamma\lambda^*v=a_1p$, while the latter means $\gamma\lambda^*\in\{ \tau_{\leq m}(p)\tau_{\leq n}(q)^*: m,n\geq0\}$, and $\gamma\lambda^*v=a_1q$.
	
	\bigskip 
	
	\textit{\textbf{Case 2}: $q=\lambda\beta$ and $\gamma\beta\in V$}. It follows that $\gamma\beta=p$ or $\gamma\beta=q$. The same arguments show that $\gamma\lambda^*$ belongs to either $\{ \tau_{\leq m}(q)\tau_{\leq n}(p)^*: m,n\geq0\}$ and $\gamma\lambda^*v=a_2p$, or else $ \{ \tau_{\leq m}(q)\tau_{\leq m}(q)^*: m\geq0\}$, and $\gamma\lambda^*v=a_2q$.
	
	\bigskip 
	
	The computations just made allow us to conclude that $\gamma\lambda^* V\subseteq V$ if and only if $\gamma\lambda^* V=0$ (which means $\gamma\lambda^*\in I$), or else $\gamma\lambda^*$ belongs to the set
	$$
	S=\left\lbrace \tau_{\leq n}(p)\tau_{\leq n}(p)^*, \tau_{\leq m}(q)\tau_{\leq m}(q)^*, \tau_{\leq m}(p)\tau_{\leq n}(q)^*, \tau_{\leq m}(q)\tau_{\leq n}(p)^*\right\rbrace, 
	$$
	where $m,n\geq 0$. 
	This means that  
	$$
	R/I=\left\lbrace \sum_i k_i\gamma_i\lambda_i^*+I\;|\; k_i\in K,  \gamma_i\lambda_i^*\in S \text{ for all } 1\leq i\leq l\right\rbrace .
	$$
	Set 
	$$ 
	\begin{aligned}
		& S_1=\{ \tau_{\leq n}(p)\tau_{\leq n}(p)^*, n\geq0 \},	& S_2=\{ \tau_{\leq m}(q)\tau_{\leq m}(q)^*, m\geq0 \},
	\end{aligned}
	$$
	$$
	\begin{aligned}
		& S_3=\{ \tau_{\leq m}(p)\tau_{\leq n}(q)^*, m,n\geq0 \},	&  S_4=\{ \tau_{\leq m}(q)\tau_{\leq n}(p)^*, m,n\geq0 \},
	\end{aligned}
	$$
	Next, we claim that for each element $\sum_{i=1}^l k_i\gamma_i\lambda_i^*\in L_K(E)$, the following statements hold:
	\begin{enumerate}[font=\normalfont]
		\item[(1)] If $\gamma_i\lambda_i^*\in S_1$ for all $1\leq i\leq n$, then $\sum_{i=1}^{l} k_i\gamma_i\lambda_i^*-kr(f)\in I$ for some $k\in K$. (Note that $s(p)=r(f)$.)
		\item[(2)] If  $\gamma_i\lambda_i^*\in S_2$ for all $1\leq i\leq n$, then $\sum_{i=1}^{l} k_i\gamma_i\lambda_i^*-ks(f)\in I$ for some $k\in K$. (Note that $s(q)=s(f)$.)
		\item[(3)] If $\gamma_i\lambda_i^*\in S_3$ for all $1\leq i\leq n$, then $\sum_{i=1}^{l} k_i\gamma_i\lambda_i^*-kf^*\in I$, for some $k\in K$.
		\item[(4)] If  $\gamma_i\lambda_i^*\in S_4$ for all $1\leq i\leq n$, then $\sum_{i=1}^{l} k_i\gamma_i\lambda_i^*-kf\in I$ for some $k\in K$.
	\end{enumerate}
	First, we prove (1). Note that, as we have showed in Case 1, $\gamma_i\lambda_i^*v=a_1p$ for every $\gamma_i\lambda_i^*\in S_1$. The proof will be done by induction on $k$. Take an arbitrary element $v = a_1p+a_2q\in V$, where $a_1, a_2\in K$. If $l=1$, then
	$$
	\left( k_1\gamma_1\lambda_1^*-k_1s(p)\right) v=a_1k_1p-a_1k_1p=0,
	$$
	which implies that $k_1\gamma_1\lambda_1^*-k_1s(p)\in I$. Now, assume that $l>1$ and the statement is true for the sum of length less than $l$. Now, there exists $k'\in K$ such that 
	$\sum_{i=1}^{l-1} k_i\gamma_i\lambda_i^*-k's(q)\in I$, which implies that $\left( \sum_{i=1}^{l-1} k_i\gamma_i\lambda_i^*\right)=k's(q) v$. If we set $k=-k'-k_l$, then 
	$$
	\left( \sum_{i=1}^l k_i\gamma_i\lambda_i^*-ks(q)\right) v=\left( \sum_{i=1}^{l-1} k_i\gamma_i\lambda_i^*+k_l\gamma_l\lambda_l^*-ks(q)\right) v
	$$
	$$
	=k's(q)v+k_l\gamma_l\lambda_l^*v-ks(q) v=a_1k'p+a_1k_lp+a_1kp=0.
	$$
	This says that $\sum_{i=1}^{l} k_i\gamma_i\lambda_i^*-ks(q)\in I$. The proofs of (2), (3) and (4) follow similarly. Therefore, $R/I$ is described precisely as follows:
	$$
	R/I=\left\lbrace k_1s(f)+k_2 r(f)+k_3f+k_4f^*+I\;|\; k_1, k_2, k_3, k_4\in K\right\rbrace .
	$$
	In view of Lemma \ref{lemma_4.1}, one can easily check that $\{s(f)+I, r(f)+I, f+I, f^*+I\}$ is a $K$-basis of $R/I$. Therefore, the mapping $\varphi: R \to {\rm M}_2(K)$ given by 
	$$
	\begin{aligned}
		\varphi(s(f))=   f_{11}, & \hspace{1cm}
		\varphi(r(f))= f_{22},\\
		\varphi(f)= f_{12}, &  \hspace{1cm}
		\varphi(f^*)= f_{21}\\
	\end{aligned}
	$$
	is a surjective ring homomorphism with $\ker \varphi =I$. Consider the following matrices in ${\rm M}_2(K)$:
	\begin{align*}
		&	A             =    {\rm I}_2+\alpha f_{21},\\
		&	B             =    {\rm I}_2+\alpha f_{12},\\
		&	C             =    A+(\beta-1)f_{11}+(\beta^{-1}-1)f_{22},\\
		&	D             =    B+(\beta-1)f_{11}+(\beta^{-1}-1)f_{22}.
	\end{align*}
	Then, it is easy to see that
	\begin{align*}
		&	\varphi(a)             =    \varphi\left( 1_{L_K(E)}+\alpha f^*\right)  = A,\\
		&	\varphi(b)             =    \varphi\left( 1_{L_K(E)}+\alpha f\right)      = B,\\
		&	\varphi(c) 	            = 	\varphi\left( a +(\beta-1) s(f) + (\beta^{-1}-1) r(f)\right) = C,\\
		&	\varphi(d) 	            = 	\varphi\left( b+(\beta-1) s(f) + (\beta^{-1}-1) r(f)\right)= D.
	\end{align*}
		Moreover, the elements $a$, $b$, $c$, $d$ are invertible in $L_K(E)$ with 
		\begin{align*}
			&	a^{-1}     =   1_{L_K(E)}-\alpha f^*,\\
			&	b^{-1}     =   1_{L_K(E)}-\alpha f,\\
			&	c^{-1}     =   a^{-1} + (\beta^{-1}-1) s(f) + (\beta-1) r(f)\,\\
			&	d^{-1}     =   b^{-1} + (\beta^{-1}-1) s(f) + (\beta-1) r(f).
		\end{align*}
	
	To prove (i), we assume that ${\rm char} K =0$ and $\alpha$ is either equal to $2$ or transcendental over the prime subfield ${\mathbb{Q}}$. According to Lemma \ref{lemma_3.5}, $A$ and $B$ generate a non-cyclic free subgroup of ${\rm GL}_2(K)$. Since $\varphi: R\to {\rm M}_2(K)$ is surjective, we may apply Lemma \ref{lemma_3.1} to conclude that $\left\langle a, b\right\rangle $ is a  non-cyclic free subgroup of $R^\times \subseteq L_K(E)^\times$, proving (i). Finally, the proof of (ii) follows similarly.
\end{proof}

Now, there is a good moment to present the main theorem of this section.
\begin{theorem}\label{theorem_5.3}
	Let $E$ be a graph and $K$ a non-absolute field. Assume that $L_K(E)$ contains a primitive ideal $P$  of type II given in Theorem \ref{theorem_3.3} such that $L_K(E)/P$ is non-commutative, and set $H=P\cap E^0$. Then $E$ contains either a vertex $w$, which is either a sink or an infinite emitter, with $E^0\backslash H=M(w)$ or an infinite path $p$ with $E^0\backslash H=M(p)$. Moreover, for any $f\in E^1$ with $r(f)=w$ or $s(p)$,  the following assertions hold.
	\begin{enumerate}[font=\normalfont]
		\item[(i)] If ${\rm char}\; K =0$ and $\alpha$ is either equal to $2$ or transcendental over the prime subfield ${\mathbb{Q}}$, then 
		$$
		a:=1_{L_K(E)}+\alpha f^*,\hspace*{0,5cm} b:=1_{L_K(E)}+\alpha f
		$$
		generate a non-cyclic free subgroup in $L_K(E)$.
		\item[(ii)] If ${\rm char}\;K=t>0$ and $\alpha, \beta$ are algebraically independent over the prime subfield $\mathbb{F}_t$, then  
		$$c:=a +(\beta-1) s(f) + (\beta^{-1}-1) r(f),$$ $$d:=b+(\beta-1) s(f) + (\beta^{-1}-1) r(f) $$ generate a non-cyclic free subgroup in $L_K(E)$.
	\end{enumerate}
\end{theorem}

\begin{proof}
	If we set $F=E\backslash(H,B_H)$, then in view of Lemma~1.4, there is an epimorphism  $\varphi: L_K(E) \to L_K(F)$ with $\ker\varphi = P$. In this case, the quotient graph $F$ is described as follows:
	$$ 
	F^0= E^0\backslash H \text{ and } F^1=\{e\in E^1: r(e)\not\in H\}, \eqno(1)
	$$
	and $r$, $s$ are restricted to $F^0$. The proof of \cite[Theorem 3.9]{Pa_ara_rangas_2014} (see Part (ii)) shows that there are the following two cases to study:

	\textit{\textbf{Case 1}: $F^0$ contains a unique sink $w$ such that $F^0=M(w)$.} Consider $w$ as a vertex in $E^1$, we see that $w$ is either a sink or an infinite emitter in $E$ such that $r(s^{-1}(w)\subseteq H$. Set $A=\{e\in E^1|r(e)=w\}$. If $A=\varnothing$ then $F$ would be the single vertex $w$ which says that $L_K(F)$ is commutative, a contradiction. Therefore $A\ne\varnothing$ and we may take $f\in F^1\subseteq E^1$ with $s(f)=w$.  To prove (i) and (ii) for this case, we set
	\begin{align*}
		&	a'=1_{L_K(F)}+\alpha f^*, \\
		&	b'=1_{L_K(F)}+\alpha f,\\
		&	c'=a' +(\beta-1 )s(f) + (\beta^{-1}-1) w,\\
		&	d'=b'+(\beta-1 )s(f) + (\beta^{-1}-1) w,
	\end{align*}
	where $\alpha, \beta\in K$. Assume that ${\rm char} F=0$ and $\alpha$ is either equal to $2$ or transcendental over the prime subfield ${\mathbb{Q}}$. Now, by Lemma \ref{lemma_1.4}, we have $\varphi(f)=f$ and $\varphi(w)=w$. Since $w$ is a sink in $F$, we deduce that $s(f)\ne r(f)=w$ in $F$. Therefore, we may apply Lemma  \ref{lemma_4.2}(i) to conclude that $\left\langle a',b'\right\rangle $ is a non-cyclic free subgroup of $(L_K(F))^\times$. It is easy to check that $a$, $b$, $c$, $d$ are invertible in $L_K(E)$ with 
	\begin{align*}
		&	a^{-1}=1_{L_K(E)}-\alpha f^*, \\
		&	b^{-1}=1_{L_K(E)}-\alpha f,\\
		&	c^{-1}=a^{-1} +(\beta^{-1}-1) s(f) + (\beta-1) r(f),\\
		&	d^{-1}=b^{-1}+(\beta^{-1}-1) s(f) + (\beta-1) r(f).
	\end{align*}
	 Moreover, since $\varphi(a)=a'$, $\varphi(b)=b'$, $\varphi(c)=c'$, $\varphi(d)=d'$ and $\varphi$ is surjective, Lemma \ref{lemma_3.1} implies that $\left\langle a,b\right\rangle $ is a non-cyclic free subgroup of $(L_K(E))^\times$. Hence, (i) is established. The proof of (ii) is similar to the proof of (i), so we omit it.  
	
	\textit{\textbf{Case 2}: There is an infinite path $p_0$ in $F$ such that $F^0=M(p_0)$}: Because $F$ is a subgraph of $E$, we can consider $p_0$ as an infinite path in $E$ which does not end in an exclusive cycle in $E$ such that $E^0=M(p_0)$. Let $V_{[p_0]}$ be the faithful $L_K(F)$-simple module determined by the infinite path $p_0$; that is, 
	$$
	V_{[p_0]}=\bigoplus_{q\in [p_0]}Kq,
	$$
	a $K$-vector space whose basis consists of all infinite paths $q$ which are tail- equivalent to $p_0$. According to Lemma \ref{lemma_2.1}, we conclude that ${\rm End}_{L_K(F)}(V_{[p_0]})\cong K$. If $L_K(F)$ is artinian, then Theorem \ref{theorem_3.4} implies that $L_K(F)\cong {\rm M}_n(K)$, where $n=\dim_K(V_{[p_0]})$ which is the number of infinite paths in $F$ tail-equivalent to $p_0$. Since $L_K(F)$ is non-commutative, it follows that $n\geq 2$. If $L_K(F)$ is non-artinian, then $\dim_K(V_{[p_0]})=~\infty$. In either cases, we have $\dim_K(V_{[p_0]})\geq 2$, and so there exist in $F$ two infinite distinct paths $p'$ and $q'$ which are tail-equivalent to $p_0$. According to Lemma \ref{lemma_5.1}, we conclude that there exist two infinite paths $p$ and $fp$, where $f\in F^1\subseteq E^1$ with $s(f)\ne r(f)=s(p)$. The proof of (i) and (ii) for this case is similar to that of previous case.
\end{proof}
\begin{example}
	Consider the graph $E$ given in Example \ref{exmaple_4}.
	
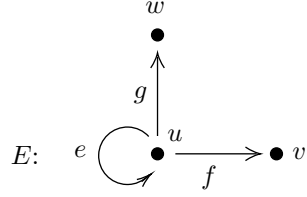
	\begin{center}
		
		\tikzset{every picture/.style={line width=0.5pt}} 
		
		\begin{tikzpicture}[x=0.75pt,y=0.75pt,yscale=-1,xscale=1]
		
		\draw  [color={rgb, 255:red, 0; green, 0; blue, 0 }  ,draw opacity=1 ][fill={rgb, 255:red, 0; green, 0; blue, 0 }  ,fill opacity=1 ] (87,135) .. controls (87,133.34) and (88.34,132) .. (90,132) .. controls (91.66,132) and (93,133.34) .. (93,135) .. controls (93,136.66) and (91.66,138) .. (90,138) .. controls (88.34,138) and (87,136.66) .. (87,135) -- cycle ;
		
		\draw  [color={rgb, 255:red, 0; green, 0; blue, 0 }  ,draw opacity=1 ][fill={rgb, 255:red, 0; green, 0; blue, 0 }  ,fill opacity=1 ] (147,135) .. controls (147,133.34) and (148.34,132) .. (150,132) .. controls (151.66,132) and (153,133.34) .. (153,135) .. controls (153,136.66) and (151.66,138) .. (150,138) .. controls (148.34,138) and (147,136.66) .. (147,135) -- cycle ;
		
		\draw    (99,135) -- (140,135) ;
		\draw [shift={(142,135)}, rotate = 180] [color={rgb, 255:red, 0; green, 0; blue, 0 }  ][line width=0.5]    (7.65,-3.43) .. controls (4.86,-1.61) and (2.31,-0.47) .. (0,0) .. controls (2.31,0.47) and (4.86,1.61) .. (7.65,3.43)   ;
		
		\draw  [color={rgb, 255:red, 0; green, 0; blue, 0 }  ,draw opacity=1 ][fill={rgb, 255:red, 0; green, 0; blue, 0 }  ,fill opacity=1 ] (87,75) .. controls (87,73.34) and (88.34,72) .. (90,72) .. controls (91.66,72) and (93,73.34) .. (93,75) .. controls (93,76.66) and (91.66,78) .. (90,78) .. controls (88.34,78) and (87,76.66) .. (87,75) -- cycle ;
		
		\draw    (90,126) -- (90,86) ;
		\draw [shift={(90,84)}, rotate = 450] [color={rgb, 255:red, 0; green, 0; blue, 0 }  ][line width=0.5]    (7.65,-3.43) .. controls (4.86,-1.61) and (2.31,-0.47) .. (0,0) .. controls (2.31,0.47) and (4.86,1.61) .. (7.65,3.43)   ;
		
		\draw  [draw opacity=0] (85.77,145.48) .. controls (83,148.67) and (78.87,150.63) .. (74.31,150.47) .. controls (66.31,150.19) and (60.05,143.48) .. (60.33,135.48) .. controls (60.61,127.48) and (67.32,121.22) .. (75.32,121.5) .. controls (79.4,121.64) and (83.03,123.46) .. (85.56,126.26) -- (74.82,135.99) -- cycle ; \draw   (85.77,145.48) .. controls (83,148.67) and (78.87,150.63) .. (74.31,150.47) .. controls (66.31,150.19) and (60.05,143.48) .. (60.33,135.48) .. controls (60.61,127.48) and (67.32,121.22) .. (75.32,121.5) .. controls (79.4,121.64) and (83.03,123.46) .. (85.56,126.26) ;
		
		\draw    (82,149) -- (85.44,146.25) ;
		\draw [shift={(87,145)}, rotate = 501.34] [color={rgb, 255:red, 0; green, 0; blue, 0 }  ][line width=0.5]    (7.65,-2.3) .. controls (4.86,-0.97) and (2.31,-0.21) .. (0,0) .. controls (2.31,0.21) and (4.86,0.98) .. (7.65,2.3)   ;
		
		\draw (-27,165) node [anchor=north west][inner sep=0.5pt]   [align=left] {$ $};
		
		\draw (15,130) node [anchor=north west][inner sep=0.5pt]   [align=left] {$\displaystyle E$:};
		
		\draw (94,122) node [anchor=north west][inner sep=0.5pt]   [align=left] {$\displaystyle u$};
		
		\draw (157,131) node [anchor=north west][inner sep=0.5pt]   [align=left] {$\displaystyle v$};
		
		\draw (83,57) node [anchor=north west][inner sep=0.5pt]   [align=left] {$\displaystyle w$};
		
		\draw (77,100) node [anchor=north west][inner sep=0.5pt]   [align=left] {$\displaystyle g$};
		
		\draw (111,140) node [anchor=north west][inner sep=0.5pt]   [align=left] {$\displaystyle f$};
		
		\draw (47,130) node [anchor=north west][inner sep=0.5pt]   [align=left] {$\displaystyle e$};
		\end{tikzpicture}
	\end{center}
	
	As we have explained in Example \ref{exmaple_4},  $H=\{w\}$  is a hereditary saturated subset of $E^0$ and $P:=I(H)$ is a primitive ideal of type II (note that $B_H=\varnothing$). Also $E^0\backslash H=\{u,v\}=M(v)$. According to Theorem \ref{theorem_5.3}, if $K$ is a field of characteristic $0$ and $\alpha$ is either equal to $2$ or transcendental over the prime subfield ${\mathbb{Q}}$, then the following elements are the generators of a free subgroup of $L_K(E)^\times$: 
	$$
	a=v+u+w+\alpha f^*,\;\;\;b=v+u+w+\alpha f.
	$$
\end{example}
\begin{example}
	Let $K$ be a non-absolute field. Consider the graphs $E$ and $F$ given bellow:
	
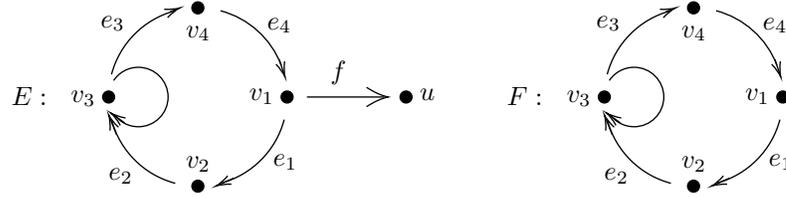
	\begin{center}
		
		\tikzset{every picture/.style={line width=0.5pt}}  
		
		\begin{tikzpicture}[x=0.75pt,y=0.75pt,yscale=-1,xscale=1]
		
		\draw  [color={rgb, 255:red, 0; green, 0; blue, 0 }  ,draw opacity=1 ][fill={rgb, 255:red, 0; green, 0; blue, 0 }  ,fill opacity=1 ] (117,90) .. controls (117,88.34) and (118.34,87) .. (120,87) .. controls (121.66,87) and (123,88.34) .. (123,90) .. controls (123,91.66) and (121.66,93) .. (120,93) .. controls (118.34,93) and (117,91.66) .. (117,90) -- cycle ;
		
		\draw  [color={rgb, 255:red, 0; green, 0; blue, 0 }  ,draw opacity=1 ][fill={rgb, 255:red, 0; green, 0; blue, 0 }  ,fill opacity=1 ] (117,180) .. controls (117,178.34) and (118.34,177) .. (120,177) .. controls (121.66,177) and (123,178.34) .. (123,180) .. controls (123,181.66) and (121.66,183) .. (120,183) .. controls (118.34,183) and (117,181.66) .. (117,180) -- cycle ;
		
		\draw  [color={rgb, 255:red, 0; green, 0; blue, 0 }  ,draw opacity=1 ][fill={rgb, 255:red, 0; green, 0; blue, 0 }  ,fill opacity=1 ] (72,135) .. controls (72,133.34) and (73.34,132) .. (75,132) .. controls (76.66,132) and (78,133.34) .. (78,135) .. controls (78,136.66) and (76.66,138) .. (75,138) .. controls (73.34,138) and (72,136.66) .. (72,135) -- cycle ;
		
		\draw  [color={rgb, 255:red, 0; green, 0; blue, 0 }  ,draw opacity=1 ][fill={rgb, 255:red, 0; green, 0; blue, 0 }  ,fill opacity=1 ] (162,135) .. controls (162,133.34) and (163.34,132) .. (165,132) .. controls (166.66,132) and (168,133.34) .. (168,135) .. controls (168,136.66) and (166.66,138) .. (165,138) .. controls (163.34,138) and (162,136.66) .. (162,135) -- cycle ;
		
		\draw  [draw opacity=0] (76.46,123.59) .. controls (80.52,108.05) and (92.68,95.79) .. (108.15,91.58) -- (120,135) -- cycle ; \draw   (76.46,123.59) .. controls (80.52,108.05) and (92.68,95.79) .. (108.15,91.58) ;
		
		\draw  [draw opacity=0] (108.59,178.54) .. controls (93.05,174.48) and (80.79,162.32) .. (76.58,146.85) -- (120,135) -- cycle ; \draw   (108.59,178.54) .. controls (93.05,174.48) and (80.79,162.32) .. (76.58,146.85) ;
		
		\draw  [draw opacity=0] (131.41,91.46) .. controls (145.52,95.15) and (156.93,105.51) .. (162.06,118.98) -- (120,135) -- cycle ; \draw   (131.41,91.46) .. controls (145.52,95.15) and (156.93,105.51) .. (162.06,118.98) ;
		
		\draw  [draw opacity=0] (163.54,146.41) .. controls (159.48,161.95) and (147.32,174.21) .. (131.85,178.42) -- (120,135) -- cycle ; \draw   (163.54,146.41) .. controls (159.48,161.95) and (147.32,174.21) .. (131.85,178.42) ;
		
		\draw  [draw opacity=0] (77.52,126.68) .. controls (80.21,122.65) and (84.79,120) .. (90,120) .. controls (98.28,120) and (105,126.72) .. (105,135) .. controls (105,143.28) and (98.28,150) .. (90,150) .. controls (85.03,150) and (80.63,147.59) .. (77.9,143.87) -- (90,135) -- cycle ; \draw   (77.52,126.68) .. controls (80.21,122.65) and (84.79,120) .. (90,120) .. controls (98.28,120) and (105,126.72) .. (105,135) .. controls (105,143.28) and (98.28,150) .. (90,150) .. controls (85.03,150) and (80.63,147.59) .. (77.9,143.87) ;
		
		\draw  [color={rgb, 255:red, 0; green, 0; blue, 0 }  ,draw opacity=1 ][fill={rgb, 255:red, 0; green, 0; blue, 0 }  ,fill opacity=1 ] (222,135) .. controls (222,133.34) and (223.34,132) .. (225,132) .. controls (226.66,132) and (228,133.34) .. (228,135) .. controls (228,136.66) and (226.66,138) .. (225,138) .. controls (223.34,138) and (222,136.66) .. (222,135) -- cycle ;
		
		\draw    (175,135) -- (214,135) ;
		\draw [shift={(216,135)}, rotate = 180] [color={rgb, 255:red, 0; green, 0; blue, 0 }  ][line width=0.5]    (10.93,-4.9) .. controls (6.95,-2.3) and (3.31,-0.67) .. (0,0) .. controls (3.31,0.67) and (6.95,2.3) .. (10.93,4.9)   ;
		
		\draw    (106,92) -- (109.14,90.74) ;
		\draw [shift={(111,90)}, rotate = 518.2] [color={rgb, 255:red, 0; green, 0; blue, 0 }  ][line width=0.5]    (7.65,-2.3) .. controls (4.86,-0.97) and (2.31,-0.21) .. (0,0) .. controls (2.31,0.21) and (4.86,0.98) .. (7.65,2.3)   ;
		
		\draw    (161,116) -- (163.37,123.1) ;
		\draw [shift={(164,125)}, rotate = 251.57] [color={rgb, 255:red, 0; green, 0; blue, 0 }  ][line width=0.5]    (7.65,-2.3) .. controls (4.86,-0.97) and (2.31,-0.21) .. (0,0) .. controls (2.31,0.21) and (4.86,0.98) .. (7.65,2.3)   ;
		
		\draw    (138,176) -- (130.83,179.19) ;
		\draw [shift={(129,180)}, rotate = 336.03999999999996] [color={rgb, 255:red, 0; green, 0; blue, 0 }  ][line width=0.5]    (7.65,-2.3) .. controls (4.86,-0.97) and (2.31,-0.21) .. (0,0) .. controls (2.31,0.21) and (4.86,0.98) .. (7.65,2.3)   ;
		
		\draw    (79,154) -- (75.74,145.86) ;
		\draw [shift={(75,144)}, rotate = 428.2] [color={rgb, 255:red, 0; green, 0; blue, 0 }  ][line width=0.5]    (7.65,-2.3) .. controls (4.86,-0.97) and (2.31,-0.21) .. (0,0) .. controls (2.31,0.21) and (4.86,0.98) .. (7.65,2.3)   ;
		
		\draw    (83,148) -- (77.63,144.16) ;
		\draw [shift={(76,143)}, rotate = 395.53999999999996] [color={rgb, 255:red, 0; green, 0; blue, 0 }  ][line width=0.5]    (7.65,-2.3) .. controls (4.86,-0.97) and (2.31,-0.21) .. (0,0) .. controls (2.31,0.21) and (4.86,0.98) .. (7.65,2.3)   ;
		
		\draw  [color={rgb, 255:red, 0; green, 0; blue, 0 }  ,draw opacity=1 ][fill={rgb, 255:red, 0; green, 0; blue, 0 }  ,fill opacity=1 ] (367,90) .. controls (367,88.34) and (368.34,87) .. (370,87) .. controls (371.66,87) and (373,88.34) .. (373,90) .. controls (373,91.66) and (371.66,93) .. (370,93) .. controls (368.34,93) and (367,91.66) .. (367,90) -- cycle ;
		
		\draw  [color={rgb, 255:red, 0; green, 0; blue, 0 }  ,draw opacity=1 ][fill={rgb, 255:red, 0; green, 0; blue, 0 }  ,fill opacity=1 ] (367,180) .. controls (367,178.34) and (368.34,177) .. (370,177) .. controls (371.66,177) and (373,178.34) .. (373,180) .. controls (373,181.66) and (371.66,183) .. (370,183) .. controls (368.34,183) and (367,181.66) .. (367,180) -- cycle ;
		
		\draw  [color={rgb, 255:red, 0; green, 0; blue, 0 }  ,draw opacity=1 ][fill={rgb, 255:red, 0; green, 0; blue, 0 }  ,fill opacity=1 ] (322,135) .. controls (322,133.34) and (323.34,132) .. (325,132) .. controls (326.66,132) and (328,133.34) .. (328,135) .. controls (328,136.66) and (326.66,138) .. (325,138) .. controls (323.34,138) and (322,136.66) .. (322,135) -- cycle ;
		
		\draw  [color={rgb, 255:red, 0; green, 0; blue, 0 }  ,draw opacity=1 ][fill={rgb, 255:red, 0; green, 0; blue, 0 }  ,fill opacity=1 ] (412,135) .. controls (412,133.34) and (413.34,132) .. (415,132) .. controls (416.66,132) and (418,133.34) .. (418,135) .. controls (418,136.66) and (416.66,138) .. (415,138) .. controls (413.34,138) and (412,136.66) .. (412,135) -- cycle ;
		
		\draw  [draw opacity=0] (326.46,123.59) .. controls (330.52,108.05) and (342.68,95.79) .. (358.15,91.58) -- (370,135) -- cycle ; \draw   (326.46,123.59) .. controls (330.52,108.05) and (342.68,95.79) .. (358.15,91.58) ;
		
		\draw  [draw opacity=0] (358.59,178.54) .. controls (343.05,174.48) and (330.79,162.32) .. (326.58,146.85) -- (370,135) -- cycle ; \draw   (358.59,178.54) .. controls (343.05,174.48) and (330.79,162.32) .. (326.58,146.85) ;
		
		\draw  [draw opacity=0] (381.41,91.46) .. controls (395.52,95.15) and (406.93,105.51) .. (412.06,118.98) -- (370,135) -- cycle ; \draw   (381.41,91.46) .. controls (395.52,95.15) and (406.93,105.51) .. (412.06,118.98) ;
		
		\draw  [draw opacity=0] (413.54,146.41) .. controls (409.48,161.95) and (397.32,174.21) .. (381.85,178.42) -- (370,135) -- cycle ; \draw   (413.54,146.41) .. controls (409.48,161.95) and (397.32,174.21) .. (381.85,178.42) ;
		
		\draw  [draw opacity=0] (327.52,126.68) .. controls (330.21,122.65) and (334.79,120) .. (340,120) .. controls (348.28,120) and (355,126.72) .. (355,135) .. controls (355,143.28) and (348.28,150) .. (340,150) .. controls (335.03,150) and (330.63,147.59) .. (327.9,143.87) -- (340,135) -- cycle ; \draw   (327.52,126.68) .. controls (330.21,122.65) and (334.79,120) .. (340,120) .. controls (348.28,120) and (355,126.72) .. (355,135) .. controls (355,143.28) and (348.28,150) .. (340,150) .. controls (335.03,150) and (330.63,147.59) .. (327.9,143.87) ;
		
		\draw    (356,92) -- (359.14,90.74) ;
		\draw [shift={(361,90)}, rotate = 518.2] [color={rgb, 255:red, 0; green, 0; blue, 0 }  ][line width=0.5]    (7.65,-2.3) .. controls (4.86,-0.97) and (2.31,-0.21) .. (0,0) .. controls (2.31,0.21) and (4.86,0.98) .. (7.65,2.3)   ;
		
		\draw    (411,116) -- (413.37,123.1) ;
		\draw [shift={(414,125)}, rotate = 251.57] [color={rgb, 255:red, 0; green, 0; blue, 0 }  ][line width=0.5]    (7.65,-2.3) .. controls (4.86,-0.97) and (2.31,-0.21) .. (0,0) .. controls (2.31,0.21) and (4.86,0.98) .. (7.65,2.3)   ;
		
		\draw    (388,176) -- (380.83,179.19) ;
		\draw [shift={(379,180)}, rotate = 336.03999999999996] [color={rgb, 255:red, 0; green, 0; blue, 0 }  ][line width=0.5]    (7.65,-2.3) .. controls (4.86,-0.97) and (2.31,-0.21) .. (0,0) .. controls (2.31,0.21) and (4.86,0.98) .. (7.65,2.3)   ;
		
		\draw    (329,154) -- (325.74,145.86) ;
		\draw [shift={(325,144)}, rotate = 428.2] [color={rgb, 255:red, 0; green, 0; blue, 0 }  ][line width=0.5]    (7.65,-2.3) .. controls (4.86,-0.97) and (2.31,-0.21) .. (0,0) .. controls (2.31,0.21) and (4.86,0.98) .. (7.65,2.3)   ;
		
		\draw    (333,148) -- (327.63,144.16) ;
		\draw [shift={(326,143)}, rotate = 395.53999999999996] [color={rgb, 255:red, 0; green, 0; blue, 0 }  ][line width=0.5]    (7.65,-2.3) .. controls (4.86,-0.97) and (2.31,-0.21) .. (0,0) .. controls (2.31,0.21) and (4.86,0.98) .. (7.65,2.3)   ;
		
		\draw (-27,165) node [anchor=north west][inner sep=0.5pt]   [align=left] {$ $};
		
		\draw (25,129) node [anchor=north west][inner sep=0.5pt]   [align=left] {$\displaystyle E:$};
		
		\draw (55,131) node [anchor=north west][inner sep=0.5pt]   [align=left] {$\displaystyle v_{3}$};
		
		\draw (231,130) node [anchor=north west][inner sep=0.5pt]   [align=left] {$\displaystyle u$};
		
		\draw (113,98) node [anchor=north west][inner sep=0.5pt]   [align=left] {$\displaystyle v_{4}$};
		
		\draw (113,164) node [anchor=north west][inner sep=0.5pt]   [align=left] {$\displaystyle v_{2}$};
		
		\draw (145,130) node [anchor=north west][inner sep=0.5pt]   [align=left] {$\displaystyle v_{1}$};
		
		\draw (70,93) node [anchor=north west][inner sep=0.5pt]   [align=left] {$\displaystyle e_{3}$};
		
		\draw (154,94) node [anchor=north west][inner sep=0.5pt]   [align=left] {$\displaystyle e_{4}$};
		
		\draw (74,170) node [anchor=north west][inner sep=0.5pt]   [align=left] {$\displaystyle e_{2}$};
		
		\draw (157,163) node [anchor=north west][inner sep=0.5pt]   [align=left] {$\displaystyle e_{1}$};
		
		\draw (186,116) node [anchor=north west][inner sep=0.5pt]   [align=left] {$\displaystyle f$};
		
		\draw (275,129) node [anchor=north west][inner sep=0.5pt]   [align=left] {$\displaystyle F:$};
		
		\draw (305,131) node [anchor=north west][inner sep=0.5pt]   [align=left] {$\displaystyle v_{3}$};
		
		\draw (363,98) node [anchor=north west][inner sep=0.5pt]   [align=left] {$\displaystyle v_{4}$};
		
		\draw (363,164) node [anchor=north west][inner sep=0.5pt]   [align=left] {$\displaystyle v_{2}$};
		
		\draw (395,130) node [anchor=north west][inner sep=0.5pt]   [align=left] {$\displaystyle v_{1}$};
		
		\draw (320,93) node [anchor=north west][inner sep=0.5pt]   [align=left] {$\displaystyle e_{3}$};
		
		\draw (404,94) node [anchor=north west][inner sep=0.5pt]   [align=left] {$\displaystyle e_{4}$};
		
		\draw (324,170) node [anchor=north west][inner sep=0.5pt]   [align=left] {$\displaystyle e_{2}$};
		
		\draw (407,163) node [anchor=north west][inner sep=0.5pt]   [align=left] {$\displaystyle e_{1}$};
		
		\end{tikzpicture}
	\end{center}
	Then $H=\{u\}$ is a hereditary saturated subset of $E^0$ and $B_H=\varnothing$. It is clear that $E\backslash H=F$ which is downward directed, satisfying the Condition (L) and the Countable Separation Property. This implies that the ideal $I(H)$ is a primitive ideal of type II, and $L_K(E)/I(H)\cong L_K(F)$ which is non-commutative. Let $c=e_1e_2e_3e_4$ be a cycle in $E$ and $p=c^\infty$. Then $p$ is an infinite path in $E$ with $E^0\backslash H=M(p)$. Also, the edge $e_4$ satisfies the condition $v_4=s(e_4)\ne r(e_4)=v_1=s(p)$. Therefore, in view of Theorem \ref{theorem_5.3}, if ${\rm char} K =0$ and $\alpha$ is either equal to $2$ or transcendental over ${\mathbb{Q}}$ then the following elements generate a non-cyclic free subgroup of $L_K(E)^\times$:
	$$
	a=u+v_1+v_2+v_3+v_4+\alpha e_4^*,\;\;\;b=u+v_1+v_2+v_3+v_4+\alpha e_4.
	$$ 	
	We note in passing that there are many choices of infinite paths in $E$ which produce different non-cyclic free subgroups. For example, we may choose $q=e^\infty$ and the edge $e_2$ that satisfies $s(e_2)\ne r(e_2)=s(q)$.
\end{example}

As it was proved in \cite[Thereom 3.10]{Pa_ara_rangas_2014}, the Leavitt path algebra $L_K(E)$ is primitive if and only if $E$ contains either a sink $w$ or an infinite path $p$ that does not end in a cycle without exits such that $E^0=M(w)$ or $M(p)$. Here, as a consequence of the previous theorem, we give in the following corollary a description of the free subgroups in a primitive Leavitt path algebra.

\begin{corollary}\label{corollary_5.4}
	Let $E$ be a graph and $K$ a non-absolute field. If $L_K(E)$ is a non-commutative primitive Leavitt path algebra, then $E$ contains either a sink $w$ or an infinite path $p$ that does not end in a cycle without exits such that $E^0=M(w)$ or $M(p)$. Moreover, for any $f\in E^1$ with $ r(f)=w$ or $s(p)$, we have the followings:
	\begin{enumerate}[font=\normalfont]
		\item[(i)] If ${\rm char}\; K =0$ and $\alpha$ is either equal to $2$ or transcendental over the prime subfield ${\mathbb{Q}}$, then 
		$$
		a:=1_{L_K(E)}+\alpha f^*,\hspace*{0,5cm} b:=1_{L_K(E)}+\alpha f
		$$
		generate a non-cyclic free subgroup in $L_K(E)$.
		\item[(ii)] If ${\rm char}\;K=t>0$ and $\alpha, \beta$ are algebraically independent over the prime subfield $\mathbb{F}_t$, then  
		$$c:=a +(\beta-1) s(f) + (\beta^{-1}-1) r(f),$$ $$d:=b+(\beta-1) s(f) + (\beta^{-1}-1) r(f) $$ generate a non-cyclic free subgroup in $L_K(E)$.
	\end{enumerate}
\end{corollary}

\begin{proof} 
	In the proof of the precedent theorem, we see that the quotient graph $F$ contains either a sink or an infinite path, namely $w$ and $p$ respectively, and an edge $f$ with $r(f)=w$ or $s(p)$. Under these conditions we have proved that $L_K(F)^{\times}$ contains a non-cyclic free subgroup. Since the graph $E$ in the present corollary satisfies all conditions as for the graph $F$ mentioned above, the proof of the existence of non-cyclic free subgroups in $L_K(F)$ is applied similarly for $L_K(E)$. Hence, we can conclude that the proof of the corollary is complete. 
\end{proof}
\begin{example}
	Consider the algebraic Toeplitz $K$-algebra $\mathscr{T}_K$ of the Toeplitz graph. Then, we have $1_{\mathscr{T}_K}=u+v$. Thus, if $K$ is a field of characteristic $t>0$ and $\alpha, \beta$ are algebraically independent over the prime subfield $\mathbb{F}_t$, then the following elements will generate a non-cyclic free subgroup of $\mathscr{T}_K^\times$:
	$$
	c=u+v+\alpha f^*+(\beta-1) u+(\beta^{-1}-1)v,
	$$
	$$
	d=u+v+\alpha f+(\beta-1)u+(\beta^{-1}-1)v.
	$$	
\end{example}

\section{Free subgroups determined by the primitive  ideals of type III}
This section is devoted to determining non-cyclic free subgroups of $L_K(E)$ by a primitive ideal of type III. Recall that an ideal $P$ of $L_K(E)$ is  primitive  of type III if $P=I(H, B_H, f(c))$, where $H=P\cap E^0$, $c$ is an exclusive cycle based at a vertex $u$, $E^0\backslash H=M(u)$, and $f(x)$ is an irreducible polynomial in $K[x, x^{-1}]$. 
\begin{theorem} 
	Let $E$ be a graph and $K$ a non-absolute field. Assume that $L_K(E)$ contains a primitive ideal $P$  of type III given in Theorem \ref{theorem_3.3} (that is, the ideal $P$ is defined by $P=I(H, B_H, f(c))$, where $H=P\cap E^0$ and $c$ is an exclusive cycle) such that $L_K(E)/P$ is non-commutative. Then $E$ contains an infinite path $p$ which is tail-equivalent to $c$ such that $E^0\backslash H=M(p)$. Moreover, for any $f\in E^1$ with $ r(f)=s(p)$, the following assertions hold.
	\begin{enumerate}[font=\normalfont]
		\item[(i)] If ${\rm char}\; K =0$ and $\alpha$ is either equal to $2$ or transcendental over the prime subfield ${\mathbb{Q}}$, then 
		$$
		a:=1_{L_K(E)}+\alpha f^*,\hspace*{0,5cm} b:=1_{L_K(E)}+\alpha f
		$$
		generate a non-cyclic free subgroup in $L_K(E)$.
		\item[(ii)] If ${\rm char}\;K=t>0$ and $\alpha, \beta$ are algebraically independent over the prime subfield $\mathbb{F}_t$, then  
		$$c:=a +(\beta-1) s(f) + (\beta^{-1}-1) r(f),$$ $$d:=b+(\beta-1) s(f) + (\beta^{-1}-1) r(f) $$ generate a non-cyclic free subgroup in $L_K(E)$.
	\end{enumerate}
\end{theorem}
\begin{proof}
	 Consider the quotient graph $F=E\backslash (H, B_H)$. By Lemma \ref{lemma_1.4}, $I(H, B_H)$  is a graded ideal of $L_K(E)$ and there is an epimorphism $\varphi: L_K(E) \to L_K(F)$ with $\ker\varphi = I(H,B_H)$. It is clear that $F$ contains $c$ which is an exclusive cycle based at $u$ and $F^0=M(u)$. Let $V_{[c]}$ be the faithful simple $L_K(F)$-simple module determined by $c$. The arguments used in Case 2 of Theorem \ref{theorem_5.3} shows that $\dim_K(V_{[c]})\geq 2$, so there exist in $F$ two infinite distinct paths $p'$ and $q'$ which are tail-equivalent to $c$. According to Lemma \ref{lemma_5.1}, we conclude that there exist two infinite paths $p$ and $fp$ where $f\in F^1\subseteq E^1$ with $s(f)\ne r(f)=s(p)$. Therefore the result follows from Lemma \ref{lemma_5.3}.
\end{proof}
\begin{example}
	Let $K$ be a non-absolute field. Consider the graph $E$ given in Example \ref{exmaple_5}. 
	
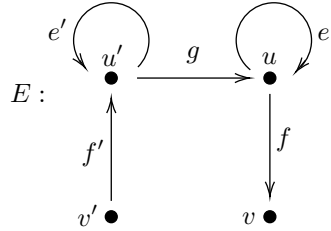
	\begin{center}
		
		\tikzset{every picture/.style={line width=0.5pt}} 
		\begin{tikzpicture}[x=0.75pt,y=0.75pt,yscale=-1,xscale=1]
		
		\draw  [color={rgb, 255:red, 0; green, 0; blue, 0 }  ,draw opacity=1 ][fill={rgb, 255:red, 0; green, 0; blue, 0 }  ,fill opacity=1 ] (182,135) .. controls (182,133.34) and (183.34,132) .. (185,132) .. controls (186.66,132) and (188,133.34) .. (188,135) .. controls (188,136.66) and (186.66,138) .. (185,138) .. controls (183.34,138) and (182,136.66) .. (182,135) -- cycle ;
		
		\draw  [draw opacity=0] (175.13,130.84) .. controls (170.78,127.36) and (168,122) .. (168,116) .. controls (168,105.51) and (176.51,97) .. (187,97) .. controls (197.49,97) and (206,105.51) .. (206,116) .. controls (206,121.25) and (203.87,126) .. (200.44,129.44) -- (187,116) -- cycle ; \draw   (175.13,130.84) .. controls (170.78,127.36) and (168,122) .. (168,116) .. controls (168,105.51) and (176.51,97) .. (187,97) .. controls (197.49,97) and (206,105.51) .. (206,116) .. controls (206,121.25) and (203.87,126) .. (200.44,129.44) ;
		
		\draw    (203,126) -- (199.28,130.46) ;
		\draw [shift={(198,132)}, rotate = 309.81] [color={rgb, 255:red, 0; green, 0; blue, 0 }  ][line width=0.5]    (7.65,-2.3) .. controls (4.86,-0.97) and (2.31,-0.21) .. (0,0) .. controls (2.31,0.21) and (4.86,0.98) .. (7.65,2.3)   ;
		
		\draw  [color={rgb, 255:red, 0; green, 0; blue, 0 }  ,draw opacity=1 ][fill={rgb, 255:red, 0; green, 0; blue, 0 }  ,fill opacity=1 ] (102,135) .. controls (102,133.34) and (103.34,132) .. (105,132) .. controls (106.66,132) and (108,133.34) .. (108,135) .. controls (108,136.66) and (106.66,138) .. (105,138) .. controls (103.34,138) and (102,136.66) .. (102,135) -- cycle ;
		
		\draw    (89,126) -- (90.8,128.4) ;
		\draw [shift={(92,130)}, rotate = 233.13] [color={rgb, 255:red, 0; green, 0; blue, 0 }  ][line width=0.5]    (7.65,-2.3) .. controls (4.86,-0.97) and (2.31,-0.21) .. (0,0) .. controls (2.31,0.21) and (4.86,0.98) .. (7.65,2.3)   ;
		
		\draw    (118,135) -- (173,135) ;
		\draw [shift={(175,135)}, rotate = 180] [color={rgb, 255:red, 0; green, 0; blue, 0 }  ][line width=0.5]    (7.65,-2.3) .. controls (4.86,-0.97) and (2.31,-0.21) .. (0,0) .. controls (2.31,0.21) and (4.86,0.98) .. (7.65,2.3)   ;
		
		\draw  [color={rgb, 255:red, 0; green, 0; blue, 0 }  ,draw opacity=1 ][fill={rgb, 255:red, 0; green, 0; blue, 0 }  ,fill opacity=1 ] (102,205) .. controls (102,203.34) and (103.34,202) .. (105,202) .. controls (106.66,202) and (108,203.34) .. (108,205) .. controls (108,206.66) and (106.66,208) .. (105,208) .. controls (103.34,208) and (102,206.66) .. (102,205) -- cycle ;
		
		\draw    (105,197) -- (105,147) ;
		\draw [shift={(105,145)}, rotate = 450] [color={rgb, 255:red, 0; green, 0; blue, 0 }  ][line width=0.5]    (7.65,-2.3) .. controls (4.86,-0.97) and (2.31,-0.21) .. (0,0) .. controls (2.31,0.21) and (4.86,0.98) .. (7.65,2.3)   ;
		
		\draw  [color={rgb, 255:red, 0; green, 0; blue, 0 }  ,draw opacity=1 ][fill={rgb, 255:red, 0; green, 0; blue, 0 }  ,fill opacity=1 ] (182,205) .. controls (182,203.34) and (183.34,202) .. (185,202) .. controls (186.66,202) and (188,203.34) .. (188,205) .. controls (188,206.66) and (186.66,208) .. (185,208) .. controls (183.34,208) and (182,206.66) .. (182,205) -- cycle ;
		
		\draw    (185,143) -- (185,193) ;
		\draw [shift={(185,195)}, rotate = 270] [color={rgb, 255:red, 0; green, 0; blue, 0 }  ][line width=0.5]    (7.65,-2.3) .. controls (4.86,-0.97) and (2.31,-0.21) .. (0,0) .. controls (2.31,0.21) and (4.86,0.98) .. (7.65,2.3)   ;
		
		\draw  [draw opacity=0] (89.68,127.24) .. controls (87.37,124.09) and (86,120.2) .. (86,116) .. controls (86,105.51) and (94.51,97) .. (105,97) .. controls (115.49,97) and (124,105.51) .. (124,116) .. controls (124,121.25) and (121.87,126) .. (118.44,129.44) -- (105,116) -- cycle ; \draw   (89.68,127.24) .. controls (87.37,124.09) and (86,120.2) .. (86,116) .. controls (86,105.51) and (94.51,97) .. (105,97) .. controls (115.49,97) and (124,105.51) .. (124,116) .. controls (124,121.25) and (121.87,126) .. (118.44,129.44) ;

		\draw (-27,165) node [anchor=north west][inner sep=0.5pt]   [align=left] {$ $};
		
		\draw (53,137) node [anchor=north west][inner sep=0.5pt]   [align=left] {$\displaystyle E:$};
		
		\draw (180,121) node [anchor=north west][inner sep=0.5pt]   [align=left] {$\displaystyle u$};
		
		\draw (170,202) node [anchor=north west][inner sep=0.5pt]   [align=left] {$\displaystyle v$};
		
		\draw (208,109) node [anchor=north west][inner sep=0.5pt]   [align=left] {$\displaystyle e$};
		
		\draw (187,159) node [anchor=north west][inner sep=0.5pt]   [align=left] {$\displaystyle f$};
		
		\draw (72,104) node [anchor=north west][inner sep=0.5pt]   [align=left] {$\displaystyle e'$};
		
		\draw (99,117) node [anchor=north west][inner sep=0.5pt]   [align=left] {$\displaystyle u'$};
		
		\draw (89,165) node [anchor=north west][inner sep=0.5pt]   [align=left] {$\displaystyle f'$};
		
		\draw (87,198) node [anchor=north west][inner sep=0.5pt]   [align=left] {$\displaystyle v'$};
		
		\draw (142,118) node [anchor=north west][inner sep=0.5pt]   [align=left] {$\displaystyle g$};
		
		\end{tikzpicture}
	\end{center}
	As we have showed, the set $H=\{v\}$ is hereditary saturated, and for each irreducible polynomial $f(x)\in K[x, x^{-1}]$, the ideal $P:=I(H,f(e))$ is primitive of type III, $L_K(E)$ and $L_K(E)/P$ is non-commutative. If we set $p=e^\infty$, then $p$ is an infinite path with $E^0\backslash H=M(p)$. Also, it is clear that the edge $g$ satisfies the condition $u'=s(g)\ne r(g)=u=s(p)$. Therefore, if ${\rm char} K=t>0$ and $\alpha, \beta$ are algebraically independent over the prime subfield $\mathbb{F}_t$, then the following elements generate a non-cyclic free subgroup of $L_K(E)^\times$:
	$$
	c=u+u'+v+v'+\alpha g^* + (\beta-1) u' +(\beta^{-1} -1)u,
	$$ 
	$$
	d=u+u'+v+v'+\alpha g + (\beta-1) u' +(\beta^{-1} -1) u.
	$$ 
\end{example}

\section{Multiplicative group of some Leavitt path algebras}
In this final section, we shall determine the multiplicative groups of some kinds of Leavitt path algebras:

\bigskip

\noindent\textbf{One-sided noetherian Leavitt path algebra.} We begin with the multiplicative groups of one-sided noetherian Leavitt path algebras.
\begin{proposition}\label{proposition_7.1}
	Let $E$ be an (arbitrary) graph and $K$ a field. 
	\begin{enumerate}[font=\normalfont]
		\item[(1)] If $L_K(E)$ is left artinian, then  
		$$
		(L_K(E))^\times\cong  \prod_{i\in I }{\rm GL}_{n_i}(K),
		$$
		where $|I|$ is the number of sinks in $E$, and for each $i\in I$, $n_i$ is the number of paths in $E$ which end in the sink corresponding to $i$.
		\item[(2)] If $L_K(E)$ is left noetherian,
		$$
		(L_K(E))^\times\cong  \prod_{i\in I}{\rm GL}_{n_j}(K)\times  \prod_{j\in J}{\rm GL}_{m_j}(K[x,x^{-1}]),
		$$
		where $|I|$ is as above and $|J|$ the number of $($necessarily disjoint$)$ cycles in $E$, and for each $j\in J$, $m_j$ is the number of paths which end in the cycle corresponding $j$.
	\end{enumerate}
\end{proposition}
\begin{proof}
	For the proof of (1), we apply \cite[Corollary 4.2.13]{Bo_abrams_2017} to conclude that $E$ is a finite graph and that 
	$$
	L_K(E)\cong \bigoplus_{i\in I}{\rm M}_{n_i}(K),
	$$
	where $I$ is the number of sinks in $E$. It follows that $L_K(E)^\times\cong \prod_{i\in I}{\rm GL}_{n_i}(K).$ Similarly, the proof of (2) follows immediately from \cite[Corollary 4.2.14]{Bo_abrams_2017}.
\end{proof}

\begin{corollary}
	If $K$ is an absolute field and $L_K(E)$ is left artinian, then $\left( L_K(E)\right)^\times$ is a locally finite group.
\end{corollary}
\begin{proof}
	First we prove that for any $n\geq1$, the general linear group ${\rm GL}_n(K)$ is a locally finite group. Indeed, let $S$ be a finite subset of ${\rm GL}_n(K)$, and let $L$ be the subfield of $K$ generated by all entries of elements of $S$ over the prime subfield $\mathbb{F}_p$. Then, the subgroup $\left\langle S\right\rangle $ of $\mathrm{GL}_n(K)$ generated by $S$ can be considered as a subgroup of ${\rm GL}_n(L)$. Since $K$ is a locally finite field, we conclude that $L$ is finite, and so $\left\langle S\right\rangle $ is finite. Therefore ${\rm GL}_n(K)$ is a locally finite group, as claimed. 
	According to Proposition~ \ref{proposition_7.1}, we have
	$$
	\left(L_K(E)\right)^\times \cong \bigoplus_{i\in I}{\rm GL}_{n_i}(K).
	$$ 
	Since each ${\rm GL}_{n_i}(K)$ is a locally finite group, we conclude that $\left(L_K(E)\right)^*$ is locally finite too.
\end{proof}
\begin{remark}
	This corollary shows that in Theorem \ref{theorem_3.6} the condition that $K$ is a non-absolute field can not be omitted. 
\end{remark}
\noindent\textbf{Commutative Leavitt path algebras.} In \cite{Pa_pino_crow_2011}, Gonzalo Aranda Pino and Kathi Crow introduced the following notion. 
\begin{definition}
	A graph $E$ is said to be an \textit{almost disjoint union} of a collection of subgraphs $\{E\}_{i\in I}$ of $E$, for which we shall write $E=\coprod_{i\in I}E_i$, if the following conditions are satisfied:
	\begin{enumerate}[font=\normalfont]
		\item[(i)] $E^0=\bigcup_{i\in I}E_i^0$.
		\item[(ii)] $r_E^{-1}(E_i^0\cap E_j^0)=\varnothing$ for every $i\ne j$.
		\item[(iii)] $E_i^0\cap E_j^0\subseteq s_E(E_i^1)\cap s_E(E_j^1)$ for every $i\ne j$.
		\item[(iv)] $E_i^1$ is a partition of $E^1$.
	\end{enumerate}
\end{definition}
Of course if $E^0$ is finite, then so is all $E_i^0$. In this case, it is always true that $\sum_{i\in I}|E_i^0| \geq |E^0|$. However, the equality may not happen as we shall see in the following example (see \cite[Example 2.3]{Pa_pino_crow_2011}).
\begin{example} 
	Consider the graph $E$ and its two subgraphs $E_1$ and $E_2$:\\

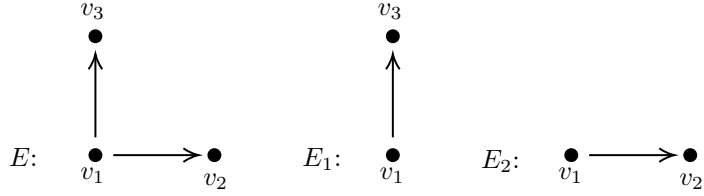
	\begin{center}
		
		\tikzset{every picture/.style={line width=0.75pt}} 
		
		\begin{tikzpicture}[x=0.75pt,y=0.75pt,yscale=-1,xscale=1]
		
		\draw  [color={rgb, 255:red, 0; green, 0; blue, 0 }  ,draw opacity=1 ][fill={rgb, 255:red, 0; green, 0; blue, 0 }  ,fill opacity=1 ] (57,135) .. controls (57,133.34) and (58.34,132) .. (60,132) .. controls (61.66,132) and (63,133.34) .. (63,135) .. controls (63,136.66) and (61.66,138) .. (60,138) .. controls (58.34,138) and (57,136.66) .. (57,135) -- cycle ;
		
		\draw  [color={rgb, 255:red, 0; green, 0; blue, 0 }  ,draw opacity=1 ][fill={rgb, 255:red, 0; green, 0; blue, 0 }  ,fill opacity=1 ] (117,135) .. controls (117,133.34) and (118.34,132) .. (120,132) .. controls (121.66,132) and (123,133.34) .. (123,135) .. controls (123,136.66) and (121.66,138) .. (120,138) .. controls (118.34,138) and (117,136.66) .. (117,135) -- cycle ;
		
		\draw    (69,135) -- (110,135) ;
		\draw [shift={(112,135)}, rotate = 180] [color={rgb, 255:red, 0; green, 0; blue, 0 }  ][line width=0.75]    (7.65,-3.43) .. controls (4.86,-1.61) and (2.31,-0.47) .. (0,0) .. controls (2.31,0.47) and (4.86,1.61) .. (7.65,3.43)   ;
		
		\draw  [color={rgb, 255:red, 0; green, 0; blue, 0 }  ,draw opacity=1 ][fill={rgb, 255:red, 0; green, 0; blue, 0 }  ,fill opacity=1 ] (57,75) .. controls (57,73.34) and (58.34,72) .. (60,72) .. controls (61.66,72) and (63,73.34) .. (63,75) .. controls (63,76.66) and (61.66,78) .. (60,78) .. controls (58.34,78) and (57,76.66) .. (57,75) -- cycle ;
		
		\draw    (60,126) -- (60,86) ;
		\draw [shift={(60,84)}, rotate = 450] [color={rgb, 255:red, 0; green, 0; blue, 0 }  ][line width=0.75]    (7.65,-3.43) .. controls (4.86,-1.61) and (2.31,-0.47) .. (0,0) .. controls (2.31,0.47) and (4.86,1.61) .. (7.65,3.43)   ;
		
		\draw  [color={rgb, 255:red, 0; green, 0; blue, 0 }  ,draw opacity=1 ][fill={rgb, 255:red, 0; green, 0; blue, 0 }  ,fill opacity=1 ] (207,135) .. controls (207,133.34) and (208.34,132) .. (210,132) .. controls (211.66,132) and (213,133.34) .. (213,135) .. controls (213,136.66) and (211.66,138) .. (210,138) .. controls (208.34,138) and (207,136.66) .. (207,135) -- cycle ;
		
		\draw  [color={rgb, 255:red, 0; green, 0; blue, 0 }  ,draw opacity=1 ][fill={rgb, 255:red, 0; green, 0; blue, 0 }  ,fill opacity=1 ] (207,75) .. controls (207,73.34) and (208.34,72) .. (210,72) .. controls (211.66,72) and (213,73.34) .. (213,75) .. controls (213,76.66) and (211.66,78) .. (210,78) .. controls (208.34,78) and (207,76.66) .. (207,75) -- cycle ;
		
		\draw    (210,126) -- (210,86) ;
		\draw [shift={(210,84)}, rotate = 450] [color={rgb, 255:red, 0; green, 0; blue, 0 }  ][line width=0.75]    (7.65,-3.43) .. controls (4.86,-1.61) and (2.31,-0.47) .. (0,0) .. controls (2.31,0.47) and (4.86,1.61) .. (7.65,3.43)   ;
		
		\draw  [color={rgb, 255:red, 0; green, 0; blue, 0 }  ,draw opacity=1 ][fill={rgb, 255:red, 0; green, 0; blue, 0 }  ,fill opacity=1 ] (297,135) .. controls (297,133.34) and (298.34,132) .. (300,132) .. controls (301.66,132) and (303,133.34) .. (303,135) .. controls (303,136.66) and (301.66,138) .. (300,138) .. controls (298.34,138) and (297,136.66) .. (297,135) -- cycle ;
		
		\draw  [color={rgb, 255:red, 0; green, 0; blue, 0 }  ,draw opacity=1 ][fill={rgb, 255:red, 0; green, 0; blue, 0 }  ,fill opacity=1 ] (357,135) .. controls (357,133.34) and (358.34,132) .. (360,132) .. controls (361.66,132) and (363,133.34) .. (363,135) .. controls (363,136.66) and (361.66,138) .. (360,138) .. controls (358.34,138) and (357,136.66) .. (357,135) -- cycle ;
		
		\draw    (309,135) -- (350,135) ;
		\draw [shift={(352,135)}, rotate = 180] [color={rgb, 255:red, 0; green, 0; blue, 0 }  ][line width=0.75]    (7.65,-3.43) .. controls (4.86,-1.61) and (2.31,-0.47) .. (0,0) .. controls (2.31,0.47) and (4.86,1.61) .. (7.65,3.43)   ;
		
		\draw (-27,165) node [anchor=north west][inner sep=0.75pt]   [align=left] {$ $};
		
		\draw (15,130) node [anchor=north west][inner sep=0.75pt]   [align=left] {$\displaystyle E$:};
		
		\draw (51,140) node [anchor=north west][inner sep=0.75pt]   [align=left] {$\displaystyle v_{1}$};
		
		\draw (113,143) node [anchor=north west][inner sep=0.75pt]   [align=left] {$\displaystyle v_{2}$};
		
		\draw (51,57) node [anchor=north west][inner sep=0.75pt]   [align=left] {$\displaystyle v_{3}$};
		
		\draw (163,130) node [anchor=north west][inner sep=0.75pt]   [align=left] {$\displaystyle E_{1}$:};
		
		\draw (202,140) node [anchor=north west][inner sep=0.75pt]   [align=left] {$\displaystyle v_{1}$};
		
		\draw (201,57) node [anchor=north west][inner sep=0.75pt]   [align=left] {$\displaystyle v_{3}$};
		
		\draw (253,131) node [anchor=north west][inner sep=0.75pt]   [align=left] {$\displaystyle E_{2}$:};
		
		\draw (292,140) node [anchor=north west][inner sep=0.75pt]   [align=left] {$\displaystyle v_{1}$};
		
		\draw (353,143) node [anchor=north west][inner sep=0.75pt]   [align=left] {$\displaystyle v_{2}$};
		\end{tikzpicture}
	\end{center}
	It is easy to see that $E=E_1\coprod E_2$ and that $|E_1^0|+|E_2^0|=4$ while $|E^0|=3$. In this following lemma, we give a certain situation in which the equality holds.
\end{example}
\begin{lemma}\label{lemma_7.4}
	Let $E$ be a graph. Assume that $E=\coprod_{i\in I}E_i$ such that $I$ is finite and for each $i$, the subgraph $E_i$ is either an isolated vertex or an isolated loop. Then $|I|=|E^0|$ and $E=\coprod_{i\in E^0}E_i$ is indeed a decomposition of $E$ into its connected components.
\end{lemma}
\begin{proof}
	It is trivial that each $E_i$ is connected. Moreover, the isolation of each $E_i$ assures us that this composition is disjoint. Also, as each $E_i$ determines and is determined by a unique vertex in $E_0$, we conclude that $|I|=|E^0|$.
\end{proof}
\begin{theorem}
	Let $E$ be a graph, and $K$ be an arbitrary field. Then, the following statements are equivalent:
	\begin{enumerate}[font=\normalfont]
		\item[(1)] $L_K(E)$ is a commutative.
		\item[(2)] $E$ is a graph such that $r$ and $s$ coincide and are injective.
		\item[(3)] $E=\coprod_{i\in I}E_i$, where $|I|$ is the number of vertices in $E$, and each subgraph $E_i$ is either an isolated vertex or the graph $R_1$.
		\item[(4)] $L_K(E)\cong \bigoplus_{i=0}^n K \oplus \bigoplus_{j=0}^m K[x,x^{-1}]$
		and that $m+n$ is equal to the number of vertices in $E$.
		\item[(5)] $(L_K(E))^\times\cong \prod _{i=0}^n K^\times \times \prod_{j=0}^m K^\times\left\langle x\right\rangle ,$
		where $m+n$ is the number of vertices in $E$, and $\left\langle x\right\rangle$ is the cyclic subgroup of $K[x,x^{-1}]$ generated by $x$.
		\item[(6)] $(L_K(E))^\times$ is an abelian group.
	\end{enumerate}
\end{theorem}
\begin{proof}
	It follows from \cite[Proposition 2.7]{Pa_pino_crow_2011} and Lemma \ref{lemma_7.4} that $(1)\Rightarrow(2)\Rightarrow(3)$. Since the set of invertible elements in $K[x,x^{-1}]$ is the product of $K^\times$ and the cyclic group $\left\langle x\right\rangle$, we have the implication $(3)\Rightarrow(4)$. The implication $(5)\Rightarrow(6)$ is clear. For the proof of $(6)\Rightarrow(1)$, we assume by contradiction that $L_K(E)$ is non-commutative. It follows from Theorem \ref{theorem_3.6} that $(L_K(E))^\times$ contains a non-cyclic free subgroup, which is a contradiction.
\end{proof}

\noindent\textbf{The Jacobson algebra}. Let $V$ be a $K$-vector space $Kv_1\oplus Kv_2\oplus \cdots$ with a (fixed) countably infinite basis $ \{v_i:i\geq 1\}$ over a field $K$ and ${\mathbb{N}}=\{1,2,3,\dots\}$ the set of natural numbers. With respect to this basis, the ring $E={\rm End}_K(V)$ of linear transformations of $V$ may be identified with the $K$-algebra of infinite  ${\mathbb{N}}\times{\mathbb{N}}$ matrices with finite number of nonzero entries in each column, and so the identity map $1_V\in E$ is identified with the infinite identity matrix ${\rm I}_\infty$. Consider the subalgebra $E_0$ of $E$ which consists of ${\mathbb{N}}\times{\mathbb{N}}$ matrices having finitely many nonzero entries in each row and each column. We denote by ${\rm M}_{\infty}(K)$ the set of finitary ${\mathbb{N}}\times{\mathbb{N}}$ matrices with only finitely many nonzero entries and by ${\rm GL}_{\infty}(K)$ the group of invertible matrices from ${\rm I}_{\infty}+{\rm M}_{\infty}(K)$. It is clear that ${\rm M}_{\infty}(K)$ is an ideal in $E_0$ and ${\rm GL}_{\infty}(K)$ is a subgroup of ${\rm Aut}(V)$.

For each natural number $n\geq 1$, the algebra of $n$-dimensional matrices ${\rm M}_n(K)$ can be written as ${\rm M}_n(K)=\bigoplus_{i,j=1}^nKf_{ij}$, where $f_{ij}$ are the matrix units. Accordingly, the following equality is obvious
$$
{\rm M}_{\infty}(K)=\lim\limits_{n\to\infty}{\rm M}_n(K)=\bigoplus_{i,j\in{\mathbb{N}}}Kf_{ij}.
$$

The \textit{Jacobson algebra}, denoted by $\mathcal{R}$, is defined to be the associative $K$-algebra on two non-commuting generators $x$ and $y$ modulo the relation $xy=1$ (\cite{Pa_abrams_2020}). In other words, 
$$
\mathcal{R}=K\left\langle xy\;|\;xy=1\right\rangle.
$$
It is clear that $\{y^ix^j|i,j\geq 0\}$ is a $K$-basis of $\mathcal{R}$; that is, $\mathcal{R}=\bigoplus_{i,j\geq 0}Ky^ix^j$. Such an algebra was first studied by Jacobson \cite{Pa_jacobson_1950} and investigated substantially by many other authors. As it was pointed out by Bavula in \cite[p.7]{Pa_bavula_2013}, the Jacobson algebra contains the ideal $F=\bigoplus_{i,j\in{\mathbb{N}}}Ke_{ij}$, where 
\begin{align}
	e_{ij}:=y^{i-1}x^{j-1}-y^ix^j, \;\; i,j\geq1.
\end{align}
The ideal $F$ is an algebra (without identity $1$) which is isomorphic to ${\rm M}_{\infty}(K)$ via $e_{ij}\to f_{ij}$. Furthermore, for all $i,j\geq 1$, we have 
\begin{align}
	xf_{ij}=f_{i-1,j}, \;\; y f_{ij}=f_{i+1,j} \;\; (f_{0,j}:=0).
\end{align}
\begin{align}
	\mathcal{R}=K\oplus xK[x]\oplus yK[y]\oplus F.
\end{align}

Now let us consider the Toeplitz $K$-algebra $\mathscr{T}_K=L_K(E_T)$; that is, the Leavitt path algebra of the Toeplitz graph $E_T$ (see the picture in the end of Introduction) over $K$. This algebra is unital with the identiy $1_{\mathscr{T}_K}=u+v$. In $\mathscr{T}_K$, if we set $X=e^*+f^*$ and $Y=e+f$, then $XY=1$ and $YX=u\ne 1$. It is straightforward to check that the map $\varphi:\mathcal{R}\to \mathscr{T}_K$ given by $\varphi(x)=X$ and $\varphi(y)=Y$ is a $K$-algebra isomorphism. This implies that $\{Y^iX^j\;|\; i,j\geq0\}$ is a $K$-basis of $\mathscr{T}_K$, and so $\mathscr{T}_K$ is isomorphic to a subalgebra of $E_0$ which consists of infinite matrices having finitely many nonzero entries in each row and each column. According to (1), $ \mathscr{T}_K$ contains the ideal $I=\bigoplus_{i,j\in{\mathbb{N}}}KF_{ij}$, where 
$$
F_{ij}:=Y^{i-1}X^{j-1}-Y^iX^j=(e+f)^{i-1}v(e^*+f^*)^{j-1}, \;\; i,j\geq1.
$$
It follows that $I$ isomorphic to ${\rm M}_{\infty}(K)$ via $F_{ij}\to f_{ij}$. From (3), we easily make the following computations: $$F_{11}=v,\;\;\;F_{21}=(e+f)v=f,\;\;\; F_{12}=v(e
^*+f^*)=f^*.$$
Therefore, we may identify $v$, $f$, $f^*$ with the matrices $f_{11}$, $f_{21}$, $f_{12}$ in ${\rm M}_{\infty}(K)$, respectively. Consequently, the vertex $u$ is identified with the matrix $\sum_{i\geq2}f_{ii}\in E_0$. Next, we shall find the corresponding matrices, say $A$ and $B$ in $E_0$, of $e$ and $e^*$.
From (2), we have 
$$
(e^*+f^*)F_{ij}=F_{i-1,j}, \;\; (e+f) F_{ij}=F_{i+1,j} \;\; (F_{0,j}:=0).
$$
This implies that 
$$
(B+f_{12})f_{ij}=f_{i-1,j}, \;\; (A+f_{21}) f_{ij}=f_{i+1,j} \;\; (f_{0,j}:=0).
$$
$$
B+f_{12}=\sum_{i\geq1}f_{i,i+1}, \;\;\;A+f_{21}=\sum_{i\geq1}f_{i+1,i}.
$$
It follows that $A=\sum_{i\geq2}f_{i+1,i}$ and $B=\sum_{i\geq2}f_{i,i+1}$. From (3), we also have 
$$
\mathscr{T}_K=K\oplus (e^*+f^*)K[e^*+f^*]\oplus (e+f)K[e+f]\oplus I.
$$
From this, it is easy to see that $\mathscr{T}_K$ is isomorphic to the subalgebra of $E_0$ generated by $A$, $B$ and ${\rm M}_{\infty}(K)$. Now, we summary what we have discussed in the following lemma.

\begin{lemma}\label{lemma_7.6}
	Let $\mathscr{T}_K$ be the Toeplitz $K$-algebra and $E_0$ the $K$-algebra of infinite matrices over $K$ with finite number of nonzero entries in each row and each column, and ${\rm M}_{\infty}(K)$ the set of finitary matrices. Then 
	$$
	\mathscr{T}_K\cong \left\langle \sum_{i\geq2}f_{i+1,i}, \sum_{i\geq2}f_{i,i+1}, {\rm M}_{\infty}(K)\right\rangle \subseteq E_0.
	$$
\end{lemma}

\begin{theorem}\label{theorem_7.7}
	Let ${\rm GL}_{\infty}(K)$ be the set of invertible matrices from ${\rm I}_{\infty}+{\rm M}_{\infty}(K)$. Then $\left( \mathscr{T}_K\right) ^\times \cong K^\times\rtimes\; {\rm GL}_{\infty}(K)$, the semidirect product of $K^\times$ and ${\rm GL}_{\infty}(K)$. In other words, the exact sequence $1\to F^\times\to \left( \mathscr{T}_K\right) ^\times \to {\rm GL}_{\infty}(K)\to 1$ is split.
\end{theorem}
\begin{proof}
	It follows from Lemma \ref{lemma_7.6} that $\mathscr{T}_K$ can be identified with $\left\langle E, E^*, {\rm M}_{\infty}(K)\right\rangle$. This means that $(\mathscr{T}_K)^\times$ contains a subgroup which is isomorphic to $K^\times\rtimes\; {\rm GL}_{\infty}(K)$. At the other extreme, for any $r\in (\mathscr{T}_K)^\times$, the left multiplication $L_r$ by $r$ is an element of ${\rm Aut}(\mathscr{T}_K)$. It is straightforward to check that the assignment  $r\to L_r$ defines an injective group homomorphism from  $(\mathscr{T}_K)^\times$ to ${\rm Aut}(\mathscr{T}_K)$. Therefore $(\mathscr{T}_K)^\times$ may be embedded in $K^\times\rtimes\; {\rm GL}_{\infty}(K)$.  It follows from \cite{Pa_alah_2013} that ${\rm Aut}(\mathscr{T}_K)\cong K^\times\rtimes\; {\rm GL}_{\infty}(K)$. This means that $(\mathscr{T}_K)^\times\cong K^\times\rtimes\; {\rm GL}_{\infty}(K)$.
\end{proof}
\begin{remark}\label{remark_2}
	At this point, it is interesting to see that there are two different ways to obtain non-cyclic free subgroups of $(\mathscr{T}_K)^\times$. The first one is obtained by applying Corollary~ \ref{corollary_5.4}. Indeed, if we assume that ${\rm char} K=0$ and $\alpha\in K$ is equal to $2$ or transcendental over the prime subfield ${\mathbb{Q}}$, then Corollary \ref{corollary_5.4}(i) implies that $G=\left\langle u+v+\alpha f^*, u+v+\alpha f \right\rangle $ is a non-cyclic free subgroup of $(\mathscr{T}_K)^\times$. The second way is due to Theorem \ref{theorem_7.7}. More precisely, from this theorem we know that $\left( \mathscr{T}_K\right) ^\times \cong K^\times\rtimes\; {\rm GL}_{\infty}(K)$ which contains the non-cyclic free subgroup $H=\left\langle {\rm I}_\infty+\alpha f_{21}, {\rm I}_\infty+\alpha f_{12}\right\rangle $. (Note that $H$ is a copy of the subgroup $ \left\langle A,B\right\rangle $ of ${\rm GL}_2(K)$ in Lemma \ref{lemma_3.5}). As we have pointed out before, the inverse images of $f_{21}$ and $f_{12}$ are $f$ and $f^*$ respectively, so we can conclude that $H\cong G$.
\end{remark}
\begin{corollary}
	If $K$ is an absolute field, then $(\mathscr{T}_K)^\times$ is a locally finite group.
\end{corollary}
\begin{proof}
	Let $S$ be a finite subset of ${\rm GL}_{\infty}(K)$. Let $L$ be the subfield of $K$ generated by all entries of all matrices in $S$. Since $K$ is absolute, the subfield $L$ is finite. It follows that ${\rm M}_{\infty}(L)$ is finite, so ${\rm GL}_{\infty}(L)\subseteq {\rm I}_{\infty}+{\rm M}_{\infty}(L)$ is finite too. Because the subgroup of ${\rm GL}_{\infty}(K)$ generated by $S$ is contained in ${\rm GL}_{\infty}(L)$, we conclude that this subgroup is finite. This means that ${\rm GL}_{\infty}(K)$ is a locally finite group. It is also easy to see that $K^\times$ is locally finite group. Finally, we conclude that $(\mathscr{T}_K)^\times \cong K^\times\rtimes\; {\rm GL}_{\infty}(K)$ is locally finite.
\end{proof}
\begin{remark}
	For a different approach to the multiplicative group of Jacobson algebra, see \cite[Theorem 4.6]{Pa_bavula_2013}.
\end{remark}

In the remaining part of this paper, we give a positive answer to Conjecture~ \ref{conjecture_3.1} for the case of the Toeplitz $K$-algebra $\mathscr{T}_K$. For this purpose, we need the concept of \textit{global determinant} of an infinite matrix which is introduced in \cite[p.15]{Pa_bavula_2013} by Bavula. For each $n\geq 1$, consider the usual determinant $\det_n=\det : {\rm I}_n+{\rm M}_n(K)\to K,u\mapsto \det(u)$. These determinants determine the global determinant
$$
\det: {\rm I}_{\infty}+{\rm M}_{\infty}(K)\to K, u\mapsto\det(u),
$$
where $\det(u)$ is the common value of all determinants of $\det_n(u)$, $n\gg1$. The global determinant has the usual properties of the determinant, and we also have 
$${\rm GL}_{\infty}(K)=\{u\in {\rm I}_\infty+{\rm M}_\infty(K)\;|\;\det(u)\ne 0\}.
$$
The kernel 
$$
{\rm SL}_{\infty}(K)=\{u\in {\rm I}_\infty+{\rm GL}_{\infty}(K)\;|\;\det(u)= 1\}.
$$
is a normal subgroup of ${\rm GL}_{\infty}(K)$. From this, we have the following:
\begin{proposition}\label{proposition_7.9}
	Let $K$ be a field which is different from $\mathbb{F}_2,  \mathbb{F}_3$. Then a non-central subgroup $N$ of ${\rm GL}_{\infty}(K)$ is normal if and only if ${\rm SL}_{\infty}(K)\subseteq N$.
\end{proposition}
\begin{proof}
	Assume that $N$ is a non-central normal subgroup of ${\rm SL}_{\infty}(K)$. Take an arbitrary matrix $u\in {\rm SL}_{\infty}(K)$. Then, there exists an integer $n\geq2$ so that $u$ can be considered as an element of $ {\rm SL}_{n}(K)\subseteq {\rm GL}_n(K)$.  Because $N$ is non-central in ${\rm GL}_{\infty}(K)$, the integer $n$ may be chosen to be large enough such that  $N\cap{\rm GL}_n(K)$ is a non-central normal subgroup of ${\rm GL}_n(K)$. We know that if $K\ne \mathbb{F}_2,  \mathbb{F}_3$ and $n\geq 2$, then a subgroup of ${\rm GL}_n(K)$ is normal if and only if it contains ${\rm SL}_n(K)$, so we have $u\in{\rm SL}_{n}(K) \subseteq N\cap{\rm GL}_n(K) \subseteq N$. Therefore ${\rm SL}_{\infty}(K)\subseteq N$. The converse is obvious. 
\end{proof}
\begin{corollary}
	Let $K$ be a non-absolute field. Then, every non-central normal subgroup of $(\mathscr{T}_K)^\times$ contains a non-cyclic free subgroup.
\end{corollary}
\begin{proof}
	Let $N$ be a non-central normal subgroup of $(\mathscr{T}_K)^\times$. From the proof of Theorem~ \ref{theorem_7.7}, we see that $\left( \mathscr{T}_K\right) ^\times \cong K^\times\rtimes\; {\rm GL}_{\infty}(K)$. By Proposition \ref{proposition_7.9}, we conclude that $N$ contains a copy of ${\rm SL}_{\infty}(K)$. It is simple to see that $N$ contained the group $H$ given in Remark \ref{remark_2} and the result follows.
\end{proof}


\begin{thebibliography}{}
	
	\bibitem{Pa_alah_2013} A. Alahmedi, H. Alsulami, S. Jain, E. Zelmanov, Structure of Leavitt Path Algebras of Polynomial Growth, Proc. Natl. Acad. Sci. USA 110 (2013) 15222--15224.
		
	\bibitem{Pa_abrams-pino_2005} G. Abrams, G. Aranda Pino, The Leavitt path algebra of a graph, J. Algebra 293 (2005) 319--334.
	
	\bibitem{Bo_abrams_2017} G. Abrams, P. Ara, and M. Siles Molina, Leavitt path algebras, Lecture Notes in Mathematics series, Vol. 2191, Springer-Verlag Inc., 2017.
		
	\bibitem{Pa_abrams_2020} G. Abrams, F. Mantese,  A. Tonolo, Injective modules over the Jacobson algebra $K\langle X, Y | XY=1\rangle $, Canad. Math. Bull. (2020) 1--17.
	
	\bibitem{Pa_ara_rangas_2014} P. Ara, K.M. Rangaswamy, Finitely presented simple modules over Leavitt path algebras, J. Algebra 417 (2014) 333--352.
	
	\bibitem{Pa_bavula_2013} V. V. Bavula, The group of automorphisms of the algebra of one-sided inverses of a polynomial algebra, M\"{u}nster J. of Math. 6 (2013) 1--51.
	
	\bibitem{Pa_chen_2012} X. W. Chen, Irreducible representations of Leavitt path algebras, Forum Math. 27 (2015) 549--574.
	
	\bibitem{Pa_jacobson_1950} N. Jacobson, Some Remarks on One-Sided Inverses, Proc. A.M.S. 1 (1950) 352--355.
	
	\bibitem{Bo_lam_2001} T. Y. Lam, A first course in noncommutative rings, volume 131 of Graduate Texts in Mathematics. Springer-Verlag, New York, second edition, 2001.
	
	\bibitem{Pa_lanski_1971} C. Lanski, Subgroups and conjugates in semiprime rings, Math. Ann. 192 (1971) 313--327.
	
	\bibitem{Pa_lanski_1981} C. Lanski, Solvable subgroups in prime rings, Proc. Amer. Math. Soc. 82 (1981) 533-537.
	
	\bibitem{Pa_pino_crow_2011} G. Aranda Pino, K. Crow, The Center of a Leavitt patha lgebra, Rev. Math. Iberoam 27(2) (2011) 621--644.
	
	\bibitem{Pa_rangas_2015} K.M. Rangaswamy, On simple modules over Leavitt path algebras, J. Algebra 431 (2015) 239--258.
	
	\bibitem{Pa_rangas_2013} K. M. Rangaswamy. The theory of prime ideals of Leavitt path algebras over arbitrary graphs, J. Algebra, 375 (2013) 73--96.
	
	\bibitem{Bo_Rotman_2015} J. J. Rotman, Advanced modern algebra, I. 3rd ed., Grad. Stud. Math., Vol. 15, Amer. Math. Soc., Providence, RI, 2015.
	
	
	\bibitem{Bo_weh_1973} B. A. F. Wehrfritz, Infinite Linear Groups (1973), Berlin: Springer-Verlag, Berlin.
	
	\bibitem{Pa_tom_2007} M. Tomforde, Uniqueness theorems and ideal structure of Leavitt path algebras, J. Algebra 318 (2007) 270--299.
		
\end{thebibliography}
\end{document}